%% file: main.tex
\documentclass[11pt]{article}
\usepackage{fullpage}
\usepackage{authblk}
\usepackage{natbib}
\usepackage{algorithm}
\usepackage{algpseudocode}
\usepackage{amsmath,amsthm,amsfonts}
\usepackage{amsmath}
\usepackage[subnum]{cases}
\usepackage{tcolorbox}
\usepackage{mathrsfs}
\usepackage{amssymb}
\usepackage{color}
\usepackage{mathrsfs}
\usepackage{enumitem}
\usepackage{bm}
\usepackage{multirow}
\usepackage{booktabs}
\usepackage{makecell}
\usepackage{graphicx}
\usepackage{subcaption}
\usepackage{comment}
\usepackage{cases}
\usepackage{appendix}
\usepackage{tikz}
\usetikzlibrary{arrows,shapes}
\usepackage[colorlinks= true, linkcolor=brown, citecolor=blue, urlcolor=blue]{hyperref}
\usepackage{cleveref}
\usepackage{marginnote}
\usepackage{diagbox} 
\usepackage{mdframed}
\usepackage{wrapfig}
\usepackage{textcomp}
\usepackage{empheq}

\numberwithin{equation}{section}

\theoremstyle{definition}

\newcommand\scalemath[2]{\scalebox{#1}{\mbox{\ensuremath{\displaystyle #2}}}}

\newcommand{\norm}[1]{\left\|#1\right\|}

\DeclareMathOperator*{\argmin}{arg\,min}

\usepackage{booktabs}
\usepackage{tikz}
\usepackage{pgfplots}
\usetikzlibrary{calc,intersections}

\usetikzlibrary{shapes.geometric, arrows, positioning}

\tikzset{
	mybox/.style  = {draw, rectangle, minimum width=4cm, minimum height=0.8cm, text centered, text width=4.4cm,   
		font=\normalsize},
	box/.style  = {draw, rectangle, minimum width=2.0cm, minimum height=0.6cm, text centered, text width=3.0cm,   
		font=\normalsize},
	myarrow/.style = {line width=0.2pt, draw=black, -triangle 60, postaction={draw, line width=0.2pt, shorten >=10pt,-}}
}

\tikzstyle{arrow} = [->, >=stealth, -triangle 60]

\makeatletter
\newcommand{\leqnomode}{\tagsleft@true}
\newcommand{\reqnomode}{\tagsleft@false}
\makeatother
\allowdisplaybreaks

\begin{document}
%\maketitle

\title{Numerical Solution for Nonlinear 4D Variational Data Assimilation (4D-Var) via ADMM}
%\title{An adjoint-free algorithm for Nonlinear 4D Variational Data Assimilation via ADMM and Sampling}

\author[1,2]{Bowen Li}
\author[1,2]{Bin Shi\thanks{Corresponding author: \url{shibin@lsec.cc.ac.cn} } }
\affil[1]{State Key Laboratory of Scientific and Engineering Computing, Academy of Mathematics and Systems Science, Chinese Academy of Sciences, Beijing 100190, China}
\affil[2]{School of Mathematical Sciences, University of Chinese Academy of Sciences, Beijing 100049, China}
\date\today

\maketitle

%\bin{Vanishing Viscosity should be included in the main part, the view is key important. At least, it is mentioned a little.}\wjs{OK.}
\begin{abstract}
The four-dimensional variational data assimilation (\texttt{4D-Var}) has emerged as an important methodology, widely used in numerical weather prediction, oceanographic modeling, and climate forecasting. Classical unconstrained gradient-based algorithms often struggle with local minima, making their numerical performance highly sensitive to the initial guess. In this study, we exploit the separable structure of the \texttt{4D-Var} problem to propose a practical variant of the alternating direction method of multipliers (\texttt{ADMM}), referred to as the linearized multi-block \texttt{ADMM} with regularization. Unlike classical first-order optimization methods that primarily focus on initial conditions, our approach derives the Euler-Lagrange equation for the entire dynamical system, enabling more comprehensive and effective utilization of observational data. When the initial condition is poorly chosen, the $\argmin$ operation steers the iteration towards the observational data, thereby reducing sensitivity to the initial guess. The quadratic subproblems further simplify the solution process, while the parallel structure enhances computational efficiency, especially when utilizing modern hardware. To validate our approach, we demonstrate its superior performance using the Lorenz system, even in the presence of noisy observational data. Furthermore, we showcase the effectiveness of the linearized multi-block \texttt{ADMM} with regularization in solving the~\texttt{4D-Var} problems for the viscous Burgers' equation, across various numerical schemes, including finite difference, finite element, and spectral methods. Finally, we illustrate the recovery of dynamics under noisy observational data in a 2D turbulence scenario, particularly focusing on vorticity concentration, highlighting the robustness of our algorithm in handling complex physical phenomena.
\end{abstract}
%%
%\thispagestyle{empty}
%\setcounter{page}{0}
%%
%%%\newpage
%\tableofcontents
%%%\newpage

\input{01_intro}

\input{02_admm}

\input{03_burgers}

\input{04_vorticity}
\input{05_conclu}
\section*{Acknowledgements}
Bowen Li was partially supported by the Hua Loo-Keng scholarship of CAS.  Bin Shi was partially supported by Grant No.12241105 of NSFC.

{\small
\bibliographystyle{abbrvnat}
\bibliography{reference}
}
\end{document}

%% file: 01_intro.tex
\section{Introduction}
\label{sec: intro}

Data assimilation is a methodology used to estimate the evolving state of a dynamical system by integrating observational data with the system's underlying dynamics. By leveraging both sources, it produces more accurate state estimates than relying solely on either observations or model predictions~\citep{asch2016data}. Unlike approaches such as machine learning, image analysis, or statistical methods, data assimilation is fundamentally anchored in the dynamics of the system under study. As a result, it has emerged as a prominent research area in mathematics~\citep{law2015data, alonso2023inverse}. A key application of data assimilation lies in determining initial conditions, where distributed observations collected over time are combined with the dynamic model, which significantly enhances forecast accuracy~\citep{kalnay2002atmospheric}. This methodology has found widespread applications across various scientific fields, including meteorology~\citep{bouttier2002data}, oceanography~\citep{munk1982observing, wunsch1996ocean}, and climatology~\citep{siedler2013ocean, wunsch2013dynamically, stammer2016ocean}.

The application of the calculus of variations to data assimilation was first pioneered by~\citet{sasaki1958objective} for meteorological analysis, with subsequent extensions that incorporated the time dimension~\citep{sasaki1969proposed, sasaki1970some, thompson1969reduction}. These advancements laid the groundwork for what is now known as~\textit{4D Variational Data Assimilation}, commonly referred to as \texttt{4D-Var}. The introduction of the adjoint method by~\citet{lions1971optimal} and~\citet{marchuk1975formulation, marchuk1975formulation2} provided an efficient framework for obtaining gradient information, enabling the use of classical first-order optimization algorithms to solve these types of problems. The adoption of~\texttt{4D-Var} in meteorological data assimilation gained significant momentum during the 1980s, with its theoretical foundations established by~\citet{le1986variational} and~\citet{talagrand1987variational}.

Consider a finite-dimensional dynamical system described by:
\begin{equation}
\label{eqn: finite-ds}
\left\{ \begin{aligned}
          & \pmb{u}_t = G(\pmb{u})                                  \\
          & \pmb{u} |_{t=0} = \pmb{u}_0,                     
           \end{aligned} \right.
\end{equation}
which $\pmb{u} = (u_1, u_2, \ldots, u_m)^{\top} \in \mathbb{R}^m$ is an $m$-dimensional vector, and its Euclidean norm is defined as:
\begin{equation}
\label{eqn: euclidean-norm}
\|\pmb{u}\| = \left( \sum_{i=1}^{m} |u_i|^2 \right)^{\frac12}. 
\end{equation}
Let $\hat{\pmb{u}}$ represent the (partial) observational data collected over time, with indices chosen from a subset $\mathcal{S} \subseteq \{1, 2, \ldots, m\}$.  Additionally, let $\hat{\pmb{u}}_0^b$ denote the data obtained from a prior prediction. The~\texttt{4D-Var} problem can be rigorously formulated in the following variational form:
\begin{equation}
\label{eqn: analytic-4dvar}
\left\{ \begin{aligned}
          & \min F(\pmb{u}) := \frac{1}{2} \int_{0}^{T} \|\pmb{u} - \hat{\pmb{u}}\|^2dt + \frac{\alpha}{2} \big\| \pmb{u}_0 - \hat{\pmb{u}}_0^b \big\|^2 \\
          & \text{s.t.} \quad \pmb{u}\;~\text{satisfies}\;~\eqref{eqn: finite-ds},
          \end{aligned} \right.
\end{equation}
where $\alpha > 0$ is a regularization coefficient. In practice, observational data are typically collected at discrete time points. For simplicity, assume that observation times are uniformly spaced with  $T_{o} = T/n$, meaning that the data is recorded at times $t = kT_o$ for $k=0,1,\ldots,n$. Let $\pmb{u}_k = \pmb{u}(kT_o)$ represent the solution at time $t=kT_0$. The dynamical system~\eqref{eqn: finite-ds} can then be expressed at these observation times as: 
\begin{equation}
\label{eqn: finite-ds-discrete}
\pmb{u}_{k+1} = H_{T_{o}}(\pmb{u}_{k}), \quad k = 0, 1, \ldots, n -1,
\end{equation}
where $H_{T_{o}}$ is an operator acting over the observational time interval $T_0$, which may be linear or nonlinear, and $\pmb{u}_0$ is the given initial condition. Consequently, the continuous~\texttt{4D-Var} problem~\eqref{eqn: analytic-4dvar} is transformed into the following discrete form:
\begin{equation}
\label{eqn: numerical-4dvar}
\left\{ \begin{aligned}
           & \min F(\pmb{u}_0, \pmb{u}_1, \ldots, \pmb{u}_n) := \frac{T_o}{2} \sum_{k=0}^{n}\|\pmb{u}_k - \hat{\pmb{u}}_k\|^2 + \frac{\alpha}{2} \big\| \pmb{u}_0 - \hat{\pmb{u}}_0^b \big\|^2 \\
           & \text{s.t.}\quad \pmb{u}_{k+1} = H_{T_o}(\pmb{u}_{k}),\quad \text{for}\;\; k = 0, 1, \ldots, n -1.
            \end{aligned} \right.
\end{equation}

%To solve the~\texttt{4D-Var} problem~\eqref{eqn: analytic-4dvar} numerically, we discretize the time interval by setting $T_{o} = T/n$. Using numerical methods such as the forward Euler or ,

%%%%%%%%%%%%%%%%%%%%%%%%%%%%%%%%%%%%%%%%%%%%%%%%%%%%%%%%%%%%%%%%%%%%%%%%%%%%%%%%%%%%%%%%%%%%%%%%%%%%%%%%%%%%%%%%%%%
\subsection{Linear dynamical system}
\label{subsec: lin-dyn}

When the dynamical system~\eqref{eqn: finite-ds} is linear, the update in the dynamical iteration~\eqref{eqn: finite-ds-discrete} can be written as $H_{T_o}(\pmb{u}_{k}) = L\pmb{u}_{k}$, where $L$ is a linear operator or matrix. This simplifies the discrete dynamical system~\eqref{eqn: finite-ds-discrete} to:
\begin{equation}
\label{eqn: lin-discrete-ds}
\pmb{u}_{k} = L^k \pmb{u}_0, \quad \text{for}\;\; k = 0, 1, \ldots, n,
\end{equation}
where $\pmb{u}_0$ is the given initial condition. In this case, the~\texttt{4D-Var} Problem~\eqref{eqn: numerical-4dvar} reduces to a function that depends solely on the initial condition $\pmb{u}_0$, rather than on all the intermediate states $\pmb{u}_0, \pmb{u}_1, \ldots, \pmb{u}_n$. The linear \texttt{4D-Var} Problem can thus be reformulated as:
\begin{equation}
\label{eqn: numerical-4dvar-linear}
\min_{\pmb{u}_0 \in \mathbb{R}^m} F(\pmb{u}_0) = \frac{T_o}{2} \sum_{k=0}^{n}\big \|L^k \pmb{u}_0 - \hat{\pmb{u}}_k \big\|^2 + \frac{\alpha}{2} \big\| \pmb{u}_0 - \hat{\pmb{u}}_0^b \big\|^2.
\end{equation}
This formulation shows that the linear~\texttt{4D-Var} Problem~\eqref{eqn: numerical-4dvar} is essentially a quadratic function of the initial condition $\pmb{u}_0$. To solve the linear~\texttt{4D-Var} problem~\eqref{eqn: numerical-4dvar-linear}, classical first-order optimization algorithms, such as the conjugate gradient method~\citep{hestenes1952methods} and the limited-memory BFGS-B method~\citep{byrd1995limited, zhu1997algorithm}, can be applied. The gradient required for these algorithms is computed as follows:
\begin{equation}
\label{eqn: gradient-4Dvar-linear}
\nabla F(\pmb{u}_0) = T_o \sum_{k=0}^{n}(L^{k})^{\top}\left(L^k\pmb{u}_0 - \hat{\pmb{u}}_k \right) + \alpha \left( \pmb{u}_0 - \hat{\pmb{u}}_0^b \right), 
\end{equation}
where it can be observed that a key aspect of this gradient computation is the adjoint operator, $(L^k)^{\top}$. This procedure, commonly known as the adjoint method, was first introduced in~\citet{le1986variational} and~\citet{talagrand1987variational}. When the system's dimension is too large to be stored $L^{k}$ explicitly, iterative methods are typically employed to compute the adjoint operator $(L^{k})^{\top}$, which poses significant challenges, particularly in extending the approach to nonlinear systems.

%%%%%%%%%%%%%%%%%%%%%%%%%%%%%%%%%%%%%%%%%%%%%%%%%%%%%%%%%%%%%%%%%%%%%%%%%%%%%%%%%%%%%%%%%%%%%%%%%%%%%%%%%%%%%%%%%%%
\subsection{The Lorenz system --- a nonlinear dynamical system}
\label{subsec: non-lin-dyn}

For the nonlinear dynamical system~\eqref{eqn: finite-ds}, where the operator $H_{T_o}$ in the dynamical iteration~\eqref{eqn: finite-ds-discrete} is nonlinear, the~\texttt{4D-Var} Problem~\eqref{eqn: numerical-4dvar} can still be reformulated in terms of the initial condition $\pmb{u}_0$ as follows: 
\begin{equation}
\label{eqn: numerical-4dvar-nonlinear}
\min_{\pmb{u}_0 \in \mathbb{R}^m} F(\pmb{u}_0) = \frac{T_o}{2} \sum_{k=0}^{n}\big \|H_{T_o}^k(\pmb{u}_0) - \hat{\pmb{u}}_k \big\|^2 + \frac{\alpha}{2} \big\| \pmb{u}_0 - \hat{\pmb{u}}_0^b \big\|^2.
\end{equation}
We can generalize the procedure from the linear case to solve the nonlinear problem.  To apply first-order optimization algorithms to solve the~\texttt{4D-Var} problem~\eqref{eqn: numerical-4dvar-nonlinear}, the gradient must be computed as follows:
\begin{equation}
\label{eqn: gradient-4Dvar-nonlinear}
\nabla F(\pmb{u}_0) = T_o \sum_{k=0}^{n}\left(\nabla H_{T_o}^{k}(\pmb{u}_0) \right)^{\top}\left(H_{T_o}^k(\pmb{u}_0) - \hat{\pmb{u}}_k \right) + \alpha \left( \pmb{u}_0 - \hat{\pmb{u}}_0^b \right).
\end{equation}
The key difference between the linear case~\eqref{eqn: gradient-4Dvar-linear} and the nonlinear case~\eqref{eqn: gradient-4Dvar-nonlinear} is that in the nonlinear case, the gradient involves the Jacobian of the nonlinear operator $H_{T_o}^{k}$, denoted as $\nabla H_{T_o}^{k}(\pmb{u}_0)$.  To compute this Jacobian efficiently, the adjoint method~\citep{le1986variational,talagrand1987variational} is extended to nonlinear systems using the tangent linear system, which  for the dynamical system~\eqref{eqn: lin-discrete-ds} is given by: 
\begin{equation}
\label{eqn: tlm}
\delta \pmb{u}_t = \nabla G(\pmb{u}) \cdot \delta \pmb{u}.
\end{equation}
From this tangent linear system~\eqref{eqn: tlm}, the Jacobian can be derived as:
\begin{equation}
\label{eqn: jacobian}
\nabla H_{T_o}^{k}(\pmb{u}_0) = \exp\left( \int_{0}^{kT_{o}} \nabla G(\pmb{u}(\pmb{u}_0, t)) dt \right).
\end{equation}
%
%Based on the same idea, the process of obtaining the gradient can be generalized from the linear case to nonlinear case
%
%The gradient of the nonlinear~\texttt{4D-Var} problem~\eqref{eqn: numerical-4dvar-nonlinear} can be computed as:
%
%
%However, comparison between

\begin{figure}[htb!]
\centering
\begin{subfigure}[t]{0.45\linewidth}
\centering
\includegraphics[scale=0.22]{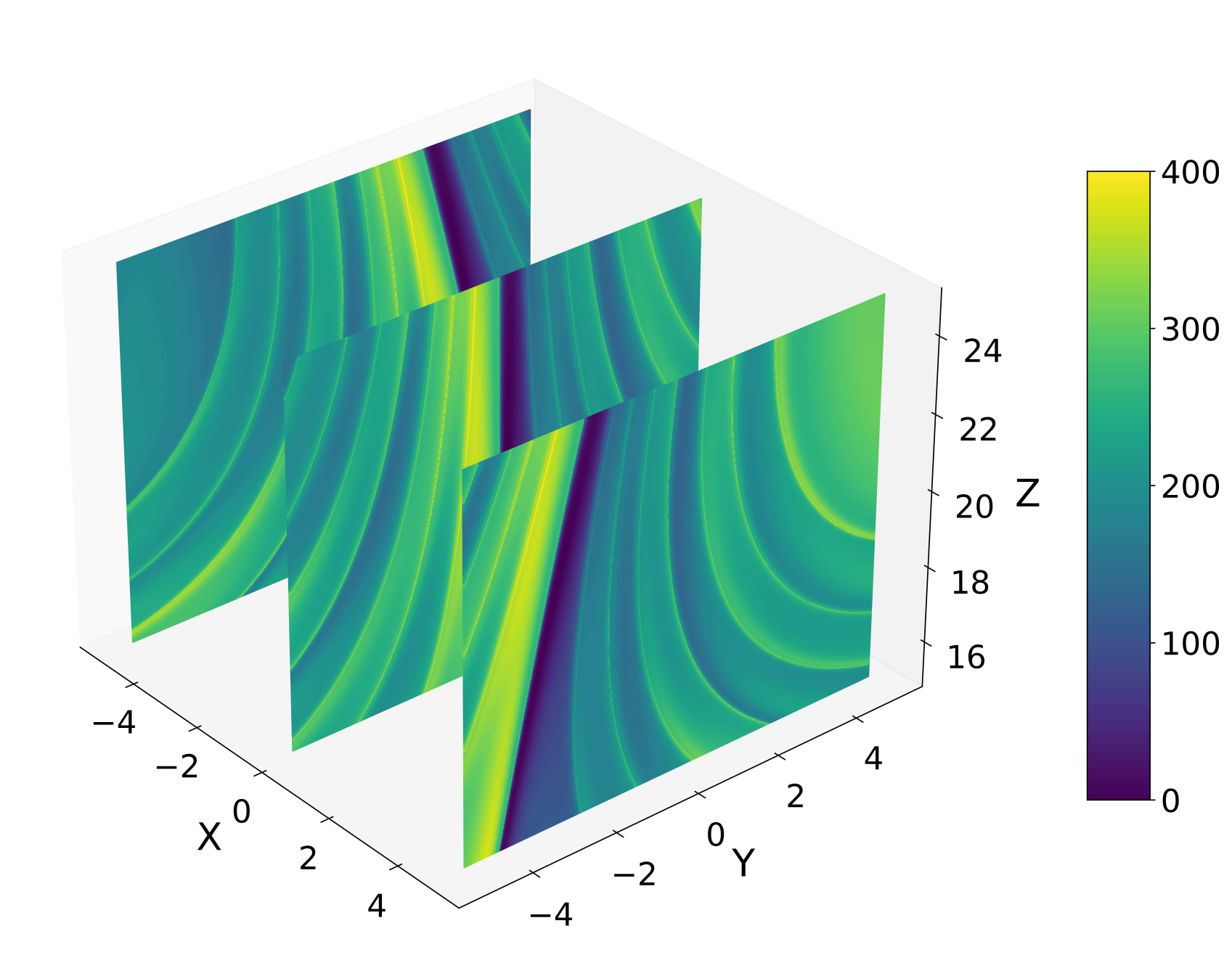}
\caption{$X$-slices}
\end{subfigure}
\begin{subfigure}[t]{0.45\linewidth}
\centering
\includegraphics[scale=0.22]{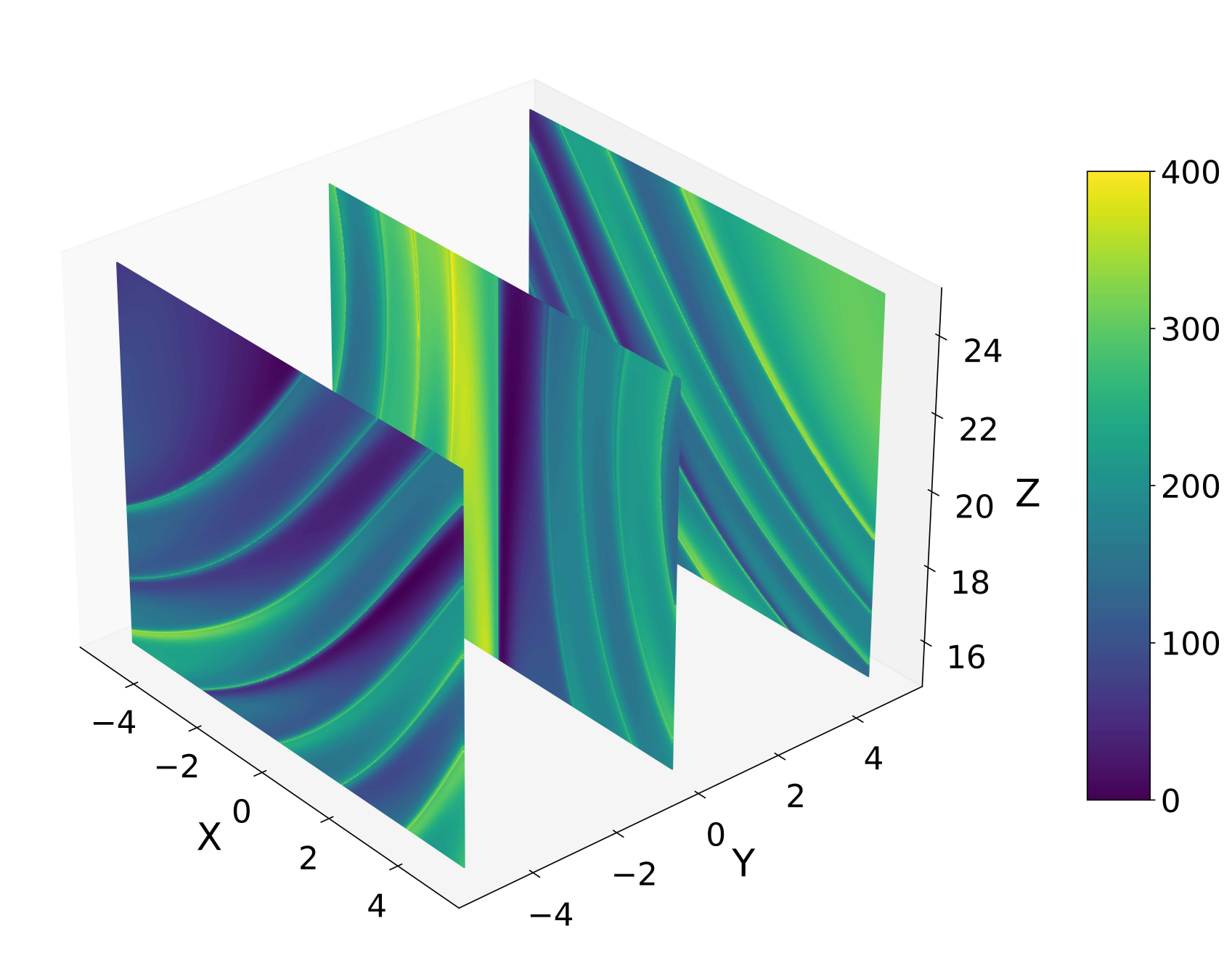}
\caption{$Y$-slices}
\end{subfigure}\\
\begin{subfigure}[t]{0.45\linewidth}
\centering
\includegraphics[scale=0.22]{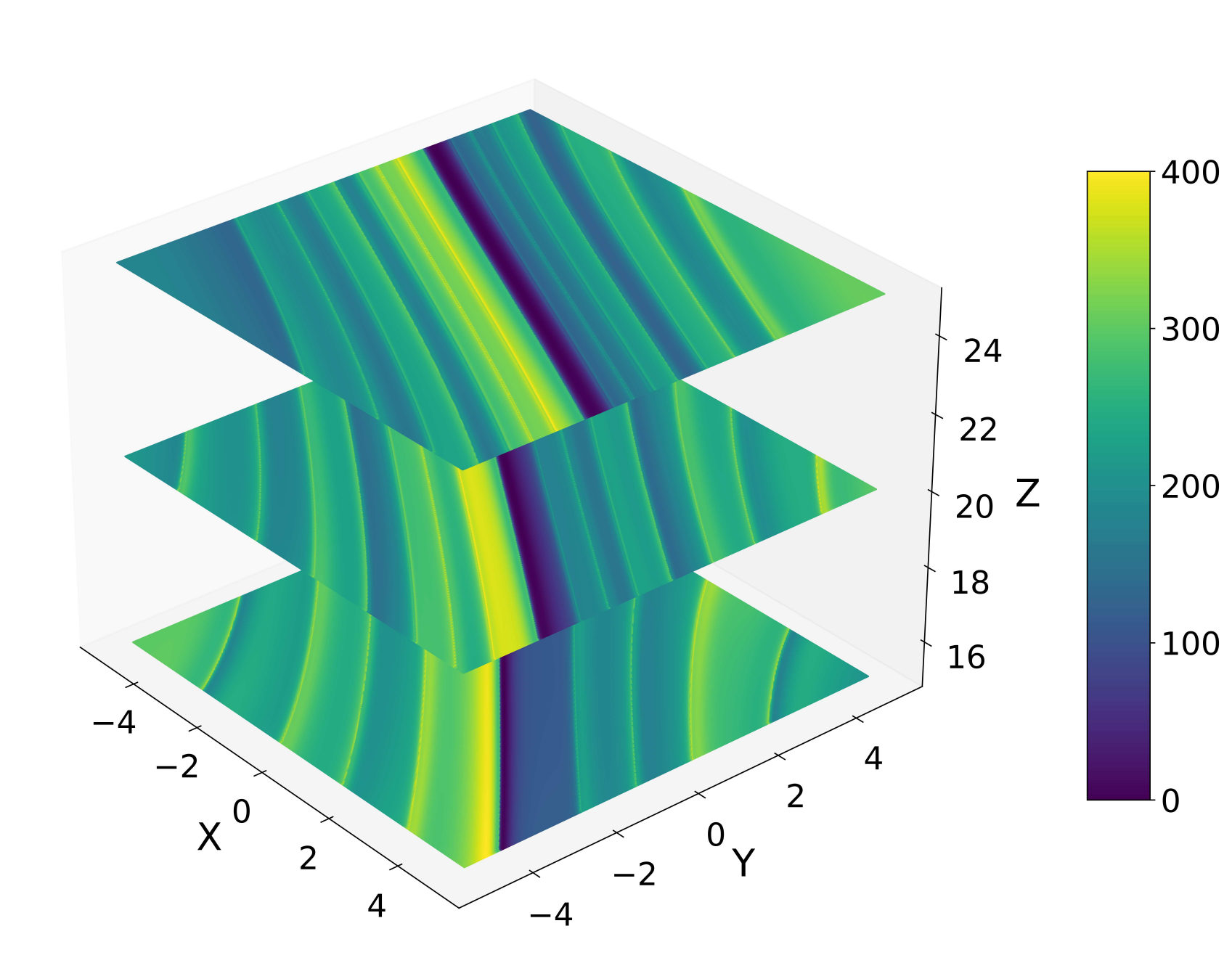}
\caption{$Z$-slices}
\end{subfigure}
\begin{subfigure}[t]{0.45\linewidth}
\centering
\includegraphics[scale=0.22]{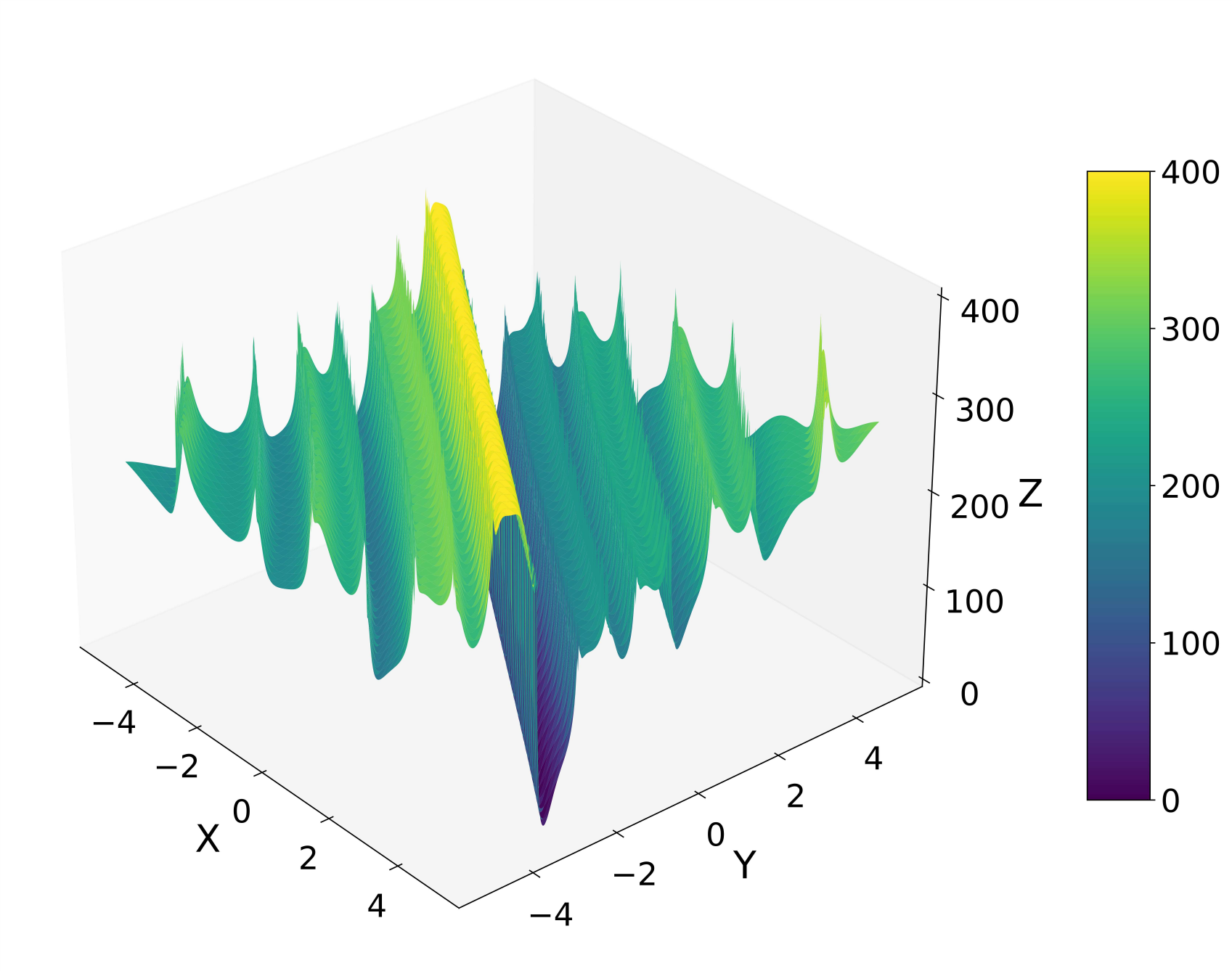}
\caption{Landscape of the second $Z$-slice}
\end{subfigure}
\caption{\small The objective function in~\eqref{eqn: numerical-4dvar-nonlinear} is evaluated at $\pmb{u}_0 = (x_0, y_0, z_0) \in [-6, 6] \times [-6, 6] \times [14, 26]$ for the Lorenz system described in~\eqref{eqn: lorenz}. In subfigures (a), (b), and (c), the colors correspond to the function values, with darker shades representing lower values and lighter shades indicating higher values. } 
\label{fig: non_cvx}
\end{figure}

Consider the classical Lorenz system, a well-known example of a nonlinear, aperiodic, and three-dimensional deterministic system known for its chaotic behavior, first studied by~\citet{lorenz1963deterministic}.  It is described by the following system of differential equations:
\begin{equation}
\label{eqn: lorenz}
\left\{ \begin{aligned}
         & \dot{x} = \sigma (y - x) \\
         & \dot{y} = x(\rho - z) - y \\
         & \dot{z} = xy - \beta z, 
         \end{aligned} \right.
\end{equation}
where the parameters are set to the classical values: $\sigma = 10$, $\rho=28$ and $\beta = 8/3$. To numerically solve the Lorenz system~\eqref{eqn: lorenz}, we employ the 4th-order Runge–Kutta method with a time step size of $\delta t = 0.01$ over a total time horizon of $T=3$. The observational interval is set to $T_o = 0.3$. Following the precise observation strategy proposed in~\citet{talagrand1987variational}, the numerical solution serves as the observational data. The observational data, denoted as $\hat{\pmb{u}}_{k} = \pmb{u}(kT_o)$ for $k=0,1,\ldots,n = T/T_o=10$, is generated using the 4th-order Runge-Kutta method for the Lorenz system~\eqref{eqn: lorenz} with the initial condition $\pmb{u}_0 = (-0.5, 0.5, 20.5)$.  In the context of the~\texttt{4D-Var} problem~\eqref{eqn: numerical-4dvar}, we set the regularization parameter $\alpha =0.1$. To illustrate the nonconvex nature of the objective function~\eqref{eqn: numerical-4dvar-nonlinear} associated with the Lorenz system~\eqref{eqn: lorenz}, we visualize slices of the objective function along the $X$-, $Y$-, and $Z$-directions, as shown in~\Cref{fig: non_cvx}. These visualizations highlight the nonconvex landscape of the objective function, which arises from the chaotic dynamics inherent to the Lorenz system~\eqref{eqn: lorenz}. As the system's nonlinearity and sensitivity to initial conditions increase, the optimization problem becomes significantly more challenging.

\begin{figure}[htb!]
\centering
\begin{subfigure}[t]{0.45\linewidth}
\centering
\includegraphics[scale=0.2]{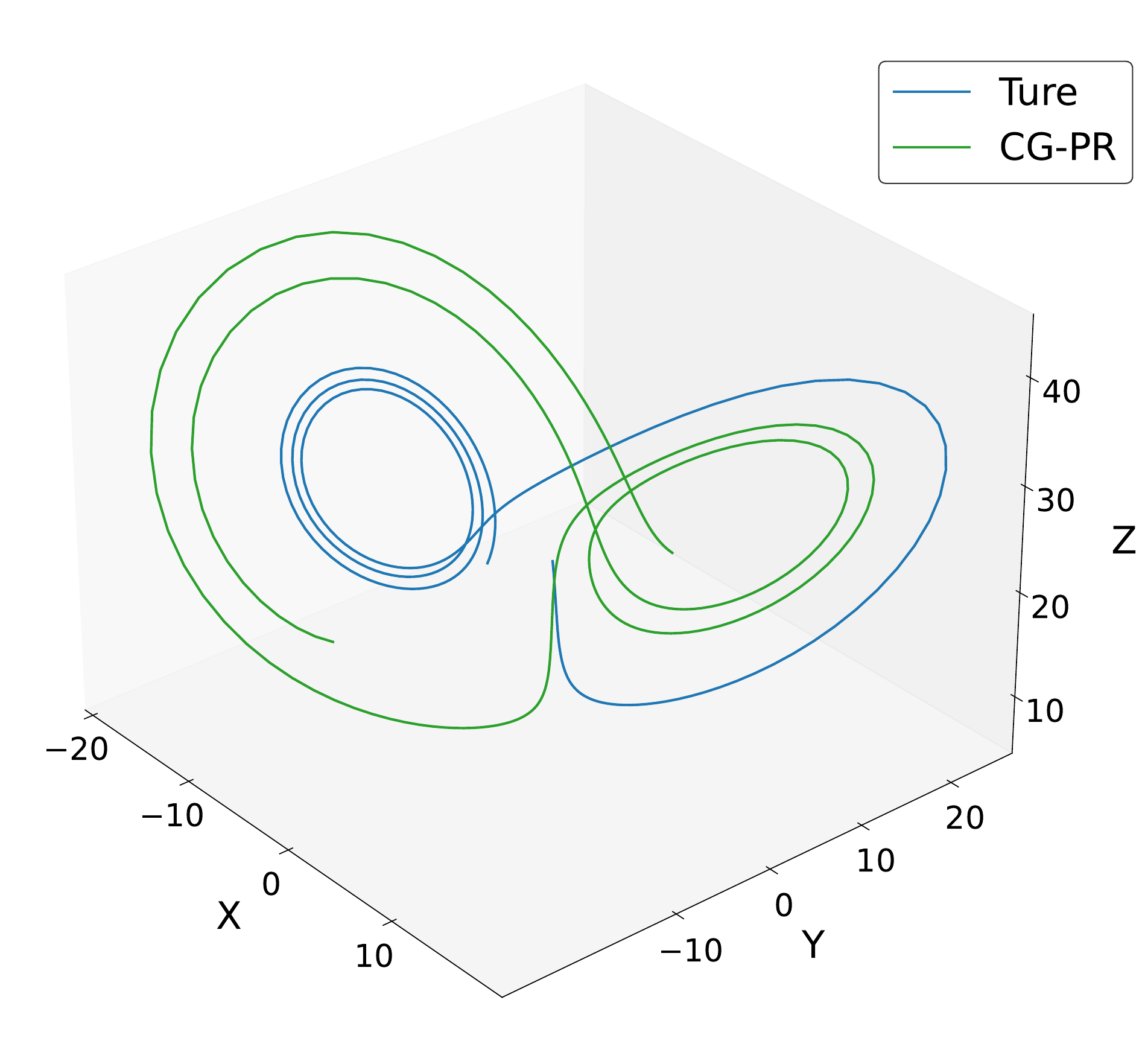}
\caption{Conjugate gradient method}
\end{subfigure}
\begin{subfigure}[t]{0.45\linewidth}
\centering
\includegraphics[scale=0.2]{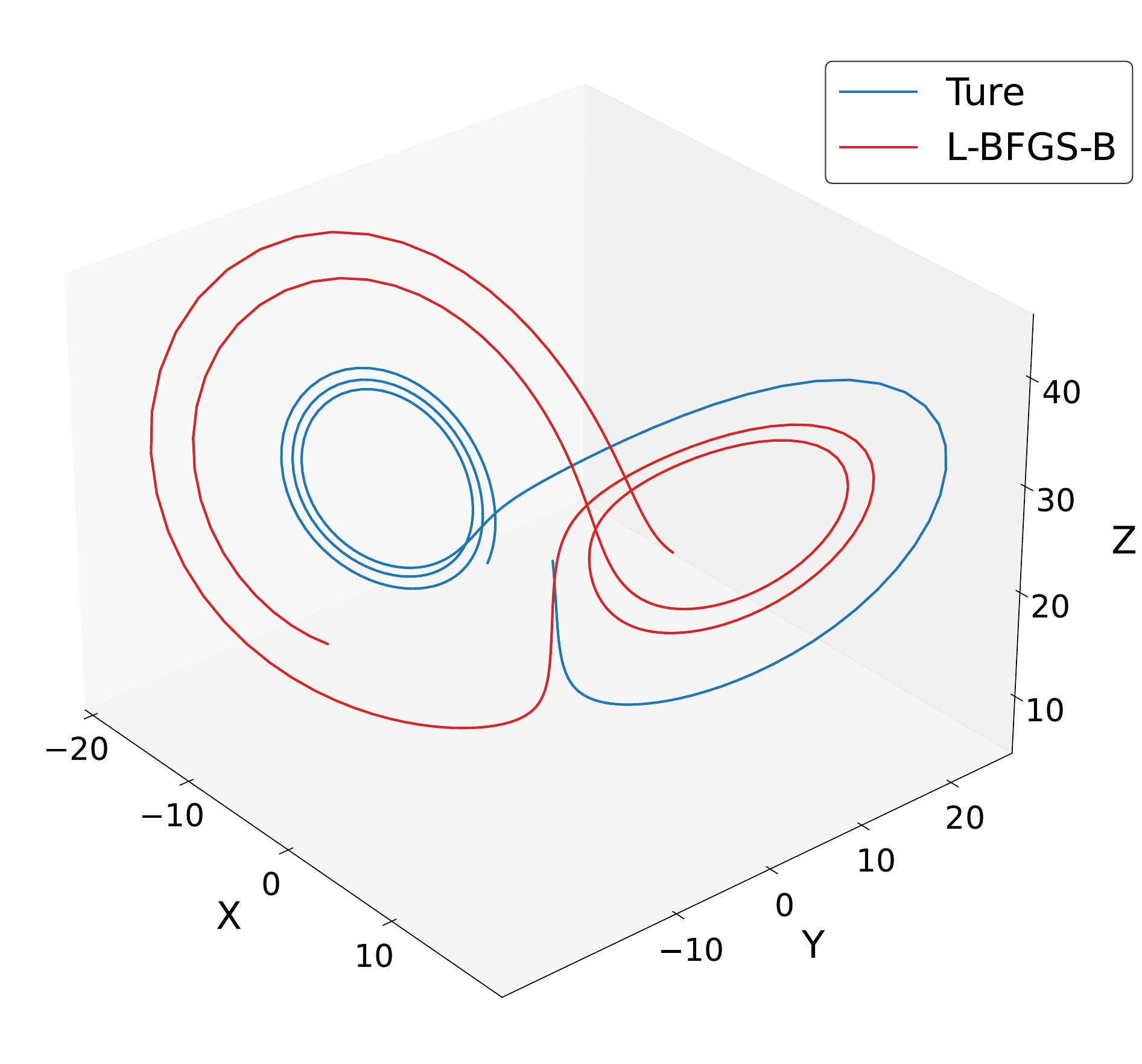}
\caption{Limited memory BFGS-B method}
\end{subfigure}
\caption{\small The optimization algorithms are applied to the~\texttt{4D-Var} problem~\eqref{eqn: numerical-4dvar-nonlinear} for the Lorenz system~\eqref{eqn: lorenz}, starting with the initial condition $\pmb{u}_0 = (x_0, y_0, z_0)  = (-3,-3, 10)$. The numerical implementation is carried out using the  Python subpackage --- SciPy 1.14.1. } 
\label{fig: lorenz_cg-bfgs}
\end{figure}

%Next, we briefly outline the process of using the so-called adjoint method to derive the gradient of the objective function in the \texttt{4D-Var} problem~\eqref{eqn: numerical-4dvar-nonlinear}, as described in~\citep{le1986variational}. 
We then apply two classical first-order optimization algorithms, the nonlinear conjugate gradient methods~\citep{hestenes1952methods, fletcher1964function, M2AN_1969__3_1_35_0, dai1999nonlinear} and the limited-memory BFGS-B method~\citep{byrd1995limited, zhu1997algorithm}, to tackle the~\texttt{4D-Var} problem~\eqref{eqn: numerical-4dvar}. For the Lorenz system~\eqref{eqn: lorenz}, its tangent linear system is derived as:
\begin{equation}
\label{eqn: tlm-lorenz}
\begin{pmatrix}
\delta\dot{x} \\
\delta\dot{y} \\
\delta\dot{z}
\end{pmatrix} = K[x,y,z] \begin{pmatrix}
\delta x \\
\delta y \\
\delta z
\end{pmatrix} =
\begin{pmatrix}
-\sigma & \sigma & 0         \\
\rho - z  & -1        & - x       \\
y           & x          & - \beta
\end{pmatrix} 
\begin{pmatrix}
\delta x \\
\delta y \\
\delta z
\end{pmatrix}. 
\end{equation}
Utilizing the expression~\eqref{eqn: jacobian} along with the 4th-order Runge-Kutta iteration, we compute the Jacobian of the nonlinear operator $\nabla H_{T_o}^{k}(\pmb{u}_0)$. Taking its transpose, $(\nabla H_{T_o}^{k}(\pmb{u}_0))^{\top}$, allows us to obtain the gradient required for optimization.  However, when these algorithms are applied to the objective function~\eqref{eqn: numerical-4dvar-nonlinear}, both fail to converge the true solution, as demonstrated in~\Cref{fig: lorenz_cg-bfgs}. This failure is attributed to the highly nonconvex nature of the objective function, which is riddled with multiple local minima.~\Cref{fig: non_cvx} illustrates the landscape of the objective function, characterized by numerous valleys and peaks. While the global minimum corresponds to the true solution, the nonconvex landscape traps both algorithms in local minima, preventing them from reaching the optimal solution. 
This is particularly evident in~\Cref{fig: lorenz_cg-bfgs}, where the final solutions produced by both methods remain far from the global minimum. The inherent nonconvexity of the problem poses significant challenges. These classical first-order optimization methods, such as the nonlinear conjugate gradient methods and the Limited memory BFGS-B method, which are typically effective for convex problems, struggle to escape local minima. Furthermore, the nonconvex nature of the~\texttt{4D-Var} problem~\eqref{eqn: numerical-4dvar-nonlinear} makes the solutions obtained by these classical first-order optimization methods highly sensitive to initial conditions. This sensitivity, as noted in~\citet{talagrand20144dvar},  is a well-known open problem that requires further in-depth investigation.

%Both algorithms are designed to handle large-scale optimization problems efficiently.
%however, due to the nonconvex nature of the landscape, both algorithms become trapped in local minima. This behavior is depicted in~\Cref{fig: lorenz_cg-bfgs}, where the final solutions reached by both methods remain far from the global minimum. 

%%%%%%%%%%%%%%%%%%%%%%%%%%%%%%%%%%%%%%%%%%%%%%%%%%%%%%%%%%%%%%%%%%%%%%%%%%%%%%%%%%%%%%%%%%%%%%%%%%%%%%%%%%%%%%%%%%%
\subsection{4D-Var for PDEs}
\label{subsec: 4DVar-PDEs}

In this section, we provide a brief outline of the \texttt{4D-Var} problem for partial differential equations (PDEs) and its numerical implementation.  Let $\Omega \subseteq \mathbb{R}^d$ be an open set with a smooth boundary $\partial \Omega$. The general form of a PDE governing the dynamics of a system can be written as:
\begin{equation}
\label{eqn: infinite-pde}
\left\{ \begin{aligned}
          & u_t = G(u)                                 \\
          & u |_{t=0} = u_0                         \\
          & u |_{\partial \Omega} = g,   
           \end{aligned} \right.
\end{equation}
where $u = u(t, x)$ is a function of time $t \in (0, +\infty)$ and the spatial variable $x \in \mathbb{R}^d$, and its Energy norm, or $L^2$-norm, is defined as:
\begin{equation}
\label{eqn: energy-norm}
\|u(t, \cdot)\|^2 = \int_{\Omega} |u(t,x)|^2dx.
\end{equation}
The numerical implementation of PDEs given by~\eqref{eqn: infinite-pde} often involves discretization techniques that transform the continuous system into a more tractable, finite-dimensional one. The most widely used methods include the finite difference method, finite element method, and spectral methods. While these methods differ in terms of theoretical error bounds and convergence rates, they share the goal of approximating the continuous PDE by converting it into a system of ordinary differential equations (ODEs). The solution to this discretized system, at the observational time points, can then be iteratively expressed through a discrete scheme in alignment with the observational data as outlined in~\eqref{eqn: finite-ds-discrete}. 

For the~\texttt{4D-Var} problem, which involves data assimilation over a specific time window, the numerical implementation focuses on minimizing a cost function that quantifies the misfit between model predictions and observations, as outlined in~\eqref{eqn: numerical-4dvar}. This minimization problem can often be reformulated as a quadratic function of the initial condition $u_0$,  as shown in~\eqref{eqn: numerical-4dvar-linear}, allowing for efficient solution of the linearized \texttt{4D-Var} problem. Classical first-order optimization algorithms, such as the conjugate gradient method~\citep{hestenes1952methods} and the limited-memory BFGS-B method~\citep{byrd1995limited, zhu1997algorithm}, are well-suited for this task.  These optimization techniques have been successfully applied to a wide range of linear PDEs, including the heat equation~\citep{burman2018data, li2024variational} and the wave equation~\citep{burman2020finite}.

%In practice, numerical discretization of the PDE~\eqref{eqn: infinite-pde} can be carried out using various methods, such as finite difference, finite element, and spectral methods. Although these methods differ in their theoretical error analysis, they all approximate the continuous PDE by converting it into a finite-dimensional system of ordinary differential equations (ODEs). 
%
%For the~\texttt{4D-Var} problem, the numerical implementation involves formulating the objective function as described in~\eqref{eqn: numerical-4dvar}. This problem can also be reformulated as a quadratic function of the initial condition $\pmb{u}_0$, as shown in~\eqref{eqn: numerical-4dvar-linear}. This reformulation enables the efficient solution of the linear \texttt{4D-Var} problem using .

%% file: 02_admm.tex
\section{Solving the 4D-Var problem via ADMM}
\label{sec: outline-admm}

In this section, we introduce a practical variant of the alternating direction method of multipliers (\texttt{ADMM}), known as multi-block~\texttt{ADMM}, to solve the~\texttt{4D-Var} problem~\eqref{eqn: numerical-4dvar}.  The use of multi-block~\texttt{ADMM} arises naturally from reformulating the \texttt{4D-Var} problem~\eqref{eqn: numerical-4dvar} as a constraint optimization problem, which facilitates the derivation of its augmented Lagrangian. By efficiently implementing the linearized version of multi-block~\texttt{ADMM} with regularization, we demonstrate strong numerical performance. Finally, we highlight several advantages of this approach over classical first-order unconstrained optimization algorithms, as discussed in~\Cref{sec: intro}.  

% being given as follows. Finally, we propose some advantages of the  linearized multi-block~\texttt{ADMM} with regularization, compared with the classical first-order unstrained optimization algorithms, as 

%%%%%%%%%%%%%%%%%%%%%%%%%%%%%%%%%%%%%%%%%%%%%%%%%%%%%%%%%%%%%%%%%%%%%%%%%%%%%%%%
\subsection{Motivation and basic idea for implementing ADMM}
\label{subsec: motivation}

To facilitate the integration of observational data with the system's dynamics, it is necessary to expand the number of variables in the objective function without alternating its values.  This expansion enables a more accurate representation of the system's state over time.  Given a total number of iterations $N = T/\delta t$ and $n = T_o/\delta t$ iterations per observational time interval, we increase the number of variables from $n$ to $N$. This expansion allows for closer tracking of the state of the system over each discrete time step. For any $k =0,1, \ldots, N$, the sub-objective function $f_k(\pmb{u}_k)$ is defined as follows: 
\begin{equation}
\label{eqn: objective-admm}
f_k(\pmb{u}_k) = \left\{ \begin{aligned} 
                                &\frac{T_o}{2}\|\pmb{u}_0 - \hat{\pmb{u}}_0\|^2 + \frac{\alpha}{2} \|\pmb{u}_0 - \hat{\pmb{u}}_0^b\|^2,  && \quad \text{for}\;\; k=0, \\
                                &\frac{T_o}{2}\big\| \pmb{u}_{k/n} - \hat{\pmb{u}}_{k/n}  \big\|^2,                                                              && \quad \text{for}\;\; k/n \in  \mathbb{N}^+, \\
                                &0,                                                                                                                                                             && \quad \text{otherwise}.
                                \end{aligned} \right.
\end{equation}
The overall~\texttt{4D-Var} problem~\eqref{eqn: numerical-4dvar} can then be reformulated as:
\begin{equation}
\label{eqn: numerical-4dvar-admm}
\left\{ \begin{aligned}
         & \min F(\pmb{u}_0, \pmb{u}_1, \ldots, \pmb{u}_N) := \sum_{k=0}^{N} f_k(\pmb{u}_k), \\
         & \text{s.t.} \quad \pmb{u}_{k+1} = H( \pmb{u}_{k}), \quad \text{for} \quad k = 0, 1, \ldots, N-1,
          \end{aligned} \right.
\end{equation}
where $H$ represents the nonlinear operator that governs the evolution of the system within some numerical scheme at each time step $\delta t$. 

\paragraph{Augmented Lagrangian Formulation} Based on the constraint optimization form of the \texttt{4D-Var} problem~\eqref{eqn: numerical-4dvar-admm}, we can then derive its augmented Lagrangian as follows:
\begin{align}
L & = L(\pmb{u}_0, \pmb{u}_1, \ldots, \pmb{u}_N; \pmb{\lambda}_0,\ldots,\pmb{\lambda}_{N-1};s) \nonumber \\ 
   & = \sum_{k=0}^{N}f_k(\pmb{u}_k) - \sum_{k=0}^{N-1}\left\langle \pmb{\lambda}_k, \pmb{u}_{k+1} - H(\pmb{u}_k) \right\rangle + \frac{1}{2s} \sum_{k=0}^{N-1} \left\| \pmb{u}_{k+1} - H(\pmb{u}_k) \right\|^2  \label{eqn: numerical-4dvar-lagrangian}, 
\end{align}
where $s>0$ is a given parameter. By completing the square on the penalty terms, this formulation~\eqref{eqn: numerical-4dvar-lagrangian} can be simplified further: 
 \begin{equation}
\label{eqn: numerical-4dvar-lagrangian-reform}
L = \sum_{k=0}^{N}f_k(\pmb{u}_k) +\frac{1}{2s} \sum_{k=0}^{N-1} \left\| \pmb{u}_{k+1} - H(\pmb{u}_k) - s\pmb{\lambda}_k\right\|^2 - \frac{s}{2} \sum_{k=0}^{N-1} \|\pmb{\lambda}_k\|^2.
\end{equation}
From the simplified expression~\eqref{eqn: numerical-4dvar-lagrangian-reform}, we observe that $\pmb{u}_k$ is only related to its neighboring time steps, $\pmb{u}_{k-1}$ and $\pmb{u}_{k+1}$, since the sub-objective function $f_k$ depends only on $\pmb{u}_k$. This local dependence simplifies the optimization process by focusing on interactions with adjacent time steps in the iterative solution. 

\paragraph{Multi-block ADMM} To demonstrate the iterative process of multi-block~\texttt{ADMM}, we apply the $\argmin$ operation to the augmented Lagrangian during the following iteration:
\begin{equation}
\label{eqn: argmin-4dvar-alm}
\left\{ \begin{aligned}
         & \pmb{u}_{k}^{\ell+1} = \argmin_{\pmb{u}_k} L(\pmb{u}_0^{\ell}, \ldots, \pmb{u}_{k-1}^{\ell},\pmb{u}_k, \pmb{u}_{k+1}^{\ell}, \ldots, \pmb{u}_N^{\ell}; \pmb{\lambda}_0^{\ell},\ldots,\pmb{\lambda}_{N-1}^{\ell};s), && \text{for}\;\; k = 0, \ldots, N,   \\
         & s \pmb{\lambda}^{\ell+1}_{k} = s \pmb{\lambda}^{\ell}_k -  \left( \pmb{u}_{k+1}^{\ell+1} - H(\pmb{u}_k^{\ell+1}) \right), && \text{for}\;\; k = 0, \ldots, N-1,
         \end{aligned} \right.
\end{equation}
where the second equation represents the update for the dual variables (Lagrange multipliers). This procedure alternates between updating the primal variables $\pmb{u}_k$ and the dual variables $\pmb{\lambda}_k$. Substituting the simplified expression~\eqref{eqn: numerical-4dvar-lagrangian-reform} into the $\argmin$ iteration~\eqref{eqn: argmin-4dvar-alm}, we can derive that the multi-block~\texttt{ADMM} iterates as follows:
\begin{equation}
\label{eqn: ADMM}
	\left\{ \begin{aligned}
		& \pmb{u}_0^{\ell+1} = \argmin_{\pmb{u}_0} \left\{f_0(\pmb{u}_0) + \frac{1}{2s}\norm{\pmb{u}_1^{\ell} - H(\pmb{u}_0) - s\pmb{\lambda}_0^{\ell}}^2  \right\},\\
%		& \pmb{u}_1^{\ell+1} = \argmin_{\pmb{u}_1} \left\{f_1(\pmb{u}_1) + \frac{1}{2s}\norm{\pmb{u}_1 - H(\pmb{u}_{0}^{\ell}) - s\pmb{\lambda}_{0}^{\ell} }^2 + \frac{1}{2s}\norm{\pmb{u}_{2}^{\ell} - H(\pmb{u}_1) - s\pmb{\lambda}_1^{\ell}}^2  \right\}, \\
%		& \qquad \qquad \quad \vdots \\
		& \pmb{u}_k^{\ell+1} = \argmin_{\pmb{u}_i} \left\{f_k(\pmb{u}_k) + \frac{1}{2s}\norm{\pmb{u}_k - H(\pmb{u}_{k-1}^{\ell}) - s\pmb{\lambda}_{k-1}^{\ell}}^2 + \frac{1}{2s}\norm{\pmb{u}_{k+1}^{\ell} - H(\pmb{u}_k) - s\pmb{\lambda}_k^{\ell}}^2  \right\}, \\
                   & \mathrel{\phantom{\pmb{u}_{N-1}^{\ell+1} = \argmin_{\pmb{u}_{N-1}} \left\{f_{N-1}(\pmb{u}_{N-1})\right.abcdefg}} \mathrm{for}\;\;k=1,2,\ldots, N-1, \\
%		& \qquad \qquad \quad \vdots \\
%		& \pmb{u}_{N-1}^{\ell+1} = \argmin_{\pmb{u}_{N-1}} \left\{f_{N-1}(\pmb{u}_{N-1}) + \frac{1}{2s}\norm{\pmb{u}_{N-1} - H(\pmb{u}_{N-2}^{\ell}) - s\pmb{\lambda}_{N-2}^{\ell} }^2\right.\\	
%		& \mathrel{\phantom{\pmb{u}_{N-1}^{\ell+1} = \argmin_{\pmb{u}_{N-1}} \left\{f_{N-1}(\pmb{u}_{N-1})\right.abcdefghijklmnopq}} \left. + \frac{1}{2s}\norm{\pmb{u}_{N}^{\ell} - H(\pmb{u}_{N-1}) - s\pmb{\lambda}_{N -1}^{\ell}}^2   \right\}, \\
		& \pmb{u}_N^{\ell+1} = \argmin_{\pmb{u}_N} \left\{f_N(\pmb{u}_N) + \frac{1}{2s}\norm{\pmb{u}_N - H(\pmb{u}_{N-1}^{\ell}) - s\pmb{\lambda}_{N-1}^{\ell}}^2  \right\},\\
		&\\
		& s \pmb{\lambda}_k^{\ell +1} = s\pmb{\lambda}_k^{\ell} - \left(\pmb{u}_{k+1}^{\ell +1} - H(\pmb{u}_k^{\ell+1})\right), \qquad \text{for}\;\; k = 0, 1,\ldots, N-1.
	\end{aligned} \right.
\end{equation}
In this formulation, it is observed that each block $\pmb{u}_k$ is updated independently, which is the key characteristic of the multi-block \texttt{ADMM} algorithm. While the classical two-block~\texttt{ADMM} was originally proposed in~\citep{glowinski1975approximation, gabay1976dual}, the multi-block~\texttt{ADMM} was first studied by~\citep{he2015full} in the context of linear constraints. It is important to note that, even in the two-block case, the algorithm described in~\eqref{eqn: ADMM} may exhibit divergence, as discussed in~\citep{he2015full}. The potential for divergence is a known issue in multi-block \texttt{ADMM}, particularly when the assumptions regarding the problem structure or linear constraints are not satisfied.

\subsection{Efficient ADMM implementation}
\label{subsec: admm-implementation}

As demonstrated in~\citep{he2015full}, the introduction of a regularization term is a common approach to guarantee the convergence of the multi-block~\texttt{ADMM}~\eqref{eqn: ADMM}. Originally proposed by~\citet{zhang2010bregmanized} for the two-block case, this regularization technique is crucial for stabilizing the iterative updates and mitigating potential divergence, particularly when the underlying problem structure does not inherently guarantee convergence. Additionally, each subproblem for $k=0,1,\ldots, N-1$ involves solving a nonlinear least square problem, $\| \pmb{u}_{k+1}^{\ell} - H(\pmb{u}_k) - s\pmb{\lambda}_k^{\ell} \|^2$, which incurs significant computational overhead. 

\paragraph{Linearized multi-block ADMM with regularization}
To reduce computational costs and enhance efficiency, a linearization technique introduced by~\citet{deng2016global} is employed. This technique, combined with regularization, allows for a modified version of the multi-block~\texttt{ADMM}, expressed as follows:
\begin{equation}
\label{eqn: lin-ADMM-reg}
	\left\{ \begin{aligned}
		& \pmb{u}_0^{\ell+1} = \argmin_{\pmb{u}_0} \left\{f_0(\pmb{u}_0) - \frac{1}{s} \left\langle \nabla H(\pmb{u}_0^{\ell})^{\top} \left(\pmb{u}_1^{\ell}-H(\pmb{u}_0^{\ell}) - s\pmb{\lambda}_0^{\ell}\right), \pmb{u}_0 \right\rangle+ \frac{\norm{\pmb{u}_0 - \pmb{u}_0^{\ell}}^2}{2\eta}  \right\},\\
%		& \pmb{u}_1^{\ell+1} = \argmin_{\pmb{u}_1} \left\{f_1(\pmb{u}_1) + \frac{1}{2s}\norm{\pmb{u}_1 - H(\pmb{u}_{0}^{\ell}) - s\pmb{\lambda}_{0}^{\ell} }^2 \right. \\
%		& \mathrel{\phantom{\pmb{u}_{i}^{\ell+1} = \argmin_{\pmb{u}_{i}} \left\{f_{i}(\pmb{u}_{i})\right.}} \left. -\frac{1}{s} \left\langle \nabla H(\pmb{u}_1^{\ell})^{\top}\left(\pmb{u}_2^{\ell}-H(\pmb{u}_1^{\ell}) - s\pmb{\lambda}_1^{\ell}\right), \pmb{u}_1 \right\rangle + \frac{\norm{\pmb{u}_1 - \pmb{u}_1^{\ell}}^2}{2\eta}  \right\}, \\
%		& \qquad \qquad \quad \vdots \\
		& \pmb{u}_k^{\ell+1} = \argmin_{\pmb{u}_k} \left\{f_k(\pmb{u}_k) + \frac{1}{2s}\norm{\pmb{u}_k - H(\pmb{u}_{k-1}^{\ell}) - s\pmb{\lambda}_{k-1}^{\ell} }^2 \right. \\
		& \mathrel{\phantom{\pmb{u}_{k}^{\ell+1} = \argmin_{\pmb{u}_{k}} \left\{f_{k}(\pmb{u}_{k})\right.}} \left. - \frac{1}{s} \left\langle \nabla H(\pmb{u}_k^{\ell})^{\top}\left(\pmb{u}_{k+1}^{\ell}-H(\pmb{u}_k^{\ell}) - s\pmb{\lambda}_{k}^{\ell} \right), \pmb{u}_{k} \right\rangle + \frac{\norm{\pmb{u}_k - \pmb{u}_k^{\ell}}^2}{2\eta}  \right\}, \\
                   & \mathrel{\phantom{\pmb{u}_{N-1}^{\ell+1} = \argmin_{\pmb{u}_{N-1}} \left\{f_{N-1}(\pmb{u}_{N-1})\right.abcdefg}} \mathrm{for}\;\;k=1,2,\ldots, N-1, \\
%		& \qquad \qquad \quad \vdots \\
%		& \pmb{u}_{N-1}^{\ell+1} = \argmin_{\pmb{u}_{N-1}} \left\{f_{N-1}(\pmb{u}_{N-1}) + \frac{1}{2s}\norm{\pmb{u}_{N-1} - H(\pmb{u}_{N-2}^{\ell}) - s\pmb{\lambda}_{N-2}^{\ell} }^2 \right.\\
%		& \mathrel{\phantom{\pmb{u}_{N-1}^{\ell+1} ar}} \left. -\frac{1}{s} \left\langle \nabla H(\pmb{u}_{N-1}^{\ell})^{\top} \left(\pmb{u}_{N}^{\ell} - H(\pmb{u}_{N-1}^{\ell}) - s \pmb{\lambda}_{N-1}^{\ell} \right), \pmb{u}_{N-1} \right\rangle  + \frac{\norm{\pmb{u}_{N-1} - \pmb{u}_{N-1}^{\ell}}^2}{2\eta}   \right\}, \\
		& \pmb{u}_N^{\ell+1} = \argmin_{\pmb{u}_N} \left\{f_N(\pmb{u}_N) + \frac{1}{2s}\norm{\pmb{u}_N - H(\pmb{u}_{N-1}^{\ell}) - s\pmb{\lambda}_{N-1}^{\ell} }^2 + \frac{\norm{\pmb{u}_{N} - \pmb{u}_{N}^{\ell}}^2}{2\eta}   \right\},\\
		&\\
		& s\pmb{\lambda}_k^{\ell+1} = s\pmb{\lambda}_k^{\ell} - \left(\pmb{u}_{k+1}^{\ell+1} - H(\pmb{u}_k^{\ell+1})\right), \qquad \text{for}\;\; k = 0, \ldots, n-1,
	\end{aligned} \right.
\end{equation}
where $\eta>0$ is a regularization parameter that helps control the deviation between successive iterations. In this formulation~\eqref{eqn: lin-ADMM-reg}, the linearization involves computing the gradient, which can still be efficiently obtained using the adjoint method~\citep{le1986variational,talagrand1987variational}.  The regularization term ensures that successive iterations do not deviate too far from one another, stabilizing the updates and preventing large, destabilizing steps. Thus, the linearization serves as an effective approximation for the nonlinear least square problem. Recently, several studies, including~\citet{xie2021complexity, el2023linearized}, and~\citet{hien2024inertial}, have explored the convergence properties of the linearized multi-block~\texttt{ADMM} with regularization~\eqref{eqn: lin-ADMM-reg}. The convergence analysis of this method has emerged as a compelling area of research.

%From this formulation~\eqref{eqn: lin-ADMM-reg}, we observe that the linearization technique, when coupled with regularization, effectively stabilizes the iterative updates. 

%However, solving the resulting subproblems often incurs significant computational costs. To address this issue,~\citet{deng2016global} introduced a linearization technique employed to improve computational efficiency. By combining regularization and linearization, the modified multi-block~\texttt{ADMM} can be expressed as follows:

%%%%%%%%%%%%%%%%%%%%%%%%%%%%%%%%%%%%%%%%%%%%%%%%%%%%%%%%%%%%%%%%%%%%%%%%%%%%%%%%%%%%%%%%%%%%%%%%%%%%%%%%%%%%%%%%%%%
\subsection{Numerical Performance}
\label{subsec: num-phenomena}

To demonstrate the high numerical efficiency of the linearized multi-block~\texttt{ADMM} with regularization~\eqref{eqn: lin-ADMM-reg}, we use the Lorenz system~\eqref{eqn: lorenz} as a benchmark for comparison with the classical first-order optimization algorithms, as outlined in~\Cref{sec: intro}. The numerical solution is obtained using the 4th-order Runge-Kutta method with a time step size of $\delta t = 0.01$ over a total time horizon of $T = 3$. The observational interval is set to $T_o = 0.3$, leading to a total number of $N = T/\delta t = 300$ iterations, with $n = T/T_o +1 = 11$ observational time points. We present the numerical performance of the linearized multi-block~\texttt{ADMM} with regularization~\eqref{eqn: lin-ADMM-reg} applied to the~\texttt{4D-Var} problem~\eqref{eqn: numerical-4dvar} under both precise and noisy observation conditions, highlighting its effectiveness across varying scenarios.

\begin{figure}[htb!]
\centering
\begin{subfigure}[t]{0.45\linewidth}
\centering
\includegraphics[scale=0.18]{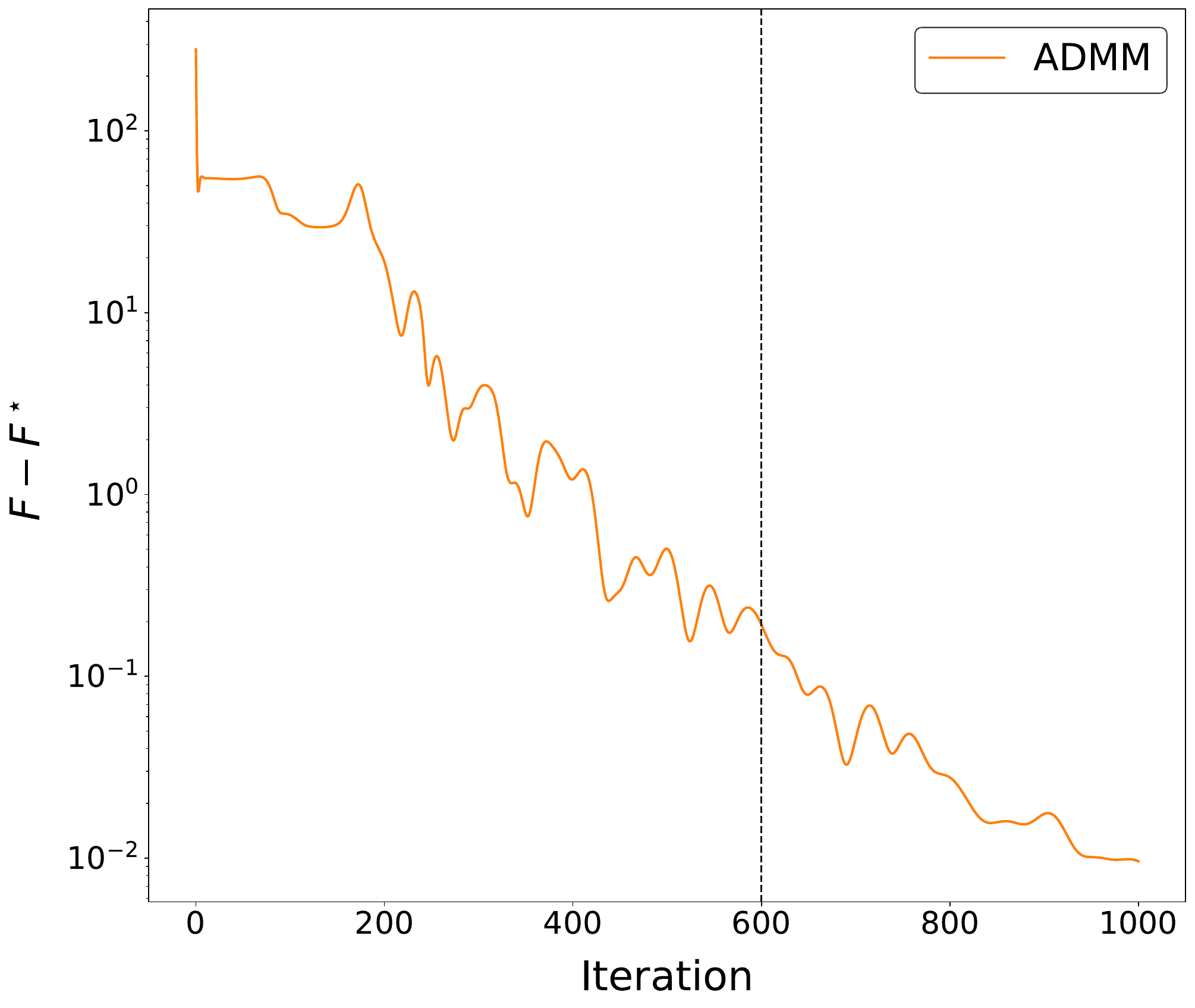}
\caption{\small Total error}
\label{subfig: precise-total-error}
\end{subfigure}
\begin{subfigure}[t]{0.45\linewidth}
\centering
\includegraphics[scale=0.18]{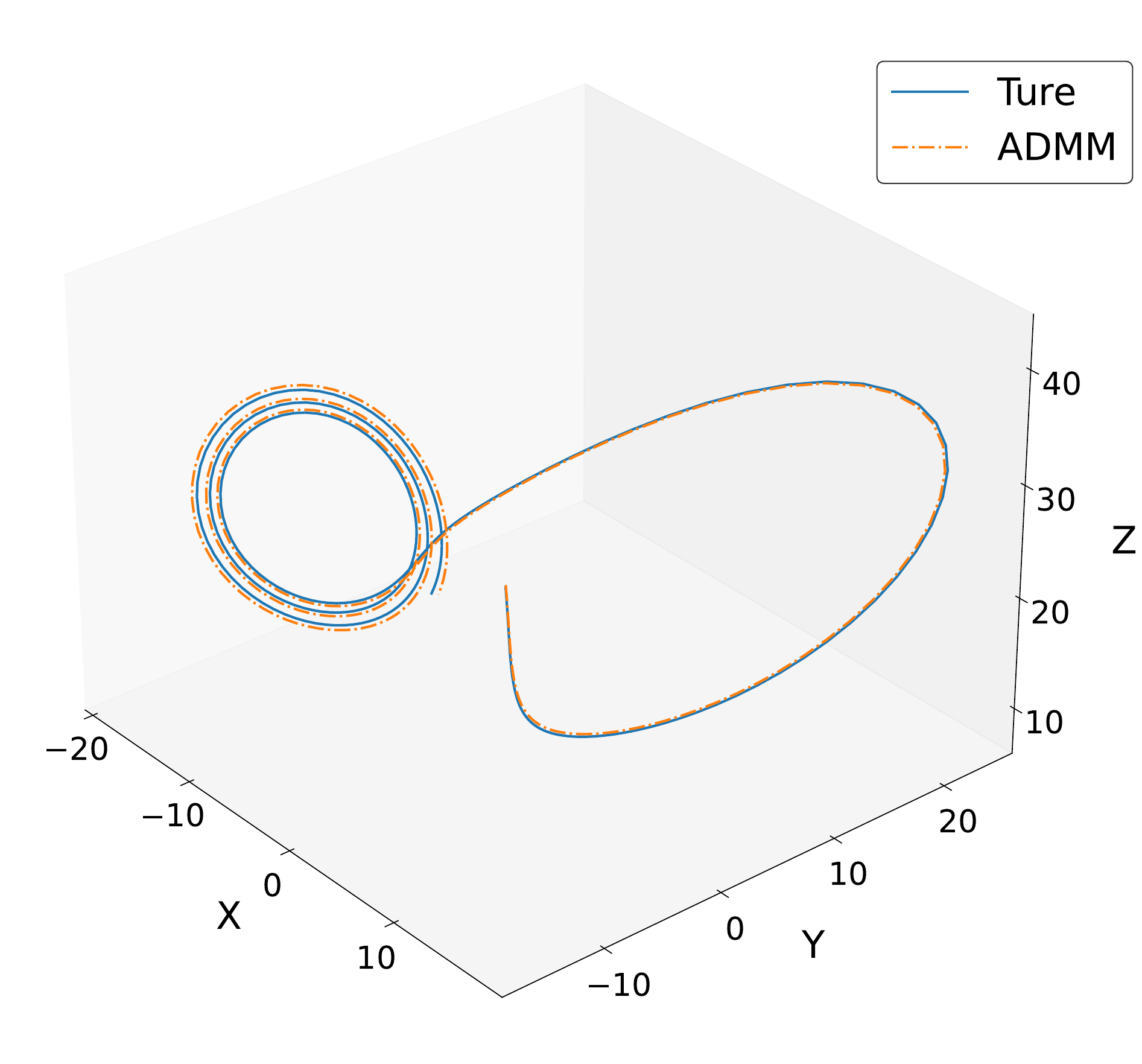}
\caption{\small Recovery solution ($\ell=600$)}
\label{subfig: precise-recover-soln}
\end{subfigure}
\caption{\small Numerical performance of the linearized multi-block~\texttt{ADMM} with regularization~\eqref{eqn: lin-ADMM-reg}, which is applied to the~\texttt{4D-Var} problem~\eqref{eqn: numerical-4dvar-nonlinear} for the Lorenz system~\eqref{eqn: lorenz} under the same settings as depicted in~\Cref{fig: lorenz_cg-bfgs}. } 
\label{fig: lorenz_admm}
\end{figure}

\paragraph{Precise observation}
The precise observational data, $\left( \hat{x}_1(k), \hat{y}_1(k), \hat{z}_1(k) \right)$ for $k=0,1,\ldots,n=10$, is generated using the numerical solution from the 4th-order Runge-Kutta method, recorded at $kT_o$. To solve the~\texttt{4D-Var} problem~\eqref{eqn: numerical-4dvar}, we apply the linearized multi-block~\texttt{ADMM} with regularization~\eqref{eqn: lin-ADMM-reg}, with the parameters set as $\mu = 100$, $\eta = 0.1$, and $s=2/3$. As outlined in~\Cref{sec: intro}, the classical first-order optimization algorithms begin with the initial condition $\pmb{u}_0 = (x_0, y_0, z_0)  = (-3,-3, 10)$.  To ensure consistency in comparison, the linearized multi-block~\texttt{ADMM} with regularization~\eqref{eqn: lin-ADMM-reg} also starts with the numerical solution recorded at $kT_o$ for $k=0,1,\ldots,n=10$, which is obtained using the 4th-order Runge-Kutta method under the same initial condition. The numerical performance is illustrated in~\Cref{fig: lorenz_admm}.
Unlike the classical first-order optimization algorithms, such as the nonlinear conjugate gradient method and the Limited memory BFGS-B method, which are prone to get trapped in local minima (see~\Cref{fig: lorenz_cg-bfgs}), the linearized multi-block~\texttt{ADMM} with regularization~\eqref{eqn: lin-ADMM-reg} exhibits robust convergence. It consistently reaches the true solution, as evidenced by the numerical performance shown in~\Cref{fig: lorenz_admm}.

\begin{figure}[htb!]
\centering
\begin{subfigure}[t]{0.32\linewidth}
\centering
\includegraphics[scale=0.15]{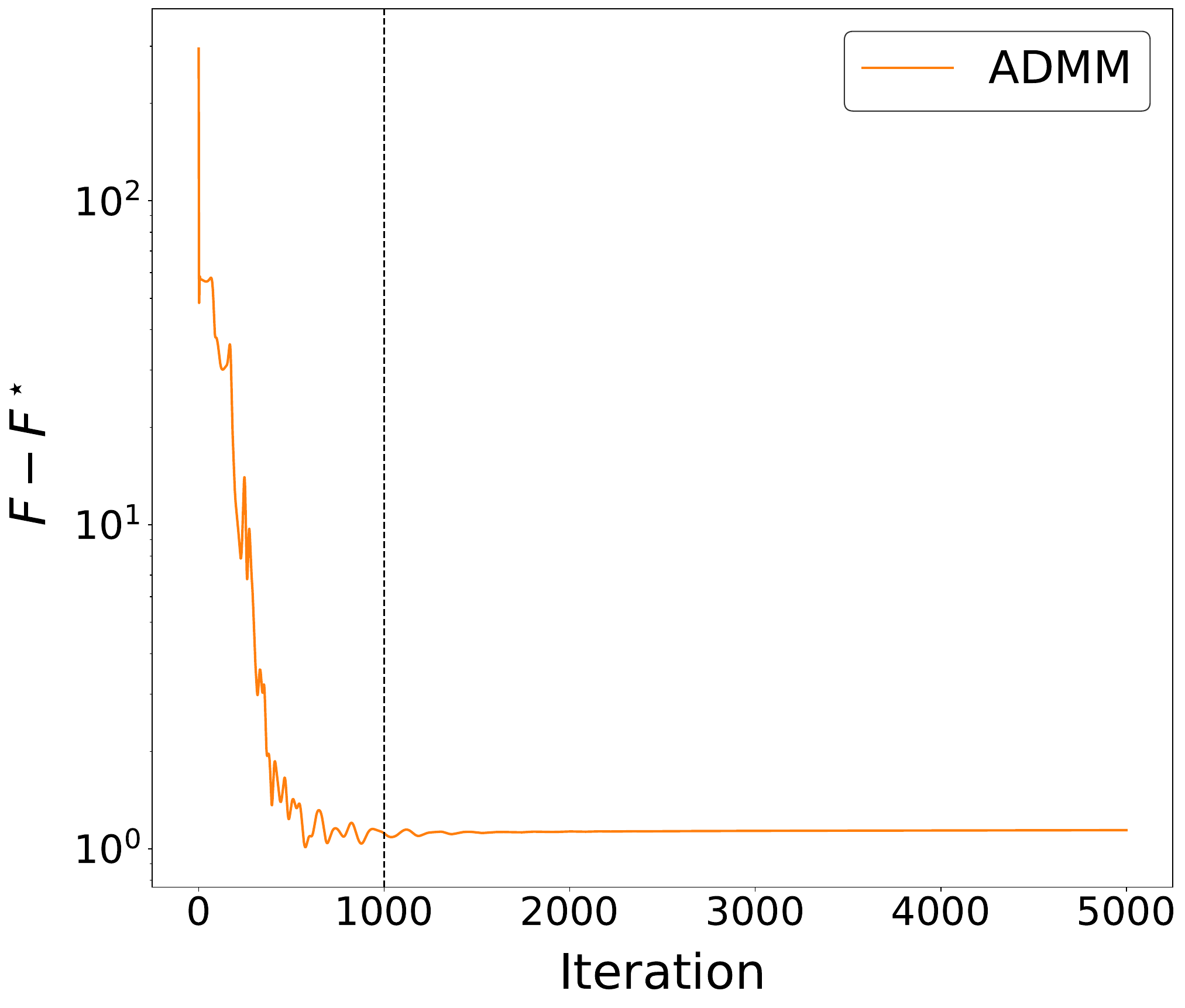}
\caption{\small Total error}
\label{subfig: noise-total-error}
\end{subfigure}
\begin{subfigure}[t]{0.32\linewidth}
\centering
\includegraphics[scale=0.15]{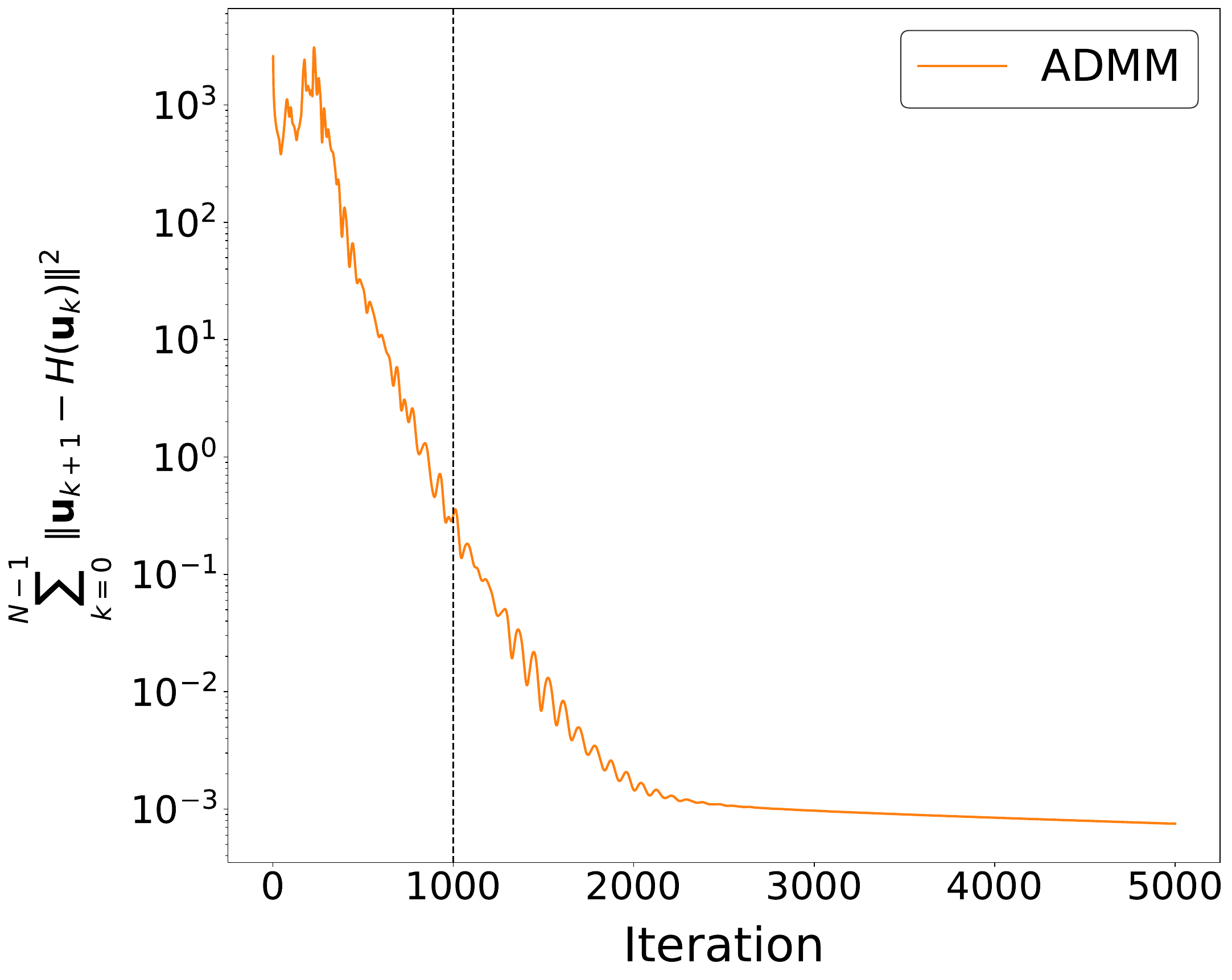}
\caption{\small Constraint error}
\label{subfig: noise-cons-error}
\end{subfigure}
\begin{subfigure}[t]{0.32\linewidth}
\centering
\includegraphics[scale=0.15]{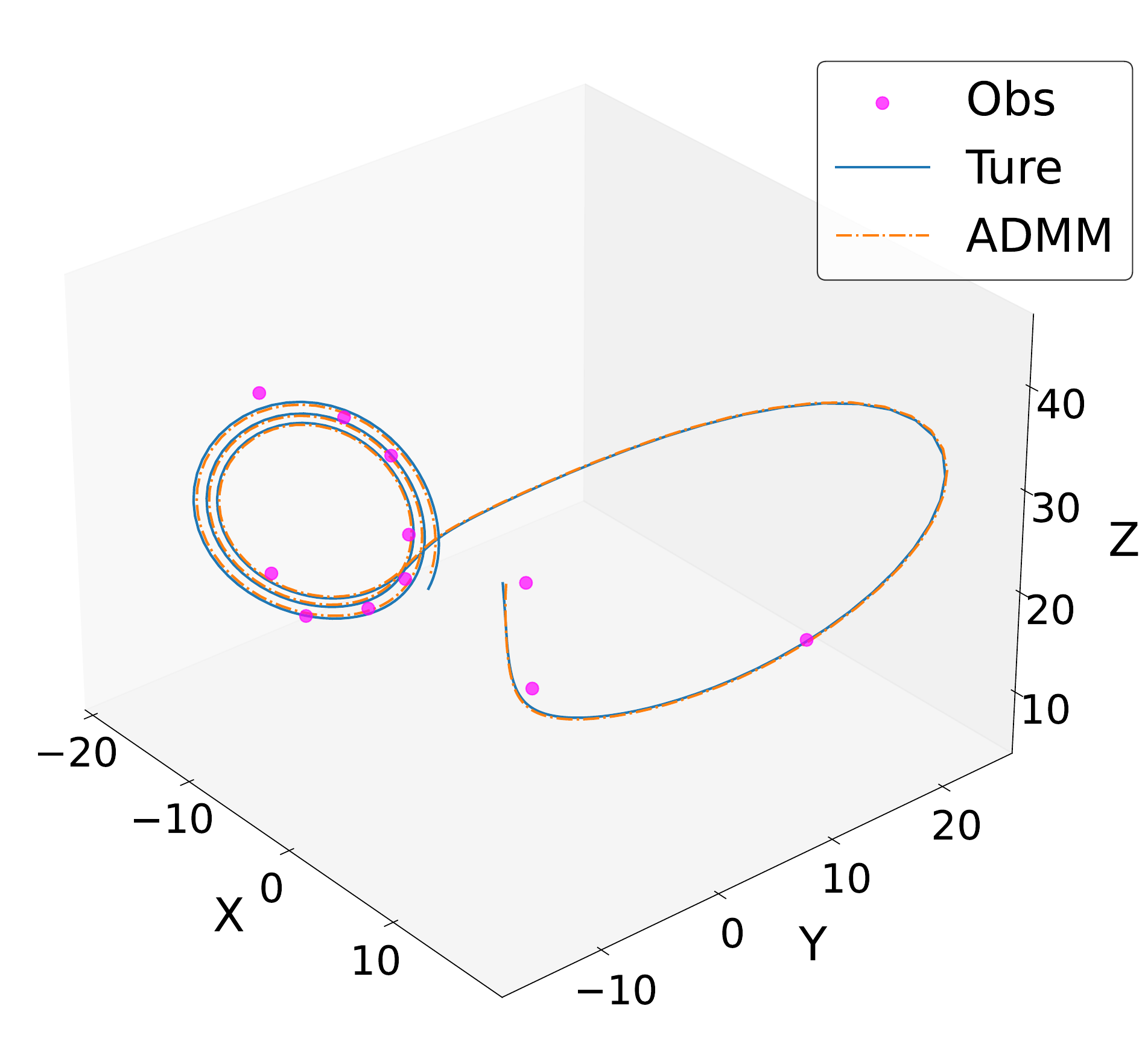}
\caption{\small Recovery solution ($\ell=1000$)}
\label{subfig: noise-recover-soln}
\end{subfigure}
\caption{\small Numerical performance of the linearized multi-block~\texttt{ADMM} with regularization~\eqref{eqn: lin-ADMM-reg}, which is applied to the~\texttt{4D-Var} problem~\eqref{eqn: numerical-4dvar-admm} for the Lorenz system~\eqref{eqn: lorenz} under noisy observation conditions, as specified in~\eqref{eqn: noise-lorenz}. } 
\label{fig: noise_lorenz_admm}
\end{figure}

\paragraph{Noisy observation}
The noisy observational data is generated by adding Gaussian noise to the precise observational data. Specifically, the noisy data is given by:
\begin{equation}
\label{eqn: noise-lorenz}
\left( \hat{x}_2(k), \hat{y}_2(k), \hat{z}_2(k) \right) = \left( \hat{x}_1(k) +\varepsilon _1, \hat{y}_1(k) + \varepsilon _2, \hat{z}_1(k) + \varepsilon _3 \right),
\end{equation}
where $\varepsilon_i \sim \mathcal{N}(0,1)$ for any $i=1,2,3$. Under the same settings as in the precise observation case, we apply the linearized multi-block~\texttt{ADMM} with regularization~\eqref{eqn: lin-ADMM-reg} to solve the~\texttt{4D-Var} problem~\eqref{eqn: numerical-4dvar}. The resulting numerical performance is depicted in~\Cref{fig: noise_lorenz_admm}.  
By comparing~\Cref{subfig: precise-total-error} (precise observation) and~\Cref{subfig: noise-total-error} (noisy observation), we observe a key difference: in the precise observation scenario, the total error converges consistently to zero, while in the noisy observation scenario, the total error converges to a non-zero value due to the influence of the noise. However, despite this noise-induced offset, the convergence remains stable, indicating robustness in the optimization process. Moreover, as illustrated in~\Cref{subfig: noise-cons-error}, we continue to monitor the constraint, $\sum_{k=0}^{N-1}\|\pmb{u}_{k+1} - H(\pmb{u}_k)\|^2$, which ensures that the solution recovered by the linearized multi-block~\texttt{ADMM} with regularization~\eqref{eqn: lin-ADMM-reg} remains accurate. Finally,~\Cref{subfig: noise-recover-soln} shows that, despite the presence of noise, the recovered solution remains close to the true numerical solution, demonstrating the method’s resilience in handling noisy data while maintaining a high degree of accuracy.

\subsection{Advantages of the linearized multi-block ADMM with regularization}
The algorithmic structure of the linearized multi-block~\texttt{ADMM} with regularization~\eqref{eqn: lin-ADMM-reg}, along with its numerical performance as demonstrated in~\Cref{fig: lorenz_admm} and~\Cref{fig: noise_lorenz_admm}, highlights several significant advantages for solving the~\texttt{4D-Var} problem~\eqref{eqn: numerical-4dvar}.  These advantages are outlined as follows:

%By comparing~\Cref{fig: lorenz_cg-bfgs} and~\Cref{fig: lorenz_admm}, 

\begin{itemize}
\item[(1)] \textbf{Effective Utilization of Observational Data}  In solving the~\texttt{4D-Var} problem with classical first-order optimization algorithms, the solution is assumed to strictly follow the governing dynamics, making it highly sensitive to the selection of initial conditions. This sensitivity, combined with noise in observational data, can significantly degrade performance. However, when using the linearized multi-block~\texttt{ADMM} with regularization~\eqref{eqn: lin-ADMM-reg}, the recovered solution capture the whole dynamics, not just the initial condition. At the outset, the iterative points do not need to strictly satisfy the constraints, or is a specific dynamical solution; they only need to converge toward the constraints. Particularly in the early stages, if the initial condition is far from the observational data, the $\argmin$ operation effectively pulls the solution closer to the data without being constrained to a specific dynamical solution. Additionally, a scaling parameter $\mu > 0$ can be introduced in the objective function to balance the solution between the observational data and constraints, improving overall performance.

%can be tuned to align each subproblem more closely with the observational data.

%\item[(2)] \textbf{Observational Data at All Time Points}\; Unlike classical optimization algorithms that focus on minimizing an objective function for initial conditions,  the linearized multi-block~\texttt{ADMM} with regularization~\eqref{eqn: lin-ADMM-reg} directly addresses the \texttt{4D-Var} problem~\eqref{eqn: numerical-4dvar} by deriving the Euler-Lagrange equation for all the intermediate states.  The quadratic structure of the objective function, with respect to the global dynamics, enables the use of a rescaling technique. This technique brings the augmented Lagrangian closer to convexity without changing the optimal solution,  thus facilitating the search for the global minimum.

\item[(2)] \textbf{Quadratic Subproblems}\; In the linearized multi-block~\texttt{ADMM} with regularization~\eqref{eqn: lin-ADMM-reg}, each subproblem is framed as a quadratic function, simplifying the optimization process.  By adjusting the regularization coefficient $\eta > 0$, the difference between consecutive iterations $\pmb{u}^{\ell+1}$ and $\pmb{u}^{\ell}$ becomes smaller, making the linearization a reasonable approximation. The flexibility in tuning the parameters $\mu > 0$ and $\eta>0$ allows us to improve the optimization performance. However, as the time step size $\delta t$ increases, the nonlinearity in the updates becomes more pronounced, which can degrade its performance.

% When the time step size $\delta t$ is sufficiently small, the iterative update $\pmb{u}_{k+1} = H(\pmb{u}_k)$ approximates a linear relationship.
%, thus simplifying s. 

\item[(3)] \textbf{Parallelizable Implementation}\; The linearized multi-block~\texttt{ADMM} with regularization~\eqref{eqn: lin-ADMM-reg} is well-suited for parallel implementation. This feature is is especially advantageous for problems involving long-term nonlinear evolutions, where computational demands can be substantial. With sufficient computing resources,  the method can be implemented effectively, making it scalable and practical for real-world applications.
\end{itemize}

Overall, these advantages highlight the efficiency and practical utility of the linearized multi-block~\texttt{ADMM} with regularization~\eqref{eqn: lin-ADMM-reg} in addressing the complex~\texttt{4D-Var} problem~\eqref{eqn: numerical-4dvar}.

%% file: 03_burgers.tex
\section{Viscous Burgers' equation}
\label{sec: burgers}

In this section, we apply the~\texttt{4D-Var} problem~\eqref{eqn: numerical-4dvar} to a classical nonlinear PDE --- the viscous Burgers' equation. The solution is assumed to satisfy $u \in L^{2}\left([0, T); H_0^1([0,\pi]) \right)$ and $u_t \in  L^{2}\left([0, T); H^{-1}([0,\pi]) \right)$. Under the Dirichlet boundary condition, the viscous Burgers' equation is given by: 
\begin{equation}
\label{eqn: burgers}
\left\{ \begin{aligned}
          & \partial_t u + u \cdot \partial_x u = \gamma \partial_{xx}u \\
          & u(0,x) = \sin x                                                             \\   
          & u(t,0) = u(t,\pi) =0,
          \end{aligned} \right.
\end{equation}
where the viscous coefficient is $\gamma = 0.05$. We demonstrate the numerical performance of the linearized multi-block~\texttt{ADMM} with regularization~\eqref{eqn: lin-ADMM-reg}, based on the~\texttt{4D-Var} reformulation~\eqref{eqn: numerical-4dvar-admm}, across various numerical discretization methods, including finite difference, (Galerkin) finite element, and spectral methods. For time discretization, we employ the forward Euler scheme.

%TLM
%\begin{equation}
%\label{eqn: burgers-tlm}
%\left\{ \begin{aligned}
%          & \partial_t \delta u +  u \cdot \partial_x \delta u + \partial_x u \cdot \delta u = \gamma \partial_{xx} \delta u \\
%          & \delta u(0,x) = \delta u_0(x)                                                                                                                                       \\   
%          & \delta u(t,0) = \delta u(t,\pi) =0
%          \end{aligned} \right.
%\end{equation}

%%%%%%%%%%%%%%%%%%%%%%%%%%%%%%%%%%%%%%%%%%%%%%%%%%%%%%%%%%%%%%%%%%%%%%%%%%%%%%%%%%%%%%%%%%%%%%%%%%%%%%%%%%%%%%%%%%%
\subsection{Finite difference method}
\label{subsec: fd-burgers}

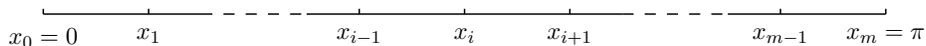
\begin{figure}[htb!]
	\centering
	\scalebox{1.4}{\begin{tikzpicture}
		\draw[-]  (0,0) -- (1.5,0);
                  \draw[dashed] (1.5,0) -- (2.5,0);
		\draw[-]  (2.5,0) -- (5.5,0);
                  \draw[dashed] (5.5,0) -- (6.5,0);
		\draw[-]  (6.5,0) -- (8,0);
%		\foreach \i in {0.0,0.1,...,5.0}
%		{\draw[very thin]
%		(\i,0)--(\i,0.15);}
		\foreach \I in {0,1,3,4,5,7,8}
		{\draw (\I,0)--(\I,0.05);}

		\node[below] at (0,0)  {$\scalemath{0.6}{x_{0} = 0}$};
		\node[below] at (1, 0) {$\scalemath{0.6}{x_{1}}$};
		\node[below] at (3, 0) {$\scalemath{0.6}{x_{i-1}}$};
		\node[below] at (4, 0) {$\scalemath{0.6}{x_{i}}$};
		\node[below] at (5, 0) {$\scalemath{0.6}{x_{i+1}}$};
		\node[below] at (7, 0) {$\scalemath{0.6}{x_{m-1}}$};
		\node[below] at (8, 0) {$\scalemath{0.6}{x_{m}=\pi}$};

		\end{tikzpicture}}
\caption{\small Uniform spatial discretization of the interval $[0,\pi]$.}
\label{fig: fd-label}
\end{figure}

Let the time step size be $\delta t = 0.02$ and the total evolution time $T = 2$.  The spatial grid size is set to $m=100$, resulting in a spatial discretization of $\delta x = \pi/m$. The spatial grid points $x_i$, for $i= 0, \ldots, m$, are illustrated in~\Cref{fig: fd-label}.  
For spatial discretization, we apply the central difference scheme. Given the Dirichlet boundary condition, $u(t,0) = u(t,\pi)=0$, we focus on the interior spatial grid points $x_i$ for $i=1,2,\ldots,m-1$. 

For $k=0,1,\ldots, N = T/\delta t = 100$, let $u_{k,i} = u_{k,i}^{(m)}$ represent the finite difference approximation of the analytic solution $u(k\delta t, i\delta x)$. For the viscous Burgers' equation~\eqref{eqn: burgers}, the central difference scheme is given by:
\begin{equation}
\label{eqn: fd-burgers}
\frac{u_{k+1,i} - u_{k,i}}{\delta t} + \frac{(u_{k,i+1})^2 - (u_{k,i-1})^2}{4\delta x} = \gamma \cdot \frac{u_{k,i+1} + u_{k,i-1} - 2u_{k,i}}{(\delta x)^2}.
\end{equation}
This scheme~\eqref{eqn: fd-burgers} can be reformulated into the iterative update form $\pmb{u}_{k+1} = H(\pmb{u}_k)$, as follows:
\begin{align}
u_{k+1,i} = & \left[ \frac{ \gamma \delta t}{(\delta x)^2}  \cdot u_{k,i-1} +  \frac{ \delta t}{4 \delta x}  \cdot (u_{k,i-1})^2  \right] \nonumber \\ 
                         & + \left( 1 -  \frac{2\gamma \delta t}{(\delta x)^2} \right) u_{k,i} +  \left[ \frac{ \gamma \delta t}{(\delta x)^2}  \cdot u_{k,i+1} -  \frac{ \delta t}{4 \delta x}  \cdot (u_{k,i+1})^2  \right],  \label{eqn: fd-burgers-numerical}
\end{align}
for $i=1,2,\ldots,m-1$.  Let the observational time be $T_o = 0.2$, resulting in $n+1$ observational time points, where $n = T/T_o = 10$. The numerical solution is generated using the central difference scheme~\eqref{eqn: fd-burgers}, with the initial condition $u_{0,i} = \sin\left( i\pi/m \right)$ for $i=1,2,\ldots,m-1$. The precise observational data is recorded at intervals of $T_o$, corresponding to the numerical solution $\pmb{u}_{10k}$ for $i=0,1,\ldots,n$, where each component is $u_{10k,i}$ for $i=1,2,\ldots,m-1$. Gaussian noise is then added to the precise observations to produce noisy observational data $\hat{\pmb{u}}_{k}$, given by:
\begin{equation}
\label{eqn: noise-fd}
\hat{u}_{k,i} = u_{10k,i} + 0.1\varepsilon _i
\end{equation}
where $\varepsilon _i \sim \mathcal{N}(0,1)$ for $i=1,2,\ldots,m-1$.

The tangent linear iteration corresponding to the iterative update~\eqref{eqn: fd-burgers-numerical} is given by:
\begin{align}
\delta u_{k+1,i} = & \left[ \frac{ \gamma \delta t}{(\delta x)^2} +  \frac{ \delta t}{2 \delta x}  \cdot u_{k,i-1}   \right] \delta u_{k,i-1}  \nonumber \\ 
                         & + \left( 1 -  \frac{2\gamma \delta t}{(\delta x)^2} \right) \delta u_{k,i} +  \left[ \frac{ \gamma \delta t}{(\delta x)^2}   -  \frac{ \delta t}{2 \delta x}  \cdot u_{k,i+1}  \right] \delta u_{k,i+1},  \label{eqn: fd-burgers-numerical-tlm}
\end{align}
for $i=1,2,\ldots,m-1$, which allows us to compute the Jacobian matrix, $\nabla H(\pmb{u}_k)$. Since $m=99$ is not large, its adjoint operator, $\nabla H(\pmb{u}_k)^{\top}$, is easily obtained. The linearized multi-block~\texttt{ADMM} with regularization~\eqref{eqn: lin-ADMM-reg} is then applied to solve the~\texttt{4D-Var} problem~\eqref{eqn: numerical-4dvar-admm}. The initial condition is set as $\pmb{u}_k^{0} = (0,0,\ldots,0)^{\top}$ for $k=0,1,\ldots,N$, with parameters $\mu = 20$, $\eta = 0.1$, and $s=2/3$.

\begin{figure}[htb!]
\centering
\begin{subfigure}[t]{0.45\linewidth}
\centering
\includegraphics[scale=0.20]{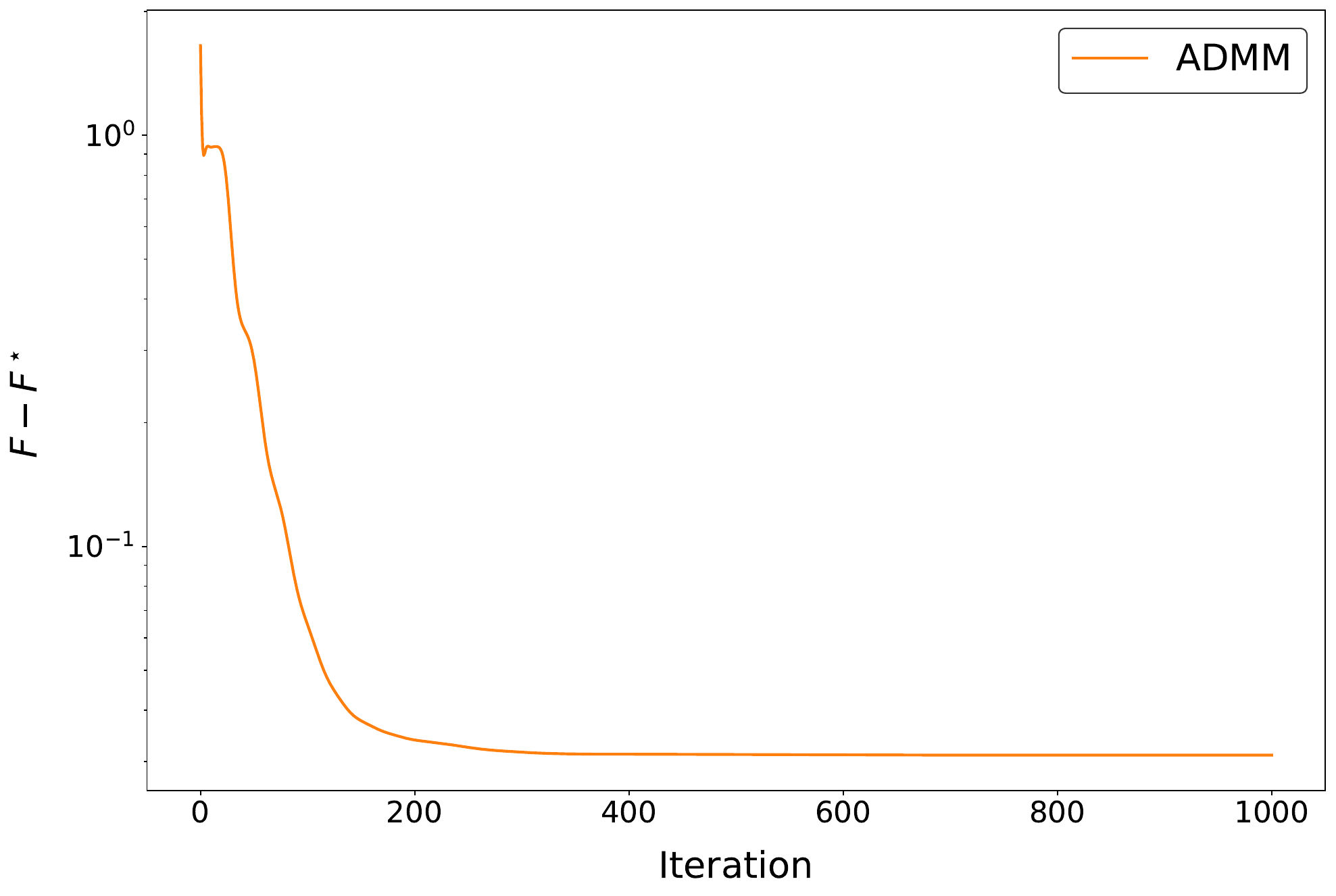}
\caption{Total error}
\end{subfigure}
\begin{subfigure}[t]{0.45\linewidth}
\centering
\includegraphics[scale=0.20]{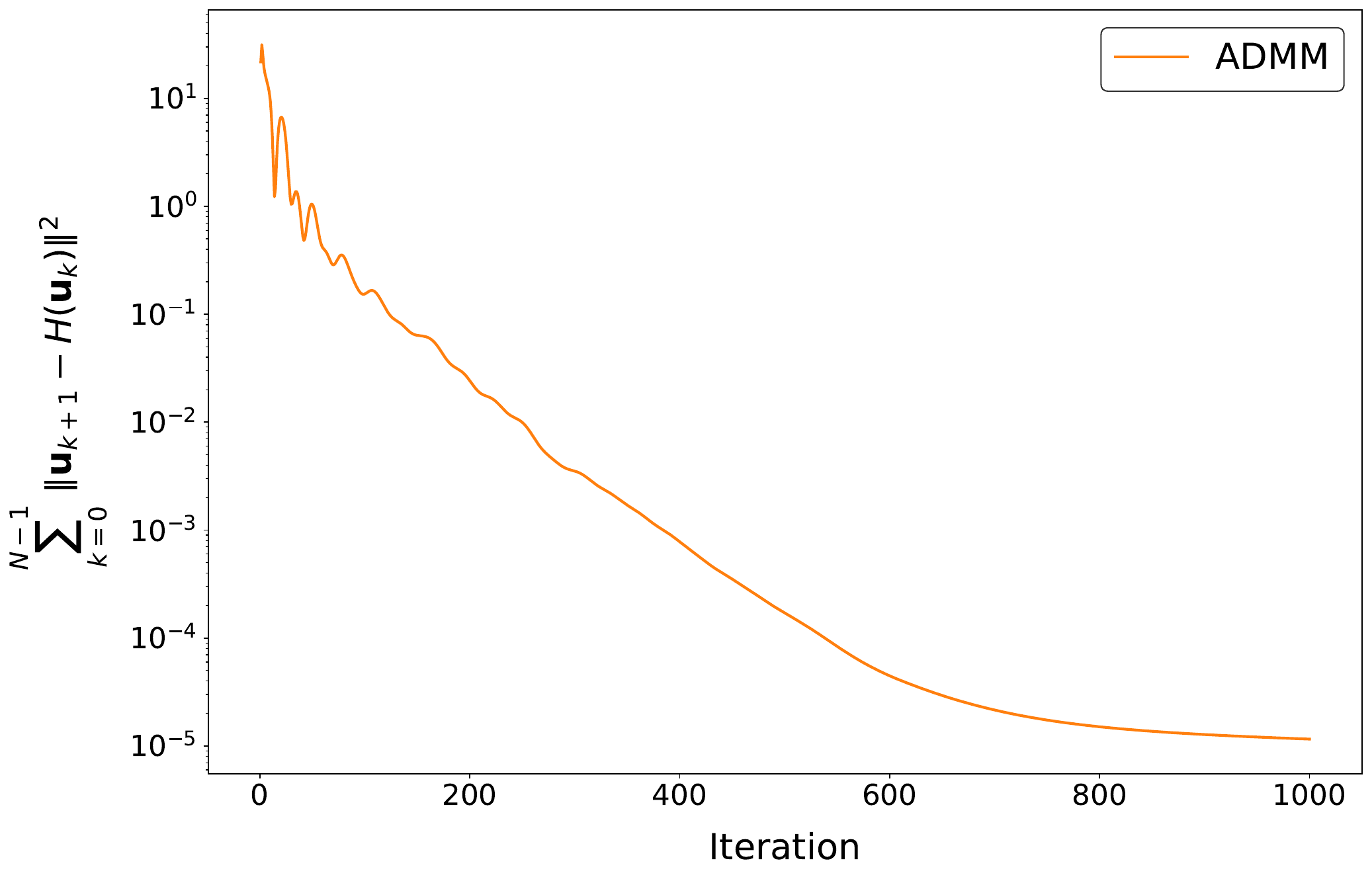}
\caption{Constraint error}
\end{subfigure}
\caption{\small The linearized multi-block~\texttt{ADMM} with regularization~\eqref{eqn: lin-ADMM-reg} is applied to  the~\texttt{4D-Var} problem~\eqref{eqn: numerical-4dvar-admm} for the viscous Burgers' equation~\eqref{eqn: burgers}, using the finite difference method~\eqref{eqn: fd-burgers}. } 
\label{fig: fd-convergence}
\end{figure}

\begin{figure}[htb!]
\centering
\includegraphics[scale=0.40]{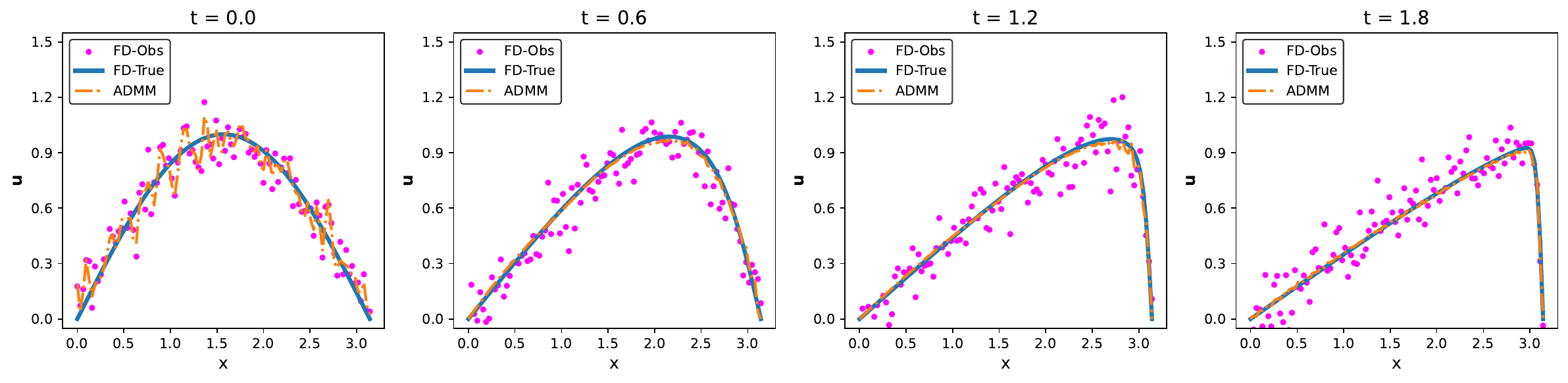}
\caption{\small Comparison of noisy observational data with the time evolution dynamics recovered via the linearized multi-block~\texttt{ADMM} with regularization~\eqref{eqn: lin-ADMM-reg}, with the true numerical solution as the reference. } 
\label{fig: fd-dynamic}
\end{figure}

\Cref{fig: fd-convergence} illustrates the numerical performance,
highlighting the convergence behavior of both the total error and the constraint error for the linearized multi-block~\texttt{ADMM} with regularization~\eqref{eqn: lin-ADMM-reg}, applied to the \texttt{4D-Var} problem~\eqref{eqn: numerical-4dvar-admm} for the viscous Burgers' equation using the finite difference method~\eqref{eqn: fd-burgers}. The total error stabilizes while the constraint error consistently converges,  exhibiting a similar pattern to that observed in the Lorenz system~(\Cref{fig: noise_lorenz_admm}). Additionally,~\Cref{fig: fd-dynamic} compares the dynamical evolution recovered via the linearized multi-block~\texttt{ADMM} with regularization~\eqref{eqn: lin-ADMM-reg} with the noisy observational data, using the true numerical solution as a reference. 
Despite the presence of noise, the recovered dynamics closely match the true numerical solution, with accuracy improving over time.

%of the objective function, $F-F^{\star}$, exhibits consistent convergence across iterations, as 
%
%
%%\begin{figure}[htb!]
%%\centering
%%\includegraphics[scale=0.30]{burgers/convergence_FD}
%%\caption{\small  The linearized multi-block~\texttt{ADMM} with regularization~\eqref{eqn: lin-ADMM-reg} is applied to  the~\texttt{4D-Var} problem~\eqref{eqn: numerical-4dvar-nonlinear} for the Burgers' equation with small viscosity, using the finite difference method~\eqref{eqn: fd-burgers}. } 
%%\label{fig: fd-convergence}
%%\end{figure}
%Furthermore, a detailed comparison between the observational data and the time evolution dynamics recovered via the linearized multi-block~\texttt{ADMM} with regularization~\eqref{eqn: lin-ADMM-reg} at iteration $\ell=200$ is presented in~\Cref{fig: fd-dynamic}. 
%
%
%This comparison highlights the effectiveness of the method in reconstructing the true dynamics alongside the observational data.

%%%%%%%%%%%%%%%%%%%%%%%%%%%%%%%%%%%%%%%%%%%%%%%%%%%%%%%%%%%%%%%%%%%%%%%%%%%%%%%%%%%%%%%%%%%%%%%%%%%%%%%%%%%%%%%%%%%
\subsection{(Galerkin) Finite element method}
\label{subsec: fe-burgers}

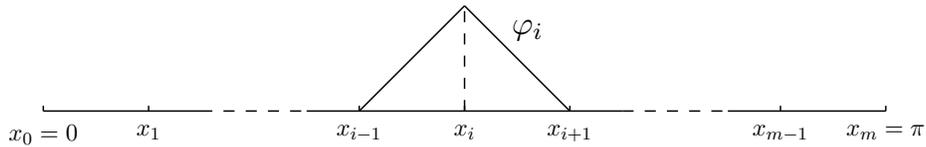
\begin{figure}[htb!]
	\centering
	\scalebox{1.4}{\begin{tikzpicture}
		\draw[-]  (0,0) -- (1.5,0);
                  \draw[dashed] (1.5,0) -- (2.5,0);
		\draw[-]  (2.5,0) -- (5.5,0);
                  \draw[dashed] (5.5,0) -- (6.5,0);
		\draw[-]  (6.5,0) -- (8,0);
                  \draw[-]  (3,0) -- (4,1);
                  \draw[-]  (4,1) -- (5,0);
                  \draw[dashed] (4,0) -- (4,1);   
		\foreach \I in {0,1,3,4,5,7,8}
		{\draw (\I,0)--(\I,0.05);}

		\node[below] at (0,0)  {$\scalemath{0.6}{x_{0} = 0}$};
		\node[below] at (1, 0) {$\scalemath{0.6}{x_{1}}$};
		\node[below] at (3, 0) {$\scalemath{0.6}{x_{i-1}}$};
		\node[below] at (4, 0) {$\scalemath{0.6}{x_{i}}$};
		\node[below] at (5, 0) {$\scalemath{0.6}{x_{i+1}}$};
		\node[below] at (7, 0) {$\scalemath{0.6}{x_{m-1}}$};
		\node[below] at (8, 0) {$\scalemath{0.6}{x_{m}=\pi}$};
                  \node[below] at (4.6,1.0)  {$\scalemath{0.8}{\varphi_{i}}$}; 

		\end{tikzpicture}}
\caption{\small Uniform spatial discretization of the interval $[0,\pi]$ and Lagrange nodal basis functions $\varphi_{i}$ for $i=1,2,\ldots,m-1$. }
\label{fig: fe-label}
\end{figure}

In the finite element method, the time step size is set to $\delta t = 0.01$ with a total evolution time of $T = 2$. The spatial grid points are identical to those used in the finite difference method, as shown in~\Cref{fig: fd-label}.  Since the solution $u=u(t,x)$ satisfies the Dirichlet boundary condition $u(t,0) = u(t,\pi)=0$, the Lagrange nodal basis functions are defined as:
\begin{equation}
\label{eqn: fe-basis}
\varphi_{i} = \left\{ \begin{aligned}
                                 & \frac{x  - x_{i-1}}{\delta x},      && x \in [x_{i-1}, x_{i}]    \\
                                 & \frac{x_{i+1} - x}{\delta x},     && x \in [x_{i}, x_{i+1}]  \\
                                 & 0,                                                      && \text{otherwise}  
                                 \end{aligned} \right.     
\end{equation}
for $i=1,2,\ldots,m-1$, as depicted in~\Cref{fig: fe-label}.  

For $k=0,1,\ldots, N=T/\delta t$ and $i=1,2,\ldots,m-1$, the numerical solution can be expressed as the following linear combination:
\[
u_{k}(x) = u_{k}^{(m)}(x) = \sum\limits_{i=1}^{m-1}u_{k,i}\varphi_i(x),
\] 
which represents the Galerkin approximation for the analytic solution $u(k\delta t,x)$. The vector of coefficients, $\pmb{u}_k =(u_{k,1}, u_{k,2}, \ldots, u_{k,m-1})^{\top}$, satisfies the following iterative scheme:
\begin{equation}
\label{eqn: fe-burgers}
\frac{\pmb{u}_{k+1} - \pmb{u}_{k}}{\delta t} + R^{-1}S_1[\pmb{u}_k] \pmb{u}_k = - \gamma R^{-1}T \pmb{u}_k,
\end{equation}
where the matrices $R$, $S_1[\pmb{u}_k]$, and $T$ are given as follows:
\[
  R = \delta x \begin{pmatrix} 
       \frac23 &\frac16    &                 &                 &                  \\
       \frac16 &\frac23    & \frac16     &                 &                  \\
                    & \ddots   & \ddots     & \ddots      &                  \\
                    &               &  \ddots    & \ddots      & \frac16    \\
                    &               &                 & \frac16     & \frac23  
       \end{pmatrix}, \quad T = \frac{1}{\delta x}\begin{pmatrix} 
       2         &-1          &                    &                  &              \\
       -1        &2           & -1                &                 &               \\
                   & \ddots & \ddots        & \ddots      &               \\
                   &             & \ddots        & \ddots      &       - 1   \\
                   &             &                    & -1            &         2  
       \end{pmatrix},  
\]
and       
\[
 S_1[\pmb{u}_k] = \begin{pmatrix} 
        \frac{u_k^2 - u_k^0}{3}      & \frac{u_k^2 - u_k^1}{6}       &                                          &                                                        &                                                      \\
        \frac{u_k^2 - u_k^1}{6}      & \frac{u_k^3 - u_k^1}{3}       & \frac{u_k^3 - u_k^2}{6}     &                                                        &                                                      \\
                                                  & \ddots                                 & \ddots                               & \ddots                                              &                                                      \\
                                                  &                                            & \ddots                               & \ddots                                              & \frac{u_k^{m-1} - u_k^{m-2}}{6}    \\
                                                  &                                            &                                          & \frac{u_k^{m} - u_k^{m-1}}{6}         & \frac{ u_k^{m}- u_k^{m-2}}{3}   
       \end{pmatrix}.
\]
The observational setup follows the same configuration described in~\Cref{subsec: fd-burgers}, with the observational time $T_o = 0.2$ and $n+1$ observational time points, where $n = T/T_o = 10$. The numerical solution is generated using the finite element scheme~\eqref{eqn: fe-burgers} with the initial condition $\pmb{u}_{0} =  R^{-1}T \pmb{\sin}[x]$, where $\pmb{\sin}[x] = (\sin x_1, \sin x_2, \ldots, \sin x_{m-1})^{\top}$. At each observational time point,  the numerical solution $\pmb{u}_{10k}$ is recorded for $i=0,1,\ldots,n$ with each component given by $u_{10k,i}$ for $i=1,2,\ldots,m-1$. Noisy observational data $\hat{\pmb{u}}_{k}$ is generated by adding Gaussian noise to the observations:
\begin{equation}
\label{eqn: noise-fe}
\hat{u}_{k,i} = u_{10k,i} + 0.1\varepsilon _i
\end{equation}
where $\varepsilon _i \sim \mathcal{N}(0,1)$ for $i=1,2,\ldots,m-1$. The tangent linear iteration corresponding to this iterative update~\eqref{eqn: fe-burgers} is:
\begin{equation}
\label{eqn: fem-burgers-numerical-tlm}
\delta \pmb{u}_{k+1} =\left[\mathbb{I} -  \delta t \left(R^{-1}S_1[\pmb{u}_k] + R^{-1}S_2[\pmb{u}_k] + \gamma R^{-1}T \right)\right] \delta \pmb{u}_{k},
\end{equation}
where the matrix $S_2[\pmb{u}_k]$ is given by:
\[
 S_2[\pmb{u}_k] =  \begin{pmatrix} 
       - \frac{u_k^2}{6}                                  & \frac{2u_k^1 +  u_k^2}{6}          &                                              &                                                         &                                                      \\
       - \frac{u_k^1 + 2u_k^2}{6}                  & \frac{u_k^1 - u_k^3}{6}              & \frac{2u_k^2 + u_k^3 }{6}     &                                                         &                                                      \\
                                                                  & \ddots                                        & \ddots                                   & \ddots                                               &                                                      \\
                                                                  &                                                    & \ddots                                   & \ddots                                              & \frac{2u_k^{m-1} + u_k^{m}}{6}    \\
                                                                  &                                                    &                                              & - \frac{u_k^{m-1} + 2u_k^{m}}{6}     & \frac{u_k^{m-2} - u_k^{m}}{3}   
       \end{pmatrix},
\]
which leads to the Jacobian matrix $\nabla H(\pmb{u}_k)$ and its adjoint $\nabla H(\pmb{u}_k)^{\top}$.

\begin{figure}[htb!]
\centering
\begin{subfigure}[t]{0.45\linewidth}
\centering
\includegraphics[scale=0.20]{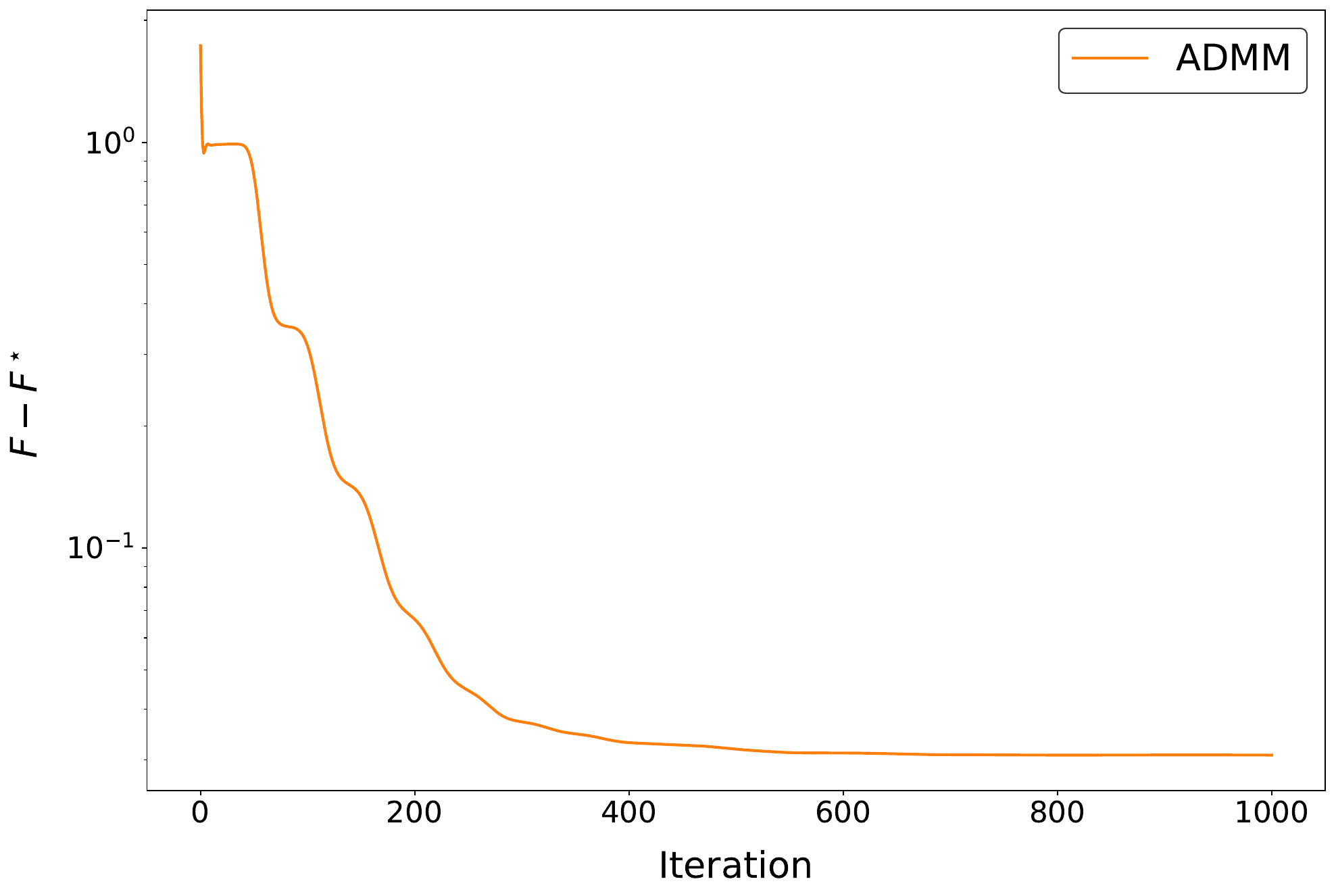}
\caption{\small Total error}
\end{subfigure}
\begin{subfigure}[t]{0.45\linewidth}
\centering
\includegraphics[scale=0.20]{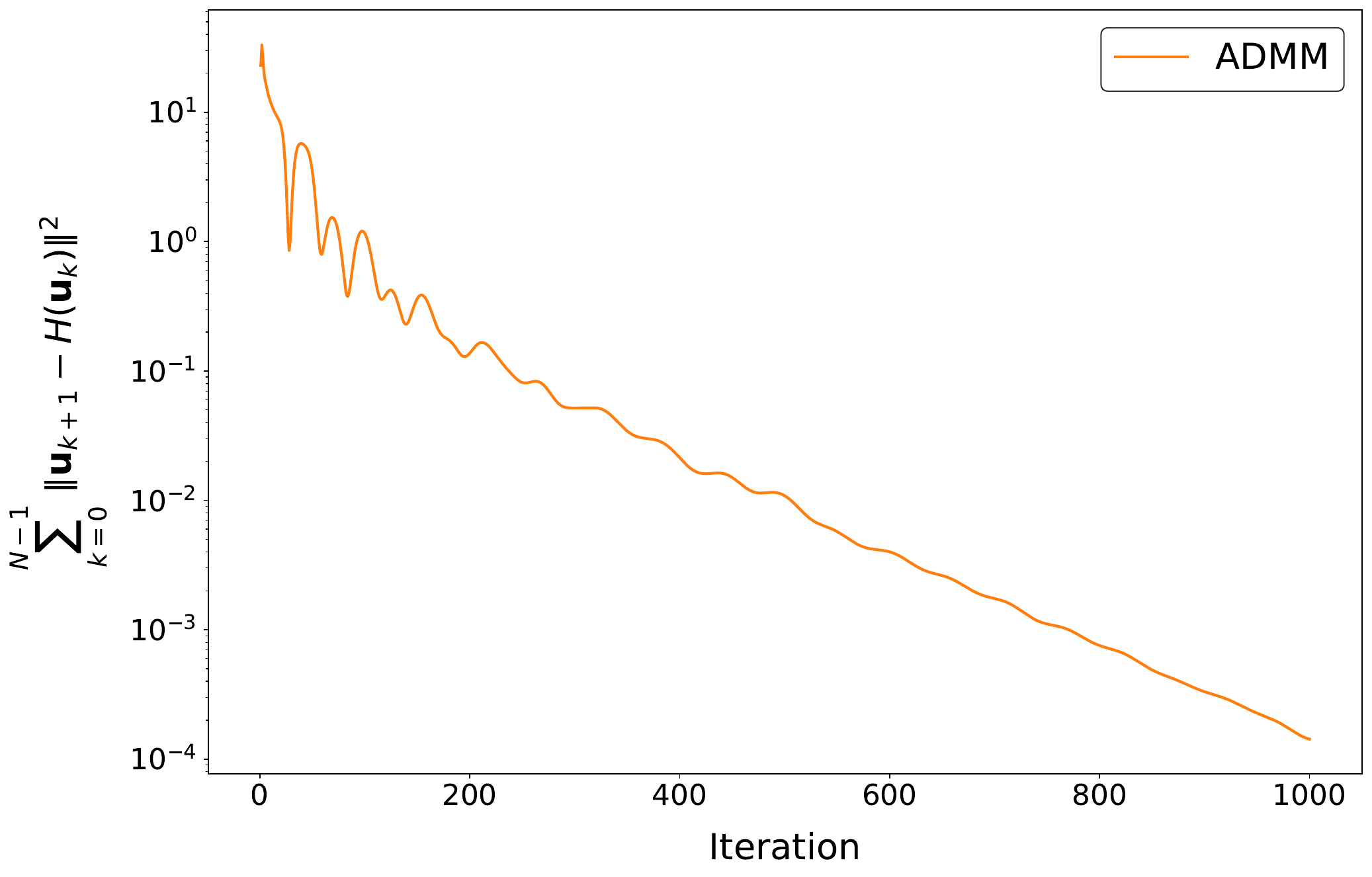}
\caption{\small Constraint error}
\end{subfigure}
\caption{\small The linearized multi-block~\texttt{ADMM} with regularization~\eqref{eqn: lin-ADMM-reg} is applied to  the~\texttt{4D-Var} problem~\eqref{eqn: numerical-4dvar-admm} for the viscous Burgers' equation~\eqref{eqn: burgers}y, using the finite element method~\eqref{eqn: fe-burgers} } 
\label{fig: fem-convergence}
\end{figure}

\begin{figure}[htb!]
\centering
\includegraphics[scale=0.40]{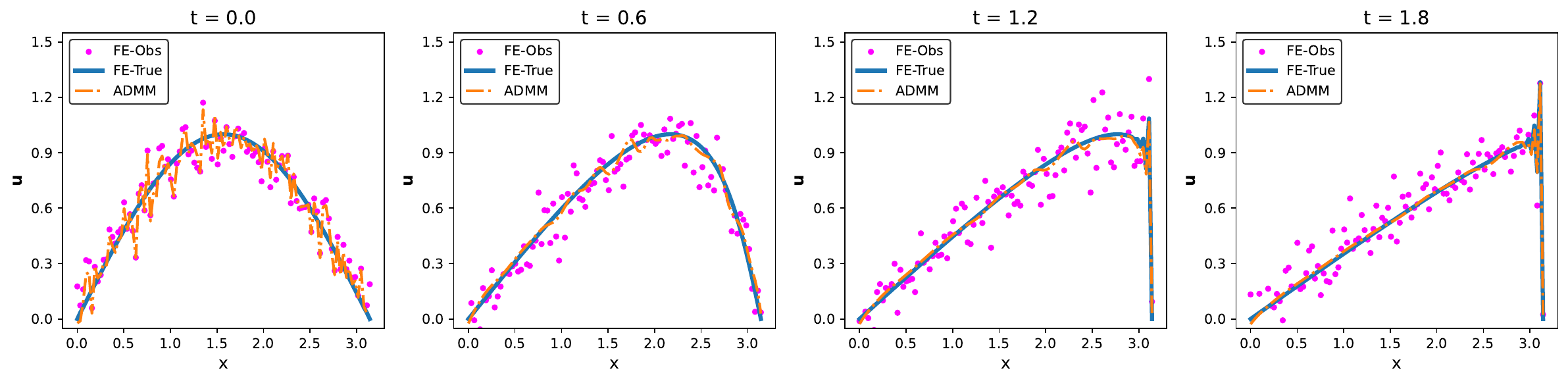}
\caption{\small Comparison of noisy observational data with the time evolution dynamics recovered via the linearized multi-block~\texttt{ADMM} with regularization~\eqref{eqn: lin-ADMM-reg}, with the true numerical solution as the reference.  } 
\label{fig: fe-dynamic}
\end{figure}

We then apply the linearized multi-block~\texttt{ADMM} with regularization~\eqref{eqn: lin-ADMM-reg} to solve the~\texttt{4D-Var} problem~\eqref{eqn: numerical-4dvar-admm}. The initial condition is set as $\pmb{u}_k^{0} = (0,0,\ldots,0)^{\top}$ for $k=0,1,\ldots,N$, with parameters $\mu =20$, $\eta = 0.1$, and $s=2/3$. The numerical performance, displayed n in~\Cref{fig: fem-convergence}, show the convergence behaviors similar to that observed in~\Cref{fig: fd-convergence}.
As the iterations increase, the total error stabilizes, while the constraint error consistently decreases.,~\Cref{fig: fe-dynamic} demonstrates the numerical performance of the dynamics recovered via the linearized multi-block~\texttt{ADMM} with regularization~\eqref{eqn: lin-ADMM-reg}.
The recovered dynamics closely track the true numerical solution, with only minor derivations due to noise in the observational data, and the accuracy improves over time. A key observation is that the finite difference method exhibits numerical dissipation without high-frequency oscillations, while the finite element method lacks numerical dissipation but contains high-frequency oscillations. The solution recovered via the multi-block~\texttt{ADMM} with regularization\eqref{eqn: lin-ADMM-reg} reflects the characteristics of the different numerical discretizations used to generate the observational data, highlighting this method's effectiveness in solving the~\texttt{4D-Var} problem~\eqref{eqn: numerical-4dvar} and its dependence on the observational data. 

\subsection{Spectral method}
\label{subsec: s-burgers}

For the spectral method, we also set the time step size $\delta t = 0.01$ and the total evolution time $T = 2$. Let $m$ be the total degrees of freedom. Given the Dirichlet boundary condition $u(t,0) = u(t,\pi)=0$, we select sine functions $\sin ix$ for $i=1,2,\ldots, m$ as the orthogonal basis.  For each time step $k=0,1,\ldots, N=T/\delta t$, the numerical solution is expressed as the following linear combination:
\[
u_{k}(x) = u_{k}^{(m)}(x) = \sum_{i=1}^{m} u_{k,i} \sin i x,
\]
which serves as the spectral approximation for the analytic solution $u(k\delta t,x)$.  The vector of coefficients $\pmb{u}_k =(u_{k,1}, u_{k,2}, \ldots, u_{k,m})^{\top}$, satisfies the following iterative scheme:
\begin{equation}
\label{eqn: spectral-burgers}
\left\{ \begin{aligned}
          & \frac{u_{k+1,1} - u_{k,1}}{\delta t} - \frac{1}{2}  \sum_{l = 1}^{m-1}u_{k,l}u_{k, l + 1} = - \gamma  u_{k,1}, \\
          & \frac{u_{k+1,i} - u_{k,i}}{\delta t} + \frac{i}{4} \left(\sum_{l =1}^{i}u_{k,l}u_{k,i-l} - 2\sum_{l = 1}^{m-i}u_{k,l}u_{k,i + l}  \right) = - \gamma i^2   u_{k,i} , \qquad \mathrm{for}\;i =2,3, \ldots,m.
          \end{aligned}\right.
\end{equation}
The observational setup remains the same, with observational time $T_o = 0.2$ and $n+1$ observational time points, where $n = T/T_o = 10$. The numerical solution is generated using the   spectral method~\eqref{eqn: spectral-burgers}, starting from the initial condition $\pmb{u}_{0} = (1,0, \ldots, 0)^{\top}$. At each time observational point,  the numerical solution $\pmb{u}_{10k}$ is recorded for $i=0,1,\ldots,n$ with each component denoted as $u_{10k,i}$ for $i=1,2,\ldots,m-1$. Noisy observational data $\hat{\pmb{u}}_{k}$ is generated by adding Gaussian noise to the observations:
\begin{equation}
\label{eqn: noise-sm}
\hat{u}_{k,i} = u_{10k,i} + 0.1\sqrt{2} \cdot 0.1\varepsilon _i
\end{equation}
where $\varepsilon _i \sim \mathcal{N}(0,1)$ for $i=1,2,\ldots,m-1$ and the scaling factor $0.1\sqrt{2}$ ensures consistency with previous cases.

\begin{figure}[htb!]
\centering
\begin{subfigure}[t]{0.45\linewidth}
\centering
\includegraphics[scale=0.20]{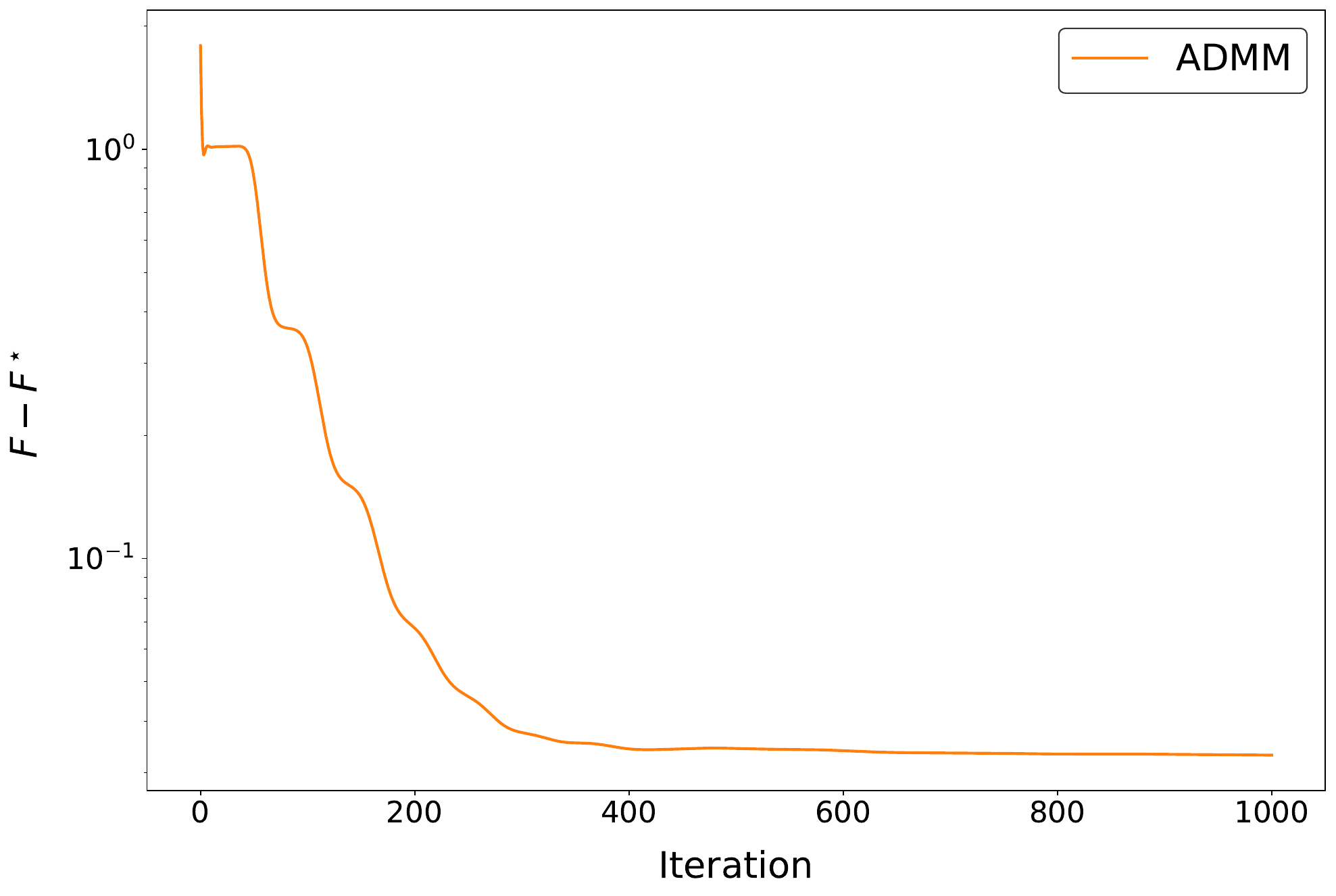}
\caption{\small Total error}
\end{subfigure}
\begin{subfigure}[t]{0.45\linewidth}
\centering
\includegraphics[scale=0.20]{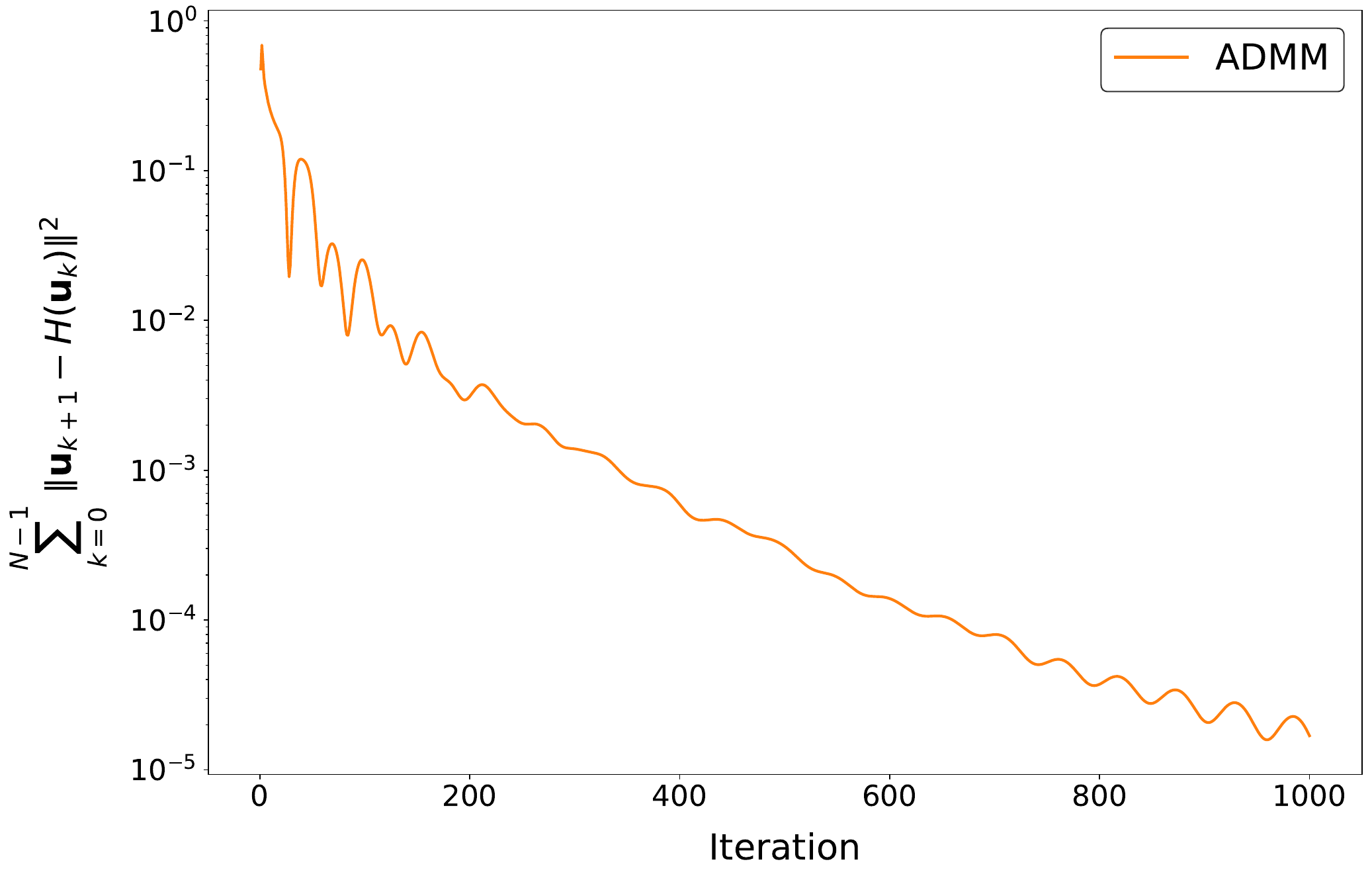}
\caption{\small Constraint error}
\end{subfigure}
\caption{\small  The linearized multi-block~\texttt{ADMM} with regularization~\eqref{eqn: lin-ADMM-reg} is applied to  the~\texttt{4D-Var} problem~\eqref{eqn: numerical-4dvar-nonlinear} for the viscous Burgers' equation~\eqref{eqn: burgers}, using the spectral method~\eqref{eqn: spectral-burgers}. } 
\label{fig: sm-convergence}
\end{figure}
\begin{figure}[htb!]
\centering
\includegraphics[scale=0.40]{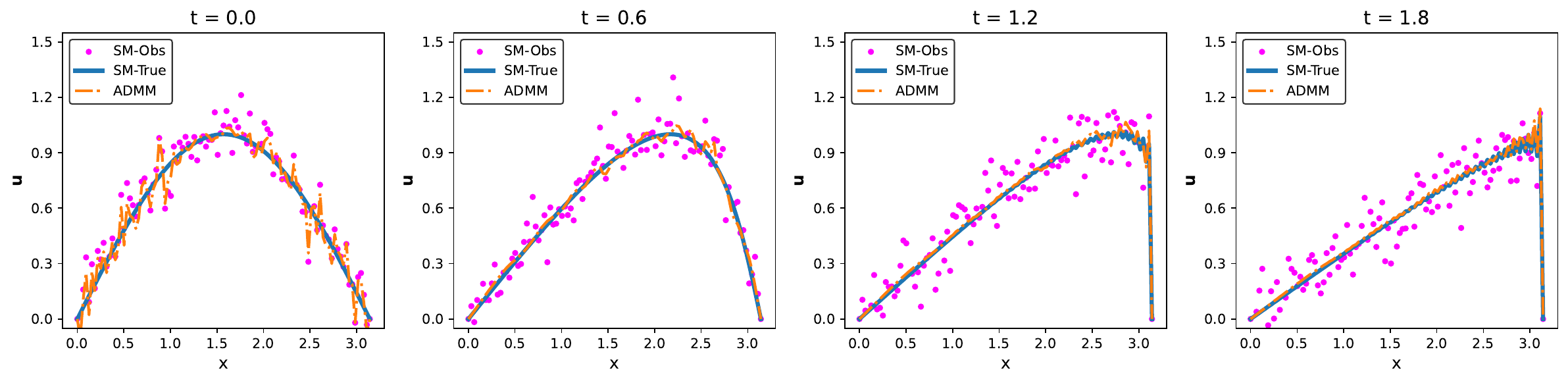}
\caption{\small Comparison of noisy observational data with the time evolution dynamics recovered via the linearized multi-block~\texttt{ADMM} with regularization~\eqref{eqn: lin-ADMM-reg}, with the true numerical solution as the reference.  } 
\label{fig: sm-dynamic}
\end{figure}

Before applying the linearized multi-block~\texttt{ADMM} with regularization~\eqref{eqn: lin-ADMM-reg} to solve the~\texttt{4D-Var} problem~\eqref{eqn: numerical-4dvar-admm},  we first derive the tangent linear iteration corresponding to the iterative scheme~\eqref{eqn: spectral-burgers} as follows: 
\begin{equation}
\label{eqn: spectral-burgers-tlm}
\left\{ \begin{aligned}
          & \delta u_{k+1,1} = \delta u_{k,1} + \delta t \left( \frac{1}{2}  \sum_{l = 1}^{m-1}\delta u_{k,l} \cdot u_{k, l + 1} + \frac{1}{2}  \sum_{l = 1}^{m-1}u_{k,l} \cdot \delta u_{k, l + 1} - \gamma  \delta u_{k,1}\right), \\
          & \delta u_{k+1,i} = \delta u_{k,i} - \delta t \gamma i^2 \cdot \delta u_{k,i}  \\
          & \;\;- \delta t  \cdot \frac{i}{4} \left(\sum_{l =1}^{i}\delta u_{k,l} \cdot u_{k,i-l} + \sum_{l =1}^{i} u_{k,l} \cdot \delta u_{k,i-l}- 2\sum_{l = 1}^{m-i}\delta u_{k,l} \cdot u_{k,i + l}  - 2\sum_{l = 1}^{m-i}u_{k,l} \cdot\delta u_{k,i + l}  \right), \\
                    & \mathrel{\phantom{abcdefghijklmnopqrstuvw}} \mathrm{for}\; i =2,3, \ldots,m,                                                           
        \end{aligned}\right.
\end{equation}
which helps us derive the Jacobian matrix $\nabla H(\pmb{u}_k)$ and its adjoint $\nabla H(\pmb{u}_k)^{\top}$. The iterative process begins with the initial condition as $\pmb{u}_k^{0} = (0,0,\ldots,0)^{\top}$ for $k=0,1,\ldots,N$, using the parameters $\mu = 20$, $\eta = 0.1$, and $s=2/3$. The numerical performance is illustrated in~\Cref{fig: sm-convergence}, showing consistency with the convergence behaviors observed in~\Cref{fig: fd-convergence} and~\Cref{fig: fem-convergence}.  Furthermore,~\Cref{fig: sm-dynamic} highlights the recovered dynamics, indicating that the solution recovered via the multi-block~\texttt{ADMM} with regularization\eqref{eqn: lin-ADMM-reg} effectively captures the characteristics of the different numerical discretizations used to generate the observational data. 

Finally, a comparison of the spectral method~\eqref{eqn: spectral-burgers} with the finite difference method~\eqref{eqn: fd-burgers} and the finite element method~\eqref{eqn: fe-burgers} reveals key distinctions in their computational complexity. The spectral method includes $m^2-1$ quadratic terms in the spectral method, while the finite difference method has $2m-2$, and the finite element method involves $6m-10$ quadratic terms. This difference leads the tangent linear iteration~\eqref{eqn: spectral-burgers-tlm} to introduce nonlinear terms that scale as $O(m^2)$ for the spectral method.
\begin{table}[htb!]
\centering
\begin{tabular}{lcc}
\toprule 
Numerical Methods                           &  $\phantom{1111111111}$      & Computational Time\\ 
\midrule 
Finite Difference               &                                               & $\phantom{11}65.5021$s                           \\
Finite Element                  &                                               & $\phantom{1}232.0715$s                           \\
Spectral Method               &                                               & $5113.0345$s                                            \\
\bottomrule 
\end{tabular}
\caption{\small Comparison of computational times for the linearized multi-block~\texttt{ADMM} with regularization~\eqref{eqn: lin-ADMM-reg} across the finite difference~\eqref{eqn: fd-burgers}, finite element~\eqref{eqn: fe-burgers}, and spectral~\eqref{eqn: spectral-burgers} discretizations. Python 3.11.5 and Numpy 1.25.2 were executed on an Intel\textsuperscript{\tiny\textregistered} Core\textsuperscript{\tiny TM} i5-12500, 3.00 GHz Processor (12th Generation). }
\label{fig: run-time}
\end{table}
As a result, the computational cost, particularly for calculating the adjoint operator, increases significantly. This higher complexity makes the spectral method substantially slower compared to the other methods, as evidenced by the computational times presented in~\Cref{fig: run-time}.
%compared to ~\eqref{eqn: spectral-burgers-tlm}, In other words, as the degree of freedom $m$ increases, the spectral discretization leads to the nonlinear term increasing as $O(m^2)$.  This leads to the computational times of the transpose of the Jacobian, $\nabla H(\pmb{u}_k)^{\top}$, from $O(m)$ to $O(m^2)$. 

%% file: 04_vorticity.tex
\section{Vorticity concentration in large-scale 2D turbulence }
\label{sec: inverse-cascade}

In this section, we apply the linearized multi-block~\texttt{ADMM} with regularization~\eqref{eqn: lin-ADMM-reg} to solve the~\texttt{4D-Var} problem~\eqref{eqn: numerical-4dvar-admm}  in the context of the large-scale 2D turbulence, with a particular focus on the phenomenon of vorticity concentration. Let $\Omega$ be an open domain, and the vorticity field is assumed to satisfy $\omega \in L^{2}\left([0, T); H_0^1(\Omega) \right)$ and $\omega_t \in  L^{2}\left([0, T); H^{-1}(\Omega) \right)$. The governing 2D vorticity equation is given by:
\begin{equation}
\label{eqn: 2d-qg-f}
\left\{ \begin{aligned}
          & \partial_t \omega + J(\psi, \omega) = -\kappa \Delta^2 \omega,                                            \\
          & \omega =  \Delta \psi ,                                                                                                                                                            
          \end{aligned} \right.
\end{equation}
where the Jacobian $J$, for any bivariate functions $u$ and $v$, is defined as:
\begin{equation}
\label{eqn: jacobian}
J(u,v) = \partial_xu \partial_yv - \partial_yu \partial_xv,
\end{equation}
and the biharmonic operator $\Delta^2$ serves as a scale-selective dissipation, filtering out high-frequency waves, as described in~\citep{mcwilliams1984emergence}. We impose a slip boundary condition:  
\begin{equation}
\label{eqn: omega-dirichlet}
\omega \big |_{\partial \Omega} = 0,
\end{equation}
which indicates the absence of tangential stress along the boundary for large-scale flow. Small-scale processes occurring near the boundary act as a buffer, allowing the large-scale flow to slip smoothly along the boundary. This condition is commonly used in large-scale ocean circulation models~\citep[Section 2.4]{pedlosky1996ocean}.

%\[
%J_1^{i,j} = \frac{( u^{i+1,j} - u^{i-1,j} )( v^{i,j+1} - v^{i,j-1} ) - ( u^{i,j+1} - u^{i,j-1} )( v^{i+1,j} - v^{i-1,j} )   }{4\delta x\delta y}  
%\]
%%%%%%%%%%%%%%%%%%%%%%%%%%%%%%%%%%%%%%%%%%%%%%%%%%%%%%%%%%%%%%%%%%%%%%%%%%%%%%%%%%%%%%%%%%%%%%%%%%%%%%%%%%%%%%%%%%
\subsection{Computational Implementation}
\label{subsec: computational-implementation}
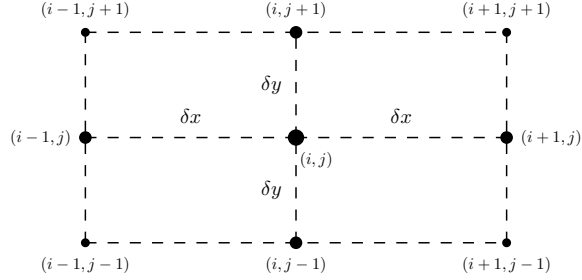
\begin{figure}[htb!]
	\centering
	\scalebox{1.4}{\begin{tikzpicture}
              
                 \filldraw [black] (0,0) circle (2pt);
                 \filldraw [black] (-2,0) circle (1.5pt);
                 \filldraw [black] (2,0) circle (1.5pt);
                 \filldraw [black] (0,1) circle (1.5pt);
                 \filldraw [black] (-2,1) circle (1pt);
                 \filldraw [black] (2,1) circle (1pt);
                 \filldraw [black] (0,-1) circle (1.5pt);
                 \filldraw [black] (-2,-1) circle (1pt);
                 \filldraw [black] (2,-1) circle (1pt);
		\node[below right, xshift=-0.1cm ] at (0,0)  {$\scalemath{0.4}{(i,j)}$};
                 \node[left] at (-2,0)  {$\scalemath{0.4}{(i-1,j)}$};
                 \node[right] at (2,0)  {$\scalemath{0.4}{(i+1,j)}$};
		\node[above] at (0,1)  {$\scalemath{0.4}{(i,j+1)}$};
                 \node[above] at (-2,1)  {$\scalemath{0.4}{(i-1,j+1)}$};
                 \node[above] at (2,1)  {$\scalemath{0.4}{(i+1,j+1)}$};
                  \node[below] at (0,-1)  {$\scalemath{0.4}{(i,j-1)}$};
                 \node[below] at (-2,-1)  {$\scalemath{0.4}{(i-1,j-1)}$};
                 \node[below] at (2,-1)  {$\scalemath{0.4}{(i+1,j-1)}$};
                 
                 \draw[dashed]  (-2,0) -- node[above]{$\scalemath{0.5}{\delta x}$} (0,0);
                 \draw[dashed]  (0,0) -- node[above]{$\scalemath{0.5}{\delta x}$} (2,0) ;      
                 \draw[dashed]  (-2,-1) --  (2,-1);
                 \draw[dashed]  (-2,1) -- (2,1) ;     
                                  
                 \draw[dashed]  (0,-1) -- node[left]{$\scalemath{0.5}{\delta y}$} (0,0);
                 \draw[dashed]  (0,0) -- node[left]{$\scalemath{0.5}{\delta y}$} (0,1) ;      
                 \draw[dashed]  (-2,-1) --  (-2,1);
                 \draw[dashed]  (2,-1) -- (2,1) ;                      
		\end{tikzpicture}}
\caption{\small Spatial grids for numerical Arakawa Jacobian $J(u,v)$ around the central point labeled $(i,j)$.}
\label{fig: arakawa}
\end{figure}

For the numerical implementation, we utilize a finite difference scheme, applying the Arakawa Jacobian, which is specifically designed to conserve energy and enstrophy in 2D fluid dynamics. In particular, we adopt the second-order Arakawa Jacobian as outlined in~\citep[Section 16.7]{cushman2011introduction}, which consists of three components:
\[
\left\{ \begin{aligned}
J_1^{i,j} = &\frac{( u^{i+1,j} - u^{i-1,j} )( v^{i,j+1} - v^{i,j-1} ) - ( u^{i,j+1} - u^{i,j-1} )( v^{i+1,j} - v^{i-1,j} )   }{4\delta x\delta y}, \\
J_2^{i,j} = & \frac{ u^{i+1,j}(v^{i+1,j+1} - v^{i+1,j-1}) - u^{i-1,j}(v^{i-1,j+1} - v^{i-1,j-1}) }{4\delta x\delta y}    \\
                 & - \frac{ u^{i,j+1}(v^{i+1,j+1} - v^{i-1,j+1}) - u^{i,j-1}(v^{i+1,j-1} - v^{i-1,j-1})}{4\delta x\delta y},  \\
J_3^{i,j} = & \frac{(u^{i+1,j+1} - u^{i-1,j+1}) v^{i,j+1} - (u^{i+1,j-1} - u^{i-1,j-1}) v^{i,j-1}}{4\delta x\delta y} \\
          & - \frac{(u^{i+1,j+1} - u^{i+1,j-1}) v^{i+1,j} - (u^{i-1,j+1} - u^{i-1,j-1}) v^{i-1,j}}{4\delta x\delta y}.                 
        \end{aligned} \right.
\]
The total Arakawa Jacobian is then computed as the average of these three components:
\begin{equation}
\label{eqn: arakawa-jacobian}
J^{i,j} = \frac{J_1^{i,j} + J_2^{i,j}  + J_3^{i,j} }{3}. 
\end{equation}
\Cref{fig: arakawa} illustrates the spatial grid around the central grid point $(i,j)$, highlighting the neighboring points used in the computation of the Jacobian. 
Numerically, both $\pmb{u}$ and $\pmb{v}$ are treated as vectors, when substituting a variable into the Arakawa Jacobian~\eqref{eqn: arakawa-jacobian}, it acts as a linear operator on the other variable, expressed as 
\begin{equation}
\label{eqn: jacobian-numerical}
J(\pmb{u},\pmb{v})= (J^{i,j})_{m^2 \times m^2} =\mathbb{J}[\pmb{u}]\pmb{v} = -\mathbb{J}[\pmb{v}]\pmb{u}.
\end{equation} 
To solve the inverse Laplacian, we utilize the successive over-relaxation (\texttt{SOR}) method~\citep{golub2013matrix}, though the conjugate gradient method~\citep{hestenes1952methods} is also a widely used approach. For time discretization, we employ a prediction-correction scheme:
\begin{subequations}
\label{eqn: pre-corr}
\begin{empheq}[left=\empheqlbrace]{align}
& \frac{\pmb{\omega}_k^p - \pmb{\omega}_k}{\delta t}   + J(\pmb{\psi}_k, \pmb{\omega}_k) = - \kappa \Delta^2 \pmb{\omega}_k,            \label{eqn: prediction} \\
& \frac{\pmb{\omega}_{k+1} - \pmb{\omega}_k}{\delta t} + J(\pmb{\psi}_k, \pmb{\omega}_k^p) = - \kappa \Delta^2 \pmb{\omega}_k^p.   \label{eqn: correction} 
\end{empheq}    
\end{subequations}
which ensures both accuracy and stability. The tangent linear iterations for the prediction-correction scheme~\eqref{eqn: pre-corr} are given by:
\begin{subequations}
\label{eqn: pre-corr-tlm}
\begin{empheq}[left=\empheqlbrace]{align}
& \frac{\delta \pmb{\omega}^p_{k} - \delta \pmb{\omega}_k}{\delta t} - \mathbb{J}[\pmb{\omega}_k] \Delta^{-1} \delta \pmb{\omega}_k + \mathbb{J}[\Delta^{-1}\pmb{\omega}_k] \delta \pmb{\omega}_k = -\kappa \Delta^2 \delta \pmb{\omega}_k,                        \label{eqn: prediction-tlm}            \\
& \frac{\delta \pmb{\omega}_{k+1} - \delta \pmb{\omega}_k}{\delta t} - \mathbb{J}[\pmb{\omega}^p_k] \Delta^{-1} \delta \pmb{\omega}_k + \mathbb{J}[\Delta^{-1}\pmb{\omega}_k] \delta \pmb{\omega}_k^p = -\kappa \Delta^2 \delta \pmb{\omega}_k^p.                  \label{eqn: correction-tlm} 
\end{empheq}    
\end{subequations}
Reformulating~\eqref{eqn: prediction} gives:
\begin{equation}
\label{eqn: prediction-reform}
\pmb{\omega}_k^p = \mathbb{P} \pmb{\omega}_k = \left[ \mathbb{I} - \delta t \left(\mathbb{J}[\Delta^{-1} \pmb{\omega}_k] - \kappa \Delta^2\right)\right] \pmb{\omega}_k,
\end{equation}
and similarly, reformulating~\eqref{eqn: prediction-tlm} gives:
\begin{equation}
\label{eqn: prediction-tlm-reform}
\delta \pmb{\omega}^p_{k} = \mathbb{\delta P}\delta \pmb{\omega}_k = \left[ \mathbb{I} + \delta t \left(  \mathbb{J}[\pmb{\omega}_k] \Delta^{-1} -  \mathbb{J}[\Delta^{-1}\pmb{\omega}_k] -\kappa \Delta^2 \right) \right] \delta \pmb{\omega}_k.
\end{equation}
Substituting~\eqref{eqn: prediction-reform} and~\eqref{eqn: prediction-tlm-reform} into~\eqref{eqn: correction-tlm} allows us to derive the tangent linear iteration: 
\begin{equation}
\label{eqn: jacobian-arakawa}
\delta \pmb{\omega}_{k+1} = \left[ \mathbb{I} + \delta t \left( \mathbb{J}[\mathbb{P}\pmb{\omega}_k]\Delta^{-1} - ( \mathbb{J}[\Delta^{-1}\pmb{\omega_k}] + \kappa \Delta^2) \mathbb{\delta P} \right) \right] \delta \pmb{\omega}_k
\end{equation}
and the adjoint operator of $\nabla H(\pmb{\omega_k})$ is given by:
\begin{equation}
\label{eqn: jacobian-arakawa-adjoint}
\nabla H(\pmb{\omega_k})^{\top} =  \left[ \mathbb{I} + \delta t \left( \Delta^{-1} \mathbb{J}[\mathbb{P}\pmb{\omega}_k]^{\top} - \mathbb{\delta P}^{\top} ( \mathbb{J}[\Delta^{-1}\pmb{\omega_k}]^{\top} + \kappa \Delta^2)  \right) \right], 
\end{equation}
which highlights the complexity of computing the adjoint operator in this context. Before implementing the linearized multi-block~\texttt{ADMM} with regularization~\eqref{eqn: lin-ADMM-reg} to solve the~\texttt{4D-Var} problem~\eqref{eqn: numerical-4dvar-admm}, we must establish the energy norm for the objective function as:
\begin{align}
 F(\pmb{\omega}_0, \pmb{\omega}_1, \ldots, \pmb{\omega}_n) & =  \frac{T_o}{2} \sum_{k=0}^{n}\big\| \nabla^{\top}(\pmb{\psi}_k - \hat{\pmb{\psi}}_k) \big\|^2 + \frac{\alpha}{2} \big\| \nabla^{\top}(\pmb{\psi}_0 - \hat{\pmb{\psi}}_0^b) \big\|^2   \nonumber            \\
                                                                                                       & = \frac{T_o}{2} \sum_{k=0}^{n}\big\| \Delta^{-\frac12}(\pmb{\omega}_k - \hat{\pmb{\omega}}_k) \big\|^2 + \frac{\alpha}{2} \big\| \Delta^{-\frac12}(\pmb{\omega}_0 - \hat{\pmb{\omega}}_0^b) \big\|^2, \label{eqn: energy-norm-vorticity}
\end{align}
where the second equation results from integration by parts.

%%%%%%%%%%%%%%%%%%%%%%%%%%%%%%%%%%%%%%%%%%%%%%%%%%%%%%%%%%%%%%%%%%%%%%%%%%%%%%%%%%%%%%%%%%%%%%%%%%%%%%%%%%%%%%%%%%
\subsection{Numerical performance}
\label{subsec: numerical-performance}

\begin{figure}[htb!]
\centering
\includegraphics[scale=0.26]{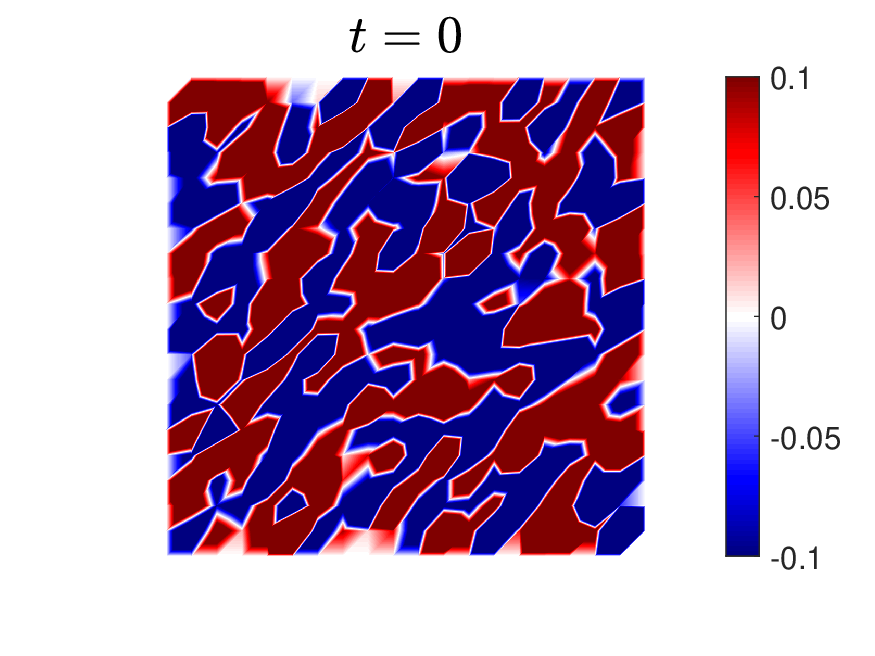}
\includegraphics[scale=0.26]{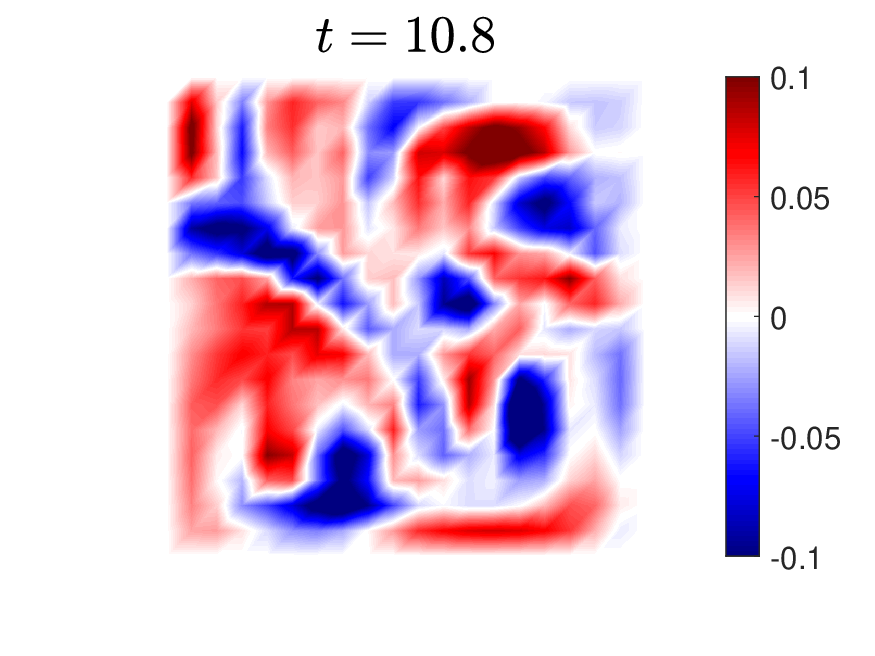}
\includegraphics[scale=0.26]{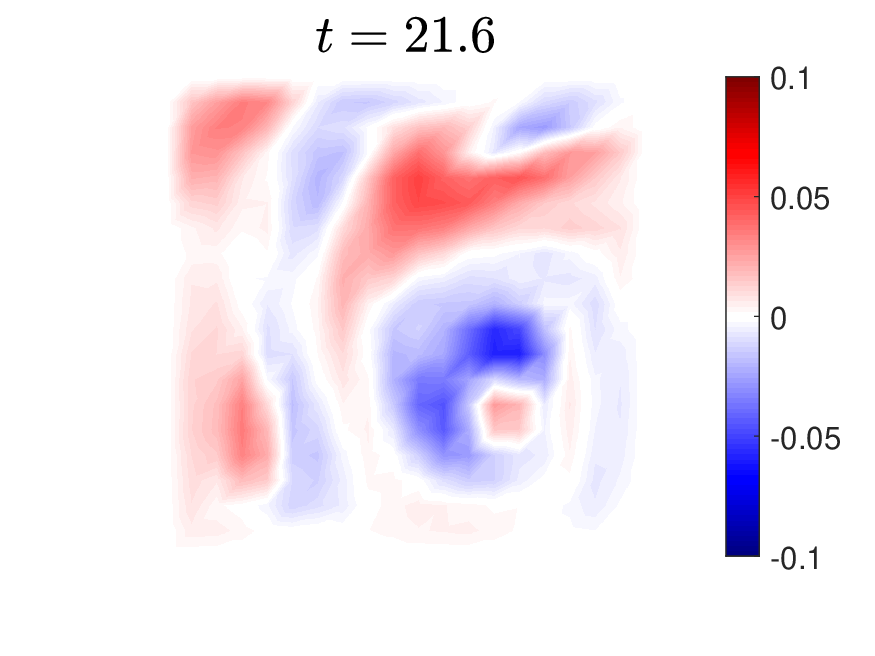}
\includegraphics[scale=0.26]{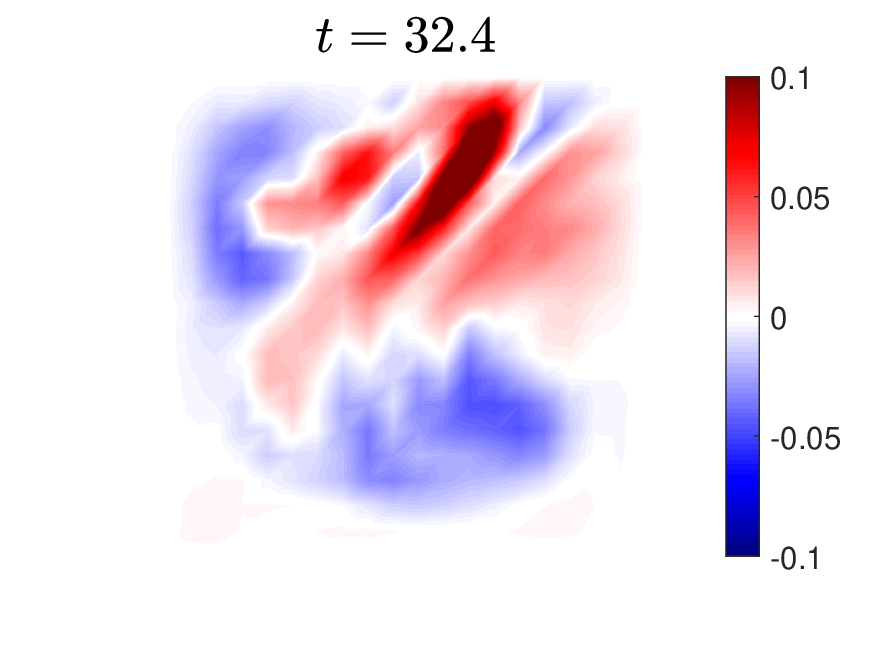}
\caption{\small The error between the solution recovered via the linearized multi-block~\texttt{ADMM} with regularization~\eqref{eqn: lin-ADMM-reg} and the ture numerical solution.} 
\label{fig: vorticity-error}
\end{figure}

We demonstrate the numerical performance within a box domain defined as $\Omega = [-2L, 2L] \times [-2L, 2L]$ where $L=1$. The spatial grids are configured with $\delta x = \delta y = 0.2$,  resulting in $4L/\delta x = 4L/\delta y = m = 20$ grid points along each axis. The time step size is chosen as $\delta t = 3  \delta x \delta y$, and the biharmonic coefficient is specified as $\kappa = 0.001 \delta x \delta y$. Notably, the equations under consideration have been nondimensionalized. With an iteration count of $N=300$, we achieve a total simulation time $T=N\delta t = 36$. The initial condition is randomly generated as $\omega_{0;ij} \sim 5\mathcal{N}(0,1)$ for $i,j=1,2,\ldots, m-1$.
\begin{figure}[htb!]
\centering
\begin{subfigure}[t]{0.45\linewidth}
\centering
\includegraphics[scale=0.18]{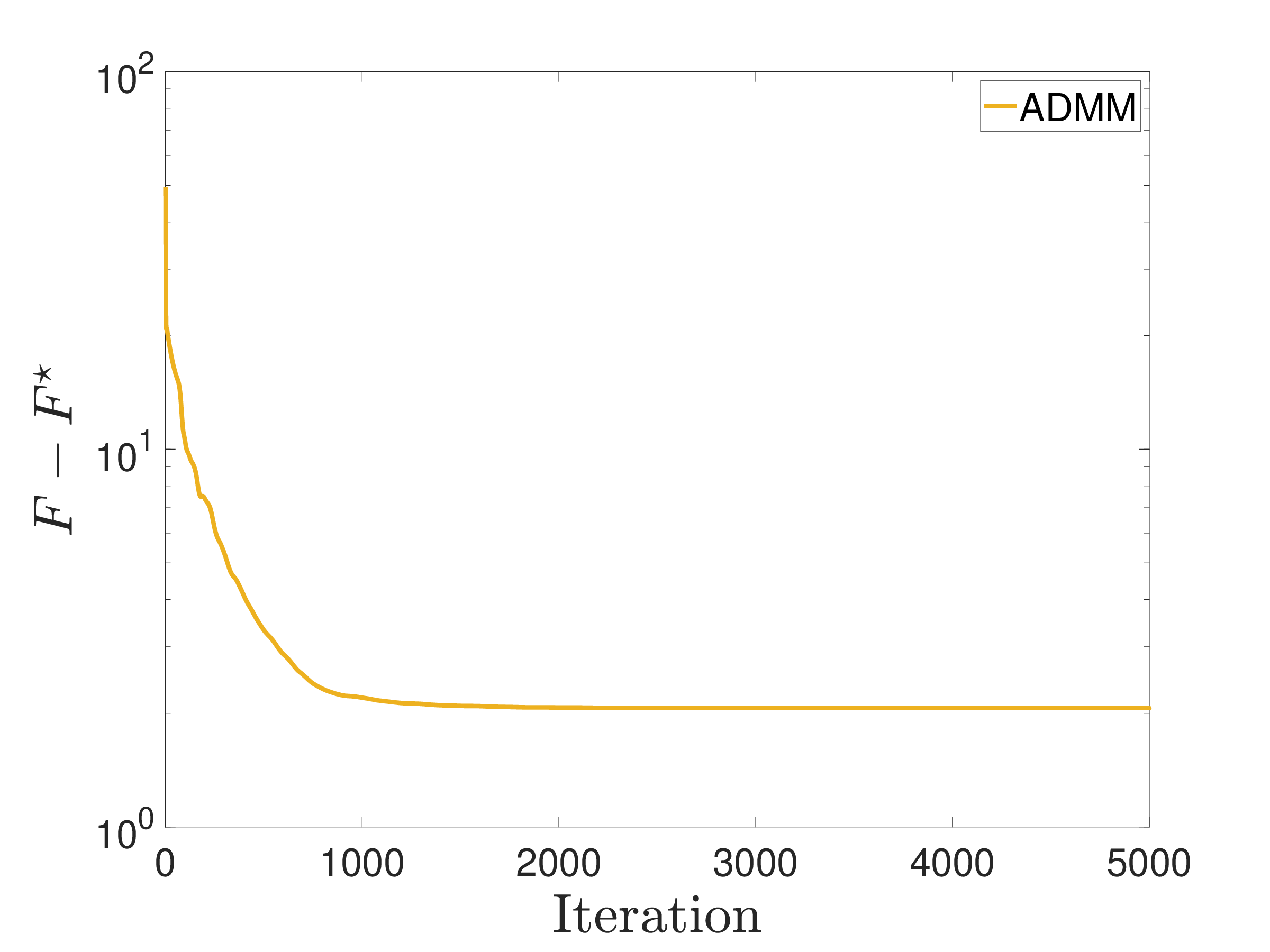}
\caption{\small Total error}
\end{subfigure}
\begin{subfigure}[t]{0.45\linewidth}
\centering
\includegraphics[scale=0.18]{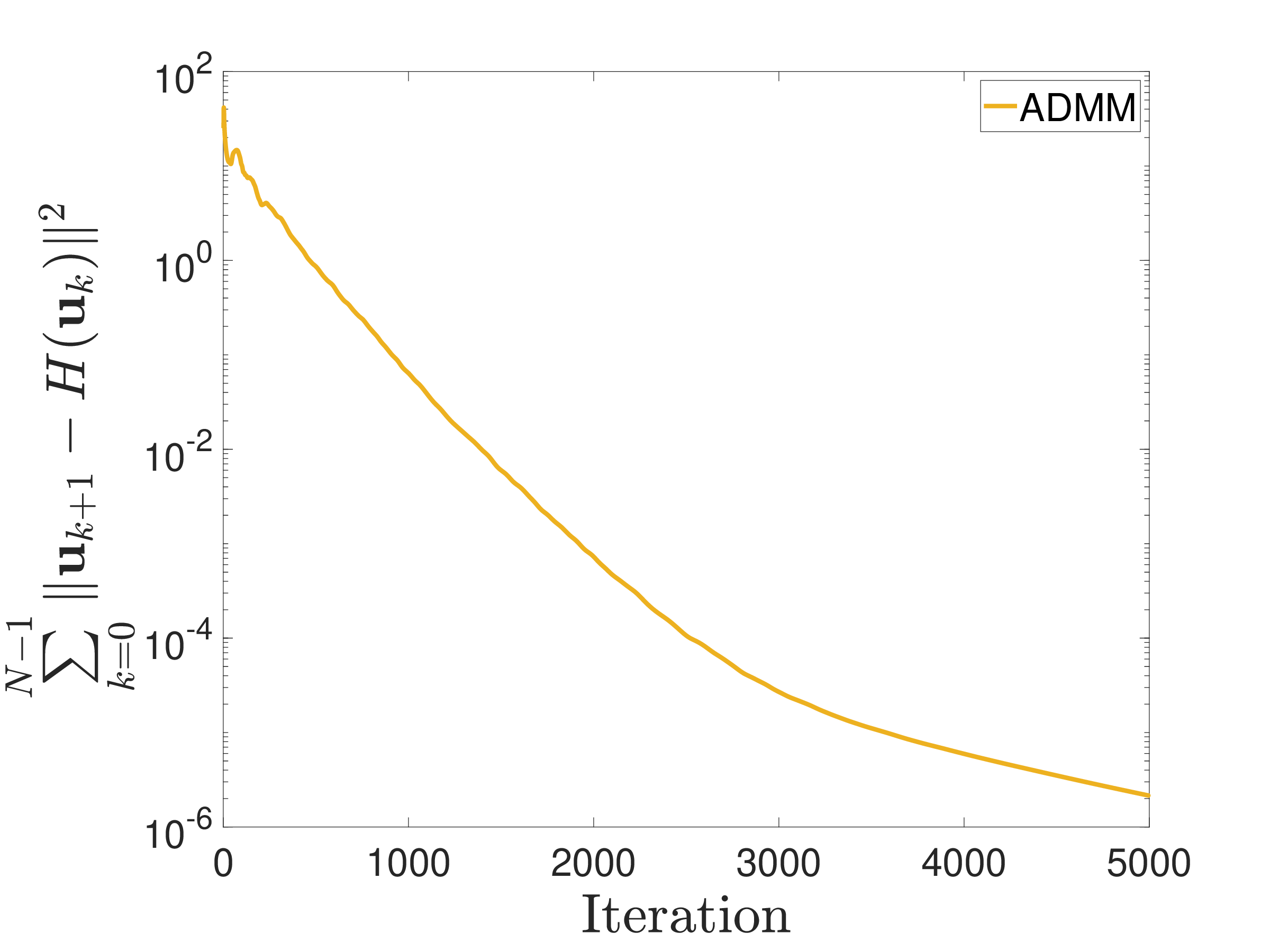}
\caption{\small Constraint error}
\end{subfigure}
\caption{\small  The linearized multi-block~\texttt{ADMM} with regularization~\eqref{eqn: lin-ADMM-reg} is applied to  the~\texttt{4D-Var} problem~\eqref{eqn: numerical-4dvar-admm} for the Burgers' equation with the 2D large-scale vorticity equation~\eqref{eqn: 2d-qg-f}, using the finite difference scheme described in~\Cref{subsec: computational-implementation}. } 
\label{fig: vorticity-convergence}
\end{figure}
The numerical solution, considered as the precise observation, is obtained following the procedures outlined in~\Cref{subsec: computational-implementation},  To generate noisy observational data, we add Gaussian noise to the observations: $\hat{\omega}_{k,ij} = \omega_{30k,ij} +  \cdot 0.5\varepsilon _{ij}$, where $\varepsilon_{ij} \sim \mathcal{N}(0,1)$ for $i.j=1,2,\ldots,m-1$. The linearized multi-block~\texttt{ADMM} with regularization~\eqref{eqn: lin-ADMM-reg} begins with the initial condition as $\pmb{\omega}_k^{0} = {0}_{(m-1)\times (m-1)}^{\top}$ for $k=0,1,\ldots,N$, using the parameters $\mu = 20$, $\eta = 0.1$, and $s=2/3$. The numerical performance is illustrated in~\Cref{fig: vorticity-convergence}, which shows consistency with the convergence behaviors observed for the Lorenz system in~\Cref{sec: outline-admm} and various numerical schemes of the viscous Burgers' equation in~\Cref{sec: burgers}.  
\begin{figure}[htb!]
\centering
\begin{subfigure}[t]{0.96\linewidth}
\centering
\includegraphics[scale=0.26]{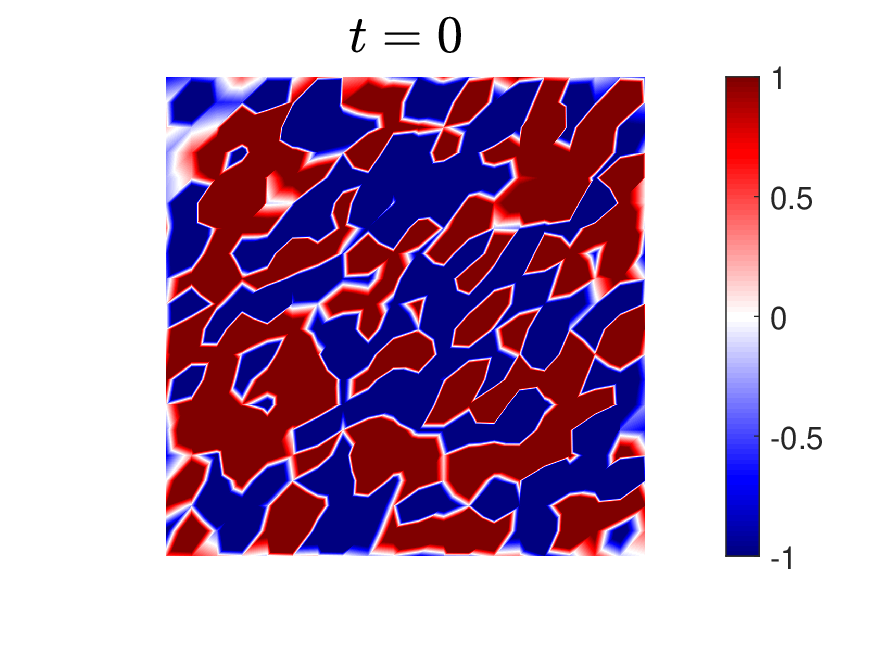}
\includegraphics[scale=0.26]{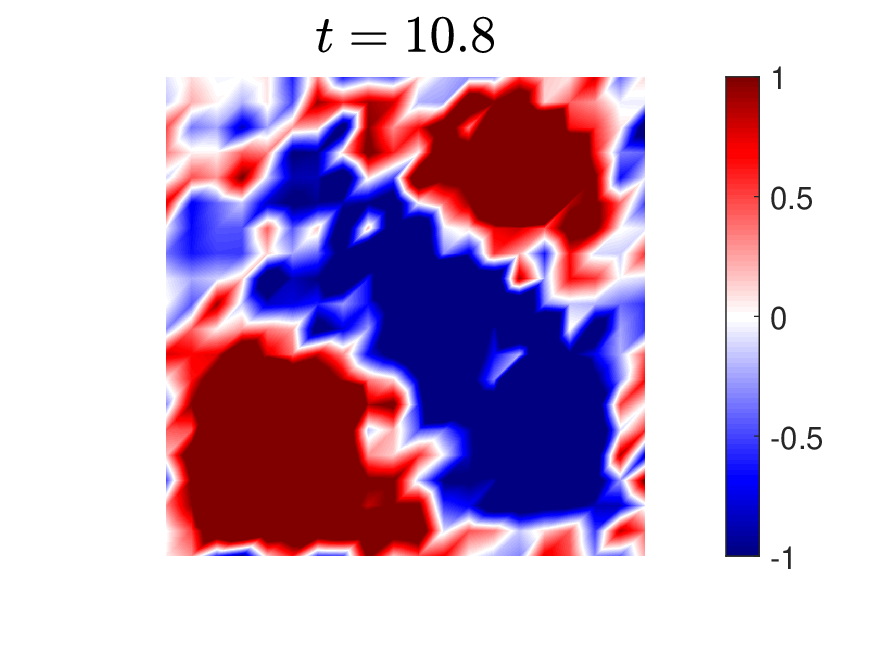}
\includegraphics[scale=0.26]{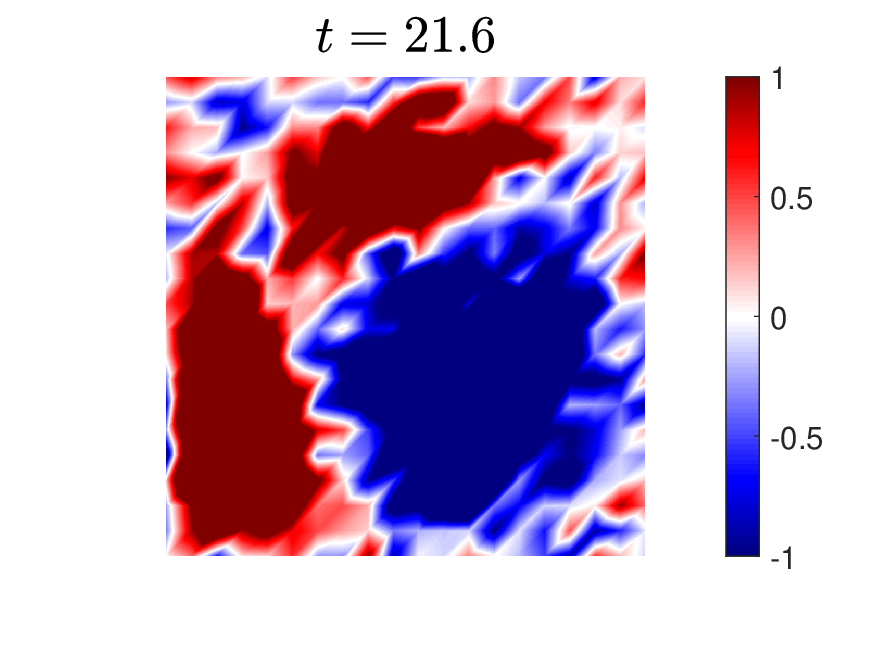}
\includegraphics[scale=0.26]{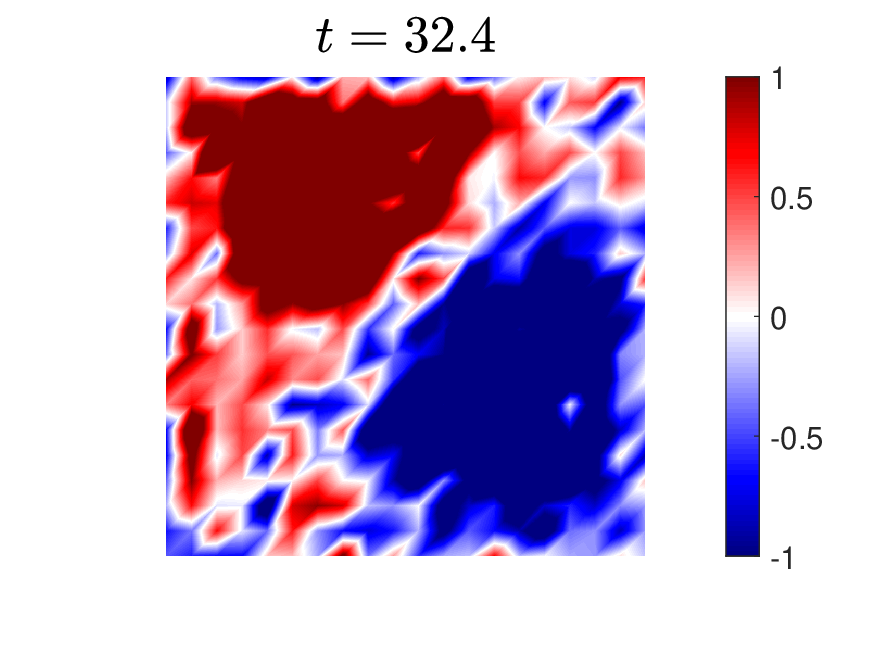}
\subcaption{\small Nosiy observational data}
\end{subfigure}
%\begin{subfigure}[t]{0.24\linewidth}
%\centering
%\includegraphics[scale=0.30]{vorticity/obs2}
%\end{subfigure}
%\begin{subfigure}[t]{0.24\linewidth}
%\centering
%\includegraphics[scale=0.30]{vorticity/obs3}
%\end{subfigure}
%\begin{subfigure}[t]{0.24\linewidth}
%\centering
%\includegraphics[scale=0.30]{vorticity/obs4}
%\end{subfigure}
\begin{subfigure}[t]{0.96\linewidth}
\centering
\includegraphics[scale=0.26]{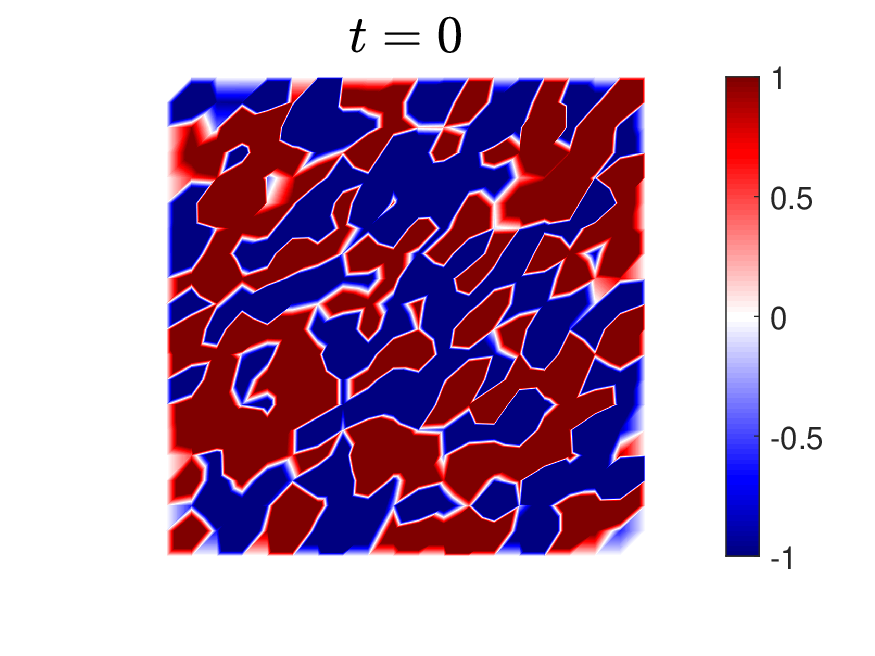}
\includegraphics[scale=0.26]{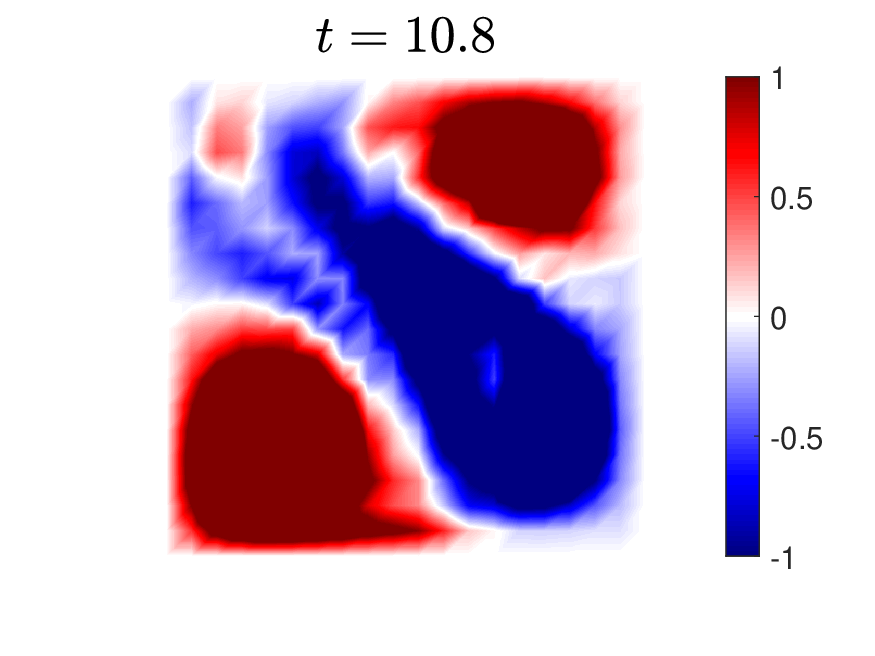}
\includegraphics[scale=0.26]{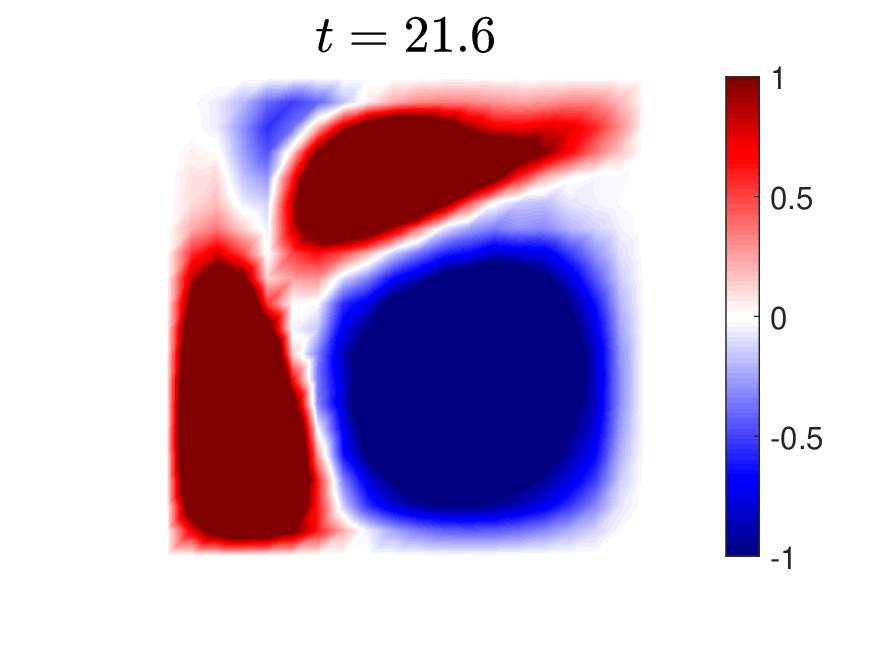}
\includegraphics[scale=0.26]{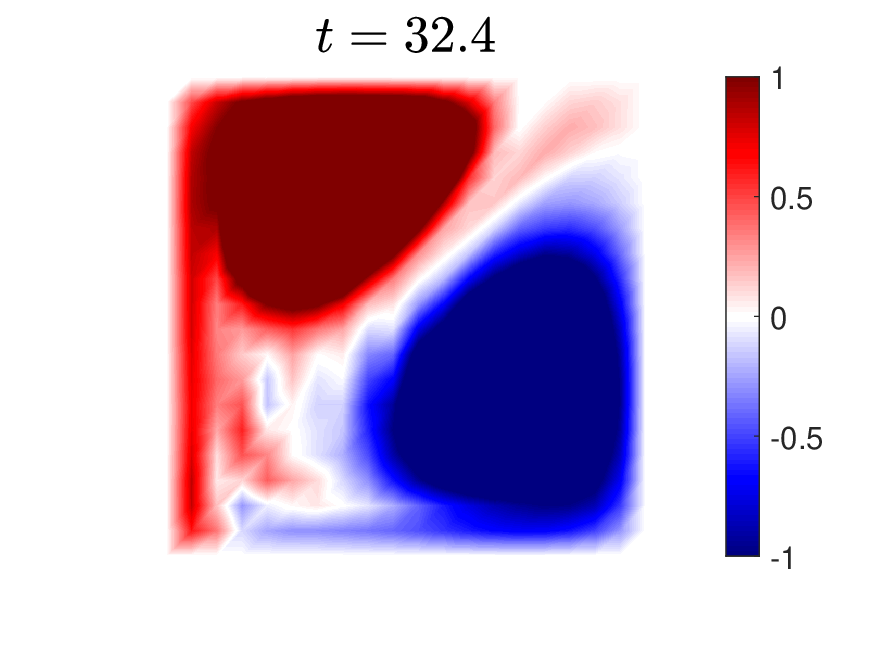}
\subcaption{\small Ture numerical solution}
\end{subfigure}
%\begin{subfigure}[t]{0.24\linewidth}
%\centering
%\includegraphics[scale=0.20]{vorticity/ture2}
%\end{subfigure}
%\begin{subfigure}[t]{0.24\linewidth}
%\centering
%\includegraphics[scale=0.20]{vorticity/ture3}
%\end{subfigure}
%\begin{subfigure}[t]{0.24\linewidth}
%\centering
%\includegraphics[scale=0.20]{vorticity/ture4}
%\end{subfigure}
\begin{subfigure}[t]{0.96\linewidth}
\centering
\includegraphics[scale=0.26]{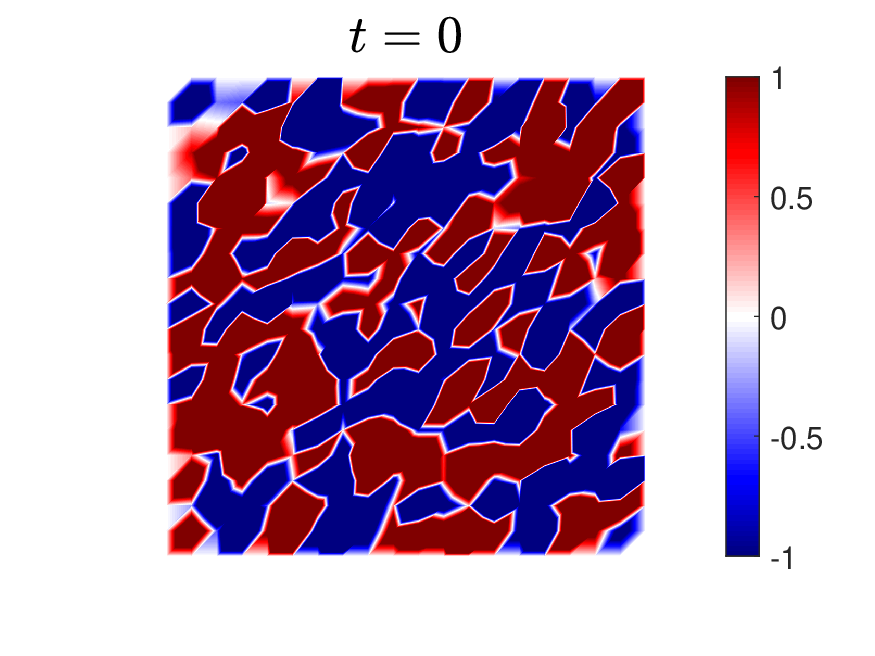}
\includegraphics[scale=0.26]{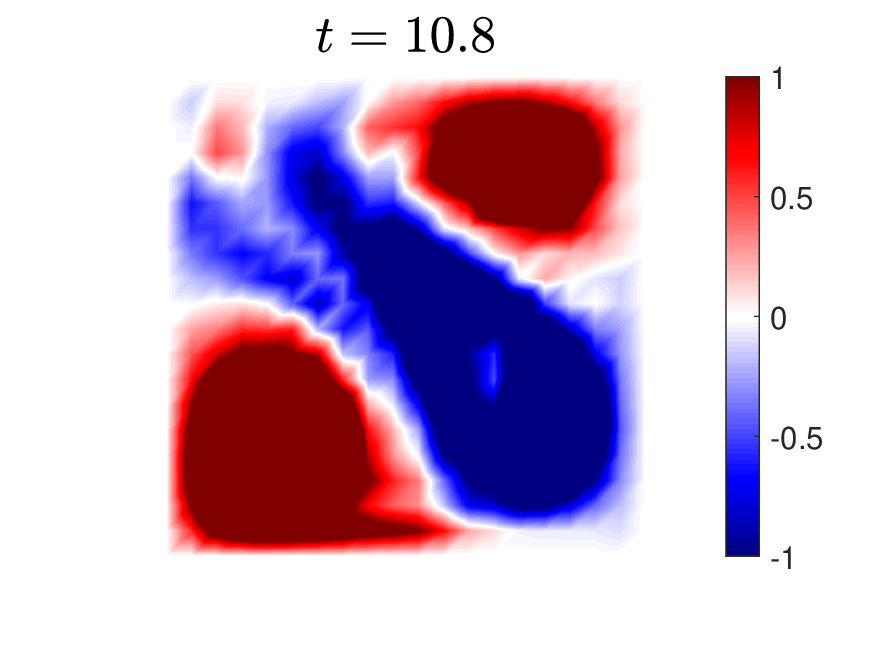}
\includegraphics[scale=0.26]{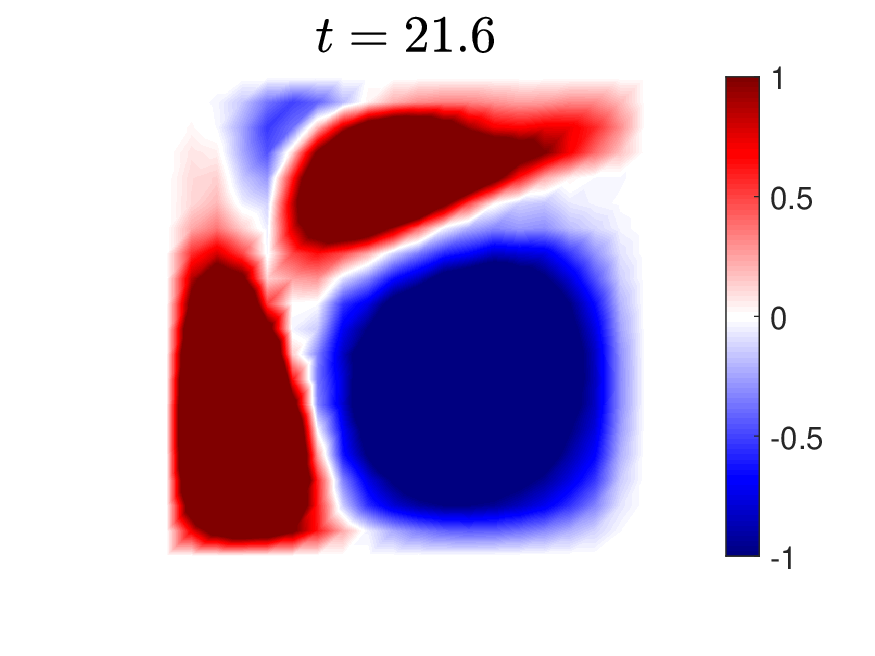}
\includegraphics[scale=0.26]{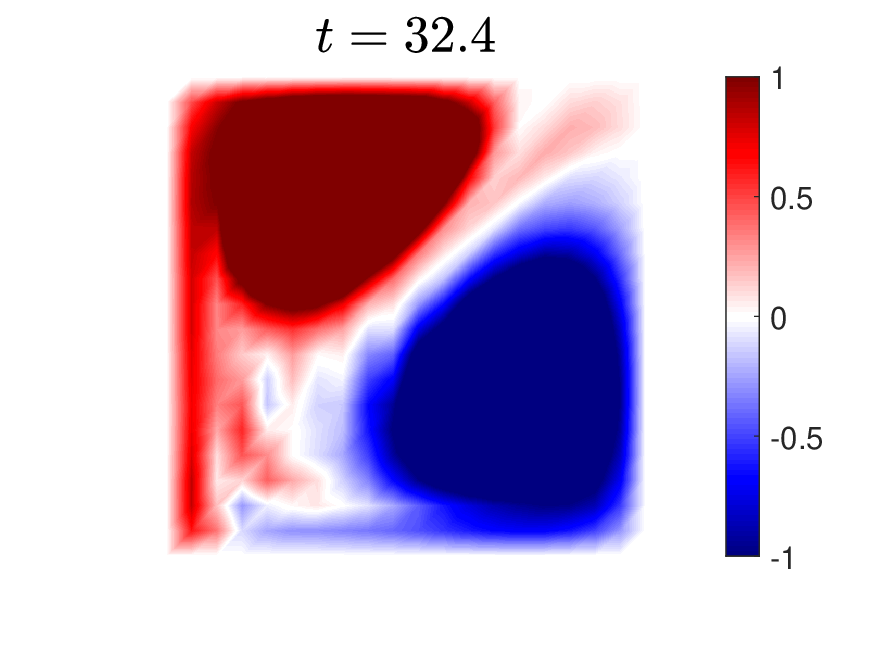}
\subcaption{\small Dynamics recovered via~\texttt{ADMM}}
\end{subfigure}
%\begin{subfigure}[t]{0.24\linewidth}
%\centering
%\includegraphics[scale=0.20]{vorticity/admm2}
%\end{subfigure}
%\begin{subfigure}[t]{0.24\linewidth}
%\centering
%\includegraphics[scale=0.20]{vorticity/admm3}
%\end{subfigure}
%\begin{subfigure}[t]{0.24\linewidth}
%\centering
%\includegraphics[scale=0.20]{vorticity/admm4}
%\end{subfigure}
\caption{\small   Comparison of noisy observational data with the time evolution dynamics recovered via the linearized multi-block~\texttt{ADMM} with regularization~\eqref{eqn: lin-ADMM-reg}, with the true numerical solution as the reference.} 
\label{fig: vorticity}
\end{figure}
Furthermore,~\Cref{fig: vorticity} highlights the recovered solution closely matches the true dynamics, even in the presence of noisy observational data.  The numerical errors, presented in~\Cref{fig: vorticity-error}, further verify that the accuracy improves over time as the simulation progresses.

By considering the scales of large-scale ocean circulation, with $L\sim2.5\mathrm{km}$ and  vorticity $\omega_0 \sim 10^{-5}\mathrm{s}^{-1}$, we derive  a characteristic time scale of $\delta t \sim 7.2\mathrm{h}$. With an iteration count of $N=300$,  the total simulation time becomes $T=N\delta t \approx 6\;\mathrm{years}$,  allowing for a comprehensive characterization of the large-scale ocean circulation. Additionally, we note that a commonly used boundary condition is the periodic boundary condition~\citep{mcwilliams1984emergence}, which necessitates the use of spectral methods rather than finite difference methods. The reason for this preference is that convergence using  \texttt{SOR} and the conjugate gradient methods  is typically poor due to the presence of a zero eigenvalue in the Laplacian operator.

%% file: 05_conclu.tex
\section{Conclusions and future work}
\label{sec: conclu-future}

In this study, we propose a linearized multi-block \texttt{ADMM} with regularization to solve the \texttt{4D-Var} problem, exploiting its separable structure. Unlike classical first-order optimization algorithms that primarily focus on initial conditions, our approach derives the Euler-Lagrange equation for the entire dynamical system, facilitating more effective use of observational data. When the initial condition is poorly chosen, the $\argmin$ operation steers the iteration towards the observational data, reducing sensitivity to the initial guess. The quadratic subproblems simplify the solution process, while the parallel structure enhances efficiency, particularly with modern computational resources. However, the computation of the adjoint operator remains a significant challenge, even in cases like 2D turbulence. Future research may explore addressing this issue using the sampling approach proposed by~\citep{shi2023adjoint, shi2024sampling}. Additionally, stochastic gradient descent (SGD) and its variants, which have proven successful in deep neural networks~\citep{shi2023learning,shi2021hyperparameters}, offer promising alternatives for solving~\texttt{4D-Var} problems. Another important direction for further research is the investigation of complex optimization problems~\citep{sorber2012unconstrained}, as many governing equations, such as the nonlinear Schr\"odinger and Ginzburg-Landau equations, involve complex variables. Further theoretical and numerical studies, including the exploration of various numerical schemes, may provide valuable insights, especially in fields such as meteorology, oceanography, and climatology.

One of the key unresolved challenges is proving the convergence of the linearized multi-block ADMM with regularization and identifying the conditions under which this convergence occurs. Establishing convergence is essential, as it would confirm that the solution to the~\texttt{4D-Var} problem~\eqref{eqn: numerical-4dvar} is independent of numerical schemes, enabling a numerical realization for the implicit scheme. Additionally, for the analytic \texttt{4D-Var} problem~\eqref{eqn: analytic-4dvar}, the issue of uniqueness, particularly over longer time horizons, remains open, although it has been addressed for short intervals~\citep{cox2015non}. A promising future direction lies in developing infinite-dimensional \texttt{ADMM}-like iterations. In the case where $T=0$, such an approach could lead to a novel method for constructing a solution to the governing differential equation. This method could draw inspiration from classical iterative techniques, such as Picard iteration and Newton iteration~\citep{arnold1992ordinary}, which have been successfully extended from finite to infinite dimensions. These ideas could also connect with advanced frameworks like KAM iteration~\citep{kolmogorov1954conservation, arnold2009proof, moser1962invariant} and Nash-Moser iteration~\citep{nash1956imbedding, moser1966rapidly1, moser1966rapidly2}, which have significant applications, including in Landau damping~\citep{mouhot2011landau}. Further exploration of these concepts could deepen our understanding of \texttt{ADMM} in infinite-dimensional spaces.